\documentclass[review,12pt]{elsarticle}

\usepackage[a4paper,width=170mm,top=30mm,bottom=25mm]{geometry}

\usepackage{lineno,hyperref}
\usepackage{multirow}
\usepackage{amsmath}

\numberwithin{figure}{section}
\numberwithin{equation}{section}
\numberwithin{table}{section}
\usepackage{color}
\usepackage{booktabs}
\usepackage{graphicx}
\usepackage{subfigure}
\usepackage{setspace}
\usepackage{amssymb}
\usepackage{bm}
\counterwithout{equation}{section}
\counterwithout{figure}{section}
\counterwithout{table}{section}
\usepackage{xcolor}
\usepackage{xfrac}

\usepackage{caption}

\journal{arXiv preprint}

\begin{document}

\begin{frontmatter}

\title{Three-dimensional DtN-FEM scattering analysis of Lamb and SH guided waves by a symmetric cavity defect in an isotropic infinite plate}


\author[mymainaddress,mysecondaryaddress]{Chen Yang\corref{mycorrespondingauthor}}
\ead{yangc@pku.edu.cn; yc@nuaa.edu.cn; yang.c.ag@m.titech.ac.jp}

\author[mysecondaryaddress]{Junichi Nakaoka}
\ead{junnakaoka0224@gmail.com}

\author[mysecondaryaddress]{Sohichi Hirose}
\ead{shirose@cv.titech.ac.jp}

\cortext[mycorrespondingauthor]{Corresponding author}
\address[mymainaddress]{State Key Lab for Turbulence and Complex Systems, College of Engineering, Peking University, Beijing 100871, China}
\address[mysecondaryaddress]{Tokyo Institute of Technology, 2-12-1-W8-22, O-okayama, Meguro, Tokyo, 152-8552, Japan}

\begin{abstract}

In this paper, a three-dimensional Dirichlet-to-Neumann (DtN) finite element method (FEM) is developed to analyze the scattering of Lamb and SH guided waves due to a symmetric cavity defect in an isotropic infinite plate. During the finite element analysis, it is necessary to determine the far-field DtN conditions at virtual boundaries where both displacements and tractions are unknown. In this study, firstly, the scattered waves at the virtual boundaries are represented by a superposition of guided waves with unknown scattered coefficients. Secondly, utilizing the mode orthogonality, the unknown tractions at virtual boundaries are expressed in terms of the unknown scattered displacements at virtual boundaries via scattered coefficients. Thirdly, this relationship at virtual boundaries can be finally assembled into the global DtN-FEM matrix to solve the problem. This method is simple and elegant, which has advantages on dimension reduction and needs no absorption medium or perfectly matched layer to suppress the reflected waves compared to traditional FEM. Furthermore, the reflection and transmission coefficients of each guided mode can be directly obtained without post-processing. This proposed DtN-FEM will be compared with boundary element method (BEM), and finally validated for several benchmark problems.

\end{abstract}


\end{frontmatter}

\section{Introduction}

Nondestructive testings demonstrate its absolute advantages and potential compared with other traditional destructive methods over the past few decades. As one of the nondestructive testings, the ultrasonic guided waves have naturally captured the attention of researchers. Although the conventional bulk wave detection technique has been widely used in many fields, the guided wave has lower preconditions, is more convenient to operate and needs to process fewer waveforms. These unique superiors cannot be replicated by bulk wave detection. Undoubtedly, the emergence of the guided wave technique has considerably reduced the cost of nondestructive testing.

Currently, many researchers have obtained huge achievements in quantitative non-destructive testing (QNDT) using ultrasonic guided waves \cite{Da20201479,Da20201923}. Wang proposed a prevalent method for shape and depth construction of plate, half space and layered semi-infinite space using whole reflection coefficients of SH waves \cite{Wang20121782,Wang2015}, Rayleigh Wave \cite{Wang2018} and Love waves \cite{Wang2019}, respectively. In addition, Da et al. \cite{Da2018181} utilize scattering coefficients of the torsional modes to reconstruct the defect shape in pipeline structures effectively. Consequently, it is enormously important to require a thorough understanding and accurate computation of forward scattering phenomenon, to acquire a database for QNDT research such as near-and-far field scattering data.

In the analysis of elastic guided wave scattering problems, two key aspects must be properly addressed and treated. First, due to the dispersive and multi-mode nature of guided and evanescent waves, the wavenumber and wave structures should be accurately determined by dispersion equations. Moreover, dispersion curves can be drawn analytically by solving transcendent equations with complex root-racking modulus-converging algorithm \cite{Yang2018}. This information is useful for mode selection, as well as generation and reception of a single mode in non-destructive evaluation. Second, Lamb wave mode interaction with defects must be accurately determined, which is the focus of this paper. The complexity of mode conversion restrains an analytical approach to very simple geometries, and for more general cases, it requires valid numerical methods like FEM and BEM to obtain an accurate scattering wave-field. In our previous work, a modified BEM was proposed to solve the three-dimensional plate scattering problem \cite{Yang2021145}, where only the interfaces and flaw boundaries need to be discretized. However, BEM needs to handle the singularity and compute the fundamental solution matrix numerically, and the final BEM matrix will be full-rank, which consume many storage and computing resources, especially when simulating high frequency scattering problems. 

Compared to BEM, FEM is particularly attractive for forward scattering analysis since FEM is more simple and popular, and the final global FEM matrix is sparse. For plate scattering computation, Koshiba \cite{Koshiba198418} implemented the analysis of the scattering problem of Lamb waves in an elastic waveguide by FEM. Al-Nassar and Datta \cite{Al-Nassar1991125} used the Lamb wave to detect a normal rectangular strip weldment and obtained the guided waves’ behavior at specific frequencies. A semi-analytical FEM proposed by Hayashi and Rose \cite{Hayashi200375} was used to analyze the behavior of guided waves in plates and pipes. This method was soon applied to the field of flaw detection, for example, rail base \cite{Ramatlo2020,Yue2020}. Gunawan and Hirose \cite{Gunawan2004996} provide a mode-exciting method to study the scattering phenomenon in an infinite plate. Moreover, the researchers also conducted related explorations in the processing approach of the scattering boundary. The orthogonality of guided wave modes was provided to solve the scattering problem in a plate \cite{Shkerdin20042089}. Other researchers, such as Moreau \cite{Moreau2006611}, and Kubrusly \cite{Kubrusly2019,Kubrusly2021} try to apply this method to study relative objects.

The goal of this paper is to propose a FEM implemented in a bi-material plate to avoid the computation error caused by spurious reflected waves in 3-D scattering problems. During the finite element analysis, it is necessary to determine the far-field DtN conditions at virtual boundaries where both displacements and tractions are unknown. In this study, firstly, the scattered waves at the virtual boundaries are represented by a superposition of guided waves with unknown scattered coefficients. Secondly, utilizing the mode orthogonality, the unknown tractions at virtual boundaries are expressed in terms of the unknown scattered displacements at virtual boundaries via scattered coefficients. Thirdly, this relationship at virtual boundaries can be finally assembled into the global DtN-FEM matrix to solve the problem. This method is simple and elegant, which has advantages on dimension reduction and needs no absorption medium or perfectly matched layer to suppress the reflected waves compared to traditional FEM. Furthermore, the reflection and transmission coefficients of each guided mode can be directly obtained without post-processing. This proposed DtN-FEM will be compared with boundary element method (BEM), and finally validated for several benchmark problems. Finally, this proposed DtN-FEM will be applied into investigating the effect on different defect shapes.

\section{Three-dimensional DtN FEM}
\subsection{General description of the problem}
\begin{figure}[h]
	\centering
	\includegraphics[scale=0.9]{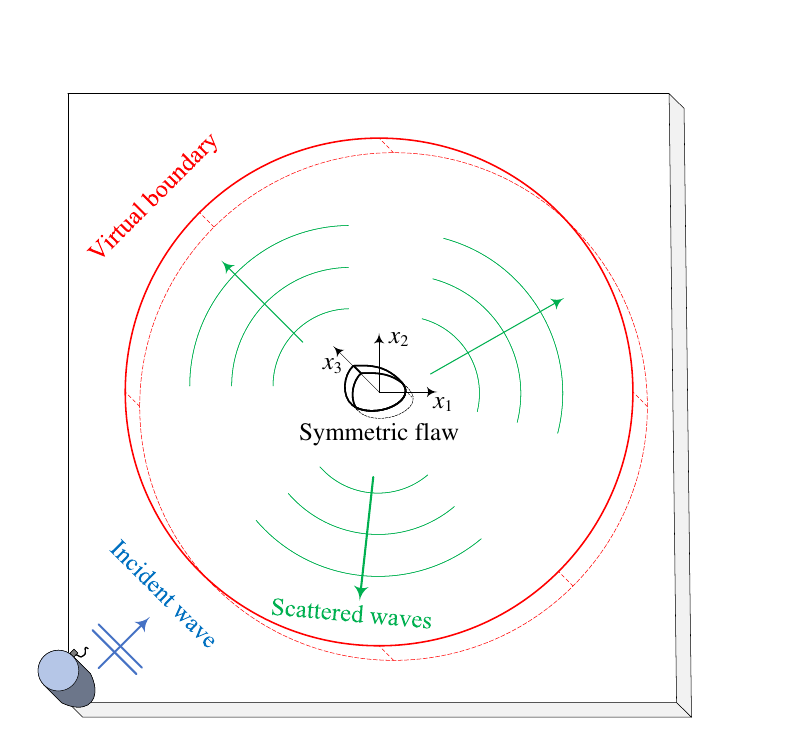}
	\caption{Boundary conditions along the horizontal direction.}
	\label{fig_problem}
 \end{figure}

A three-dimensional elasto-dynamic problem in a infinite plate of homogeneous, isotropic, linearly elastic material with thickness $2h$ is considered as shown in Fig.\ref{fig_problem},  referred to a Cartesian coordinate system such that $x_1-x_2$ plane coincides with the mid-plane of the layer. The cavity flaw is symmetric about the mid-plane of the layer. The upper and lower surfaces of the plate and the flaw surface are stress-free. Let us consider an incident plane wave excited at the far-field interacting with a cavity of arbitrary shape and generating scattered waves. By virtue of linear superposition principle, the total field in the flawed structure can be considered as the superposition of the incident and scattered fields where the incident field is known and the scattered field is unknown. The aim of this paper is to simulate the unknown scattered field. In traditional FEM, extra absorbing materials or perfectly matched layered is utilized to absorb the spurious reflected waves in order to simulate the infinite domain. However, in this paper, the domain surrounded by the red dashed surface in Fig.\ref{fig_problem} is the actual simulation area. The key to the problem is how to deal with the far-field DtN conditions at the virtual boundary since displacements and tractions at the virtual boundary are both unknown.

\subsection{Three-dimensional guided wave assumptions in an infinite isotropic plate}

In the modeling of the problem using three-dimensional scattering theory, it is crucial to express the wave fields in an appropriate way. To analyze the scattered wave field in an infinite plate, it is a natural choice to adopt a cylindrical coordinate system $(r,\theta,z)$ to express the Lamb and SH wave modes. A way to represent the following Lamb and SH wave modes is to use a potential function as shown by Achenbach, and the general equation is as follows.

For a Lamb wave mode of order $n$, the displacements can be given as
\begin{equation}
u_r^n = \frac{1}{{{k_n}}}{V_n}\left( z \right)\frac{{\partial \phi }}{{\partial r}}\left( {r,\theta } \right),\quad u_\theta ^n = \frac{1}{{{k_n}r}}{V_n}\left( z \right)\frac{{\partial \phi }}{{\partial \theta }}\left( {r,\theta } \right),\quad u_z^n = {W_n}\left( z \right)\phi \left( {r,\theta } \right).
\label{eq2.1}
\end{equation}

\noindent
where the scalar potential function $\phi \left( {r,\theta } \right)$ satisfies Helmholtz equation which in cylindrical coordinates reads
\begin{equation}
\frac{{{\partial ^2}\phi }}{{\partial {r^2}}} + \frac{1}{r}\frac{{\partial \phi }}{{\partial r}} + \frac{1}{{{r^2}}}\frac{{{\partial ^2}\phi }}{{\partial {\theta ^2}}} + k_n^2\phi  = 0
\label{eq2.2}
\end{equation}

In Eq.\ref{eq2.1}, $k_n$ represents the wavenumber and ${V_n}\left( z \right)$ and ${W_n}\left( z \right)$ are the thickness coordinate dependent functions which differ for symmetric and anti-symmetric modes. For symmetric Lamb wave, expressions for ${V_n}\left( z \right)$ and ${W_n}\left( z \right)$ are
\begin{equation}
V_{n}(z) = s_{1}{\cos\left( {p_{n}z} \right)} + s_{2}{\cos\left( {q_{n}z} \right)},\quad 
W_{n}(z) = s_{3}{\sin\left( {p_{n}z} \right)} + s_{4}{\sin\left( {q_{n}z} \right)},
\label{eq2.3}
\end{equation}

\noindent
where 
\begin{equation}
s_{1} = 2{\cos\left( {q_{n}h} \right)},\quad s_{2} = - \frac{k_{n}^{2} - q_{n}^{2}}{k_{n}^{2}}{\cos\left( {p_{n}h} \right)},
\label{eq2.4}
\end{equation}

\begin{equation}
s_{3} = - \frac{2p_{n}}{k_{n}}{\cos\left( {q_{n}h} \right)},\quad 
s_{4} = - \frac{k_{n}^{2} - q_{n}^{2}}{q_{n}k_{n}}{\cos\left( {p_{n}h} \right)}.
\label{eq2.5}
\end{equation}

$k_n$ is given by the root of the symmetric Rayleigh-Lamb equation
\begin{equation}
\frac{{\tan \left( {qh} \right)}}{{\tan \left( {ph} \right)}} =  - \frac{{4{k^2}pq}}{{{{\left( {{q^2} - {k^2}} \right)}^2}}}
\label{eq2.6}
\end{equation}

\noindent
where ${p^2} = k_L^2 - k_n^2, {q^2} = k_T^2 - k_n^2$, and $\omega$ represents the angular frequency and $k_L$ and $k_T$ indicate the longitudinal and transverse wave velocities respectively.

SH wave modes can also be generated. Similarly, these modes can be expressed using a scalar potential function $\psi \left( {r,\theta } \right)$, as shown by
\begin{equation}
u_r^n = \frac{1}{{{l_n}r}}{U_n}\left( z \right)\frac{{\partial \psi }}{{\partial \theta }}\left( {r,\theta } \right),\quad u_\theta ^n =  - \frac{1}{{{l_n}}}{U_n}\left( z \right)\frac{{\partial \psi }}{{\partial r}}\left( {r,\theta } \right),\quad u_z^n = 0.
\label{eq2.7}
\end{equation}

\noindent
where the scalar potential $\psi \left( {r,\theta } \right)$ satisfies Helmholtz equation
\begin{equation}
\frac{{{\partial ^2}\psi }}{{\partial {r^2}}} + \frac{1}{r}\frac{{\partial \psi }}{{\partial r}} + \frac{1}{{{r^2}}}\frac{{{\partial ^2}\psi }}{{\partial {\theta ^2}}} + l_n^2\psi  = 0
\label{eq2.8}
\end{equation}

\noindent
and the thickness dependent functions ${U_n}\left( z \right)$ for symmetric SH modes are given by
\begin{equation}
{U_n}\left( z \right){\rm{ = }}\cos \left( {\frac{{n\pi z}}{{2h}}} \right)\quad \left( {n = 0,2,4...} \right)
\label{eq2.9}
\end{equation}

The wavenumbers $l_n$ for the SH modes are
\begin{equation}
l_n^2 = k_T^2 - {\left( {\frac{{n\pi }}{{2h}}} \right)^2}
\label{eq2.10}
\end{equation}

 Next, we need to find proper expansions for the wave field. For this purpose, the solution of the scalar potentials $\phi \left( {r,\theta } \right)$ and $\psi \left( {r,\theta } \right)$ referred to reference for an outgoing wave are expressed as
\begin{equation}
\phi \left( {r,\theta } \right){\rm{ = }}{H_m}\left( {{k_n}r} \right){{\rm{e}}^{{\rm{i}}m\theta }},\quad \psi \left( {r,\theta } \right){\rm{ = }}{H_m}\left( {{l_n}r} \right){{\rm{e}}^{{\rm{i}}m\theta }}.
\label{eq2.11}
\end{equation}

\noindent
respectively, where ${H_m}\left( {} \right)$ indicate Hankel functions of the first kind.

Then the scattered displacement fields are obtained by expanding the wave fields in the allowed modes. For a fixed frequency, Eq.\ref{eq2.6} and Eq.\ref{eq2.10} have a finite number of real roots, corresponding to propagating modes, and an infinite number of complex roots, non-propagating modes. Thus, the scattered displacements and stresses at the far-field virtual boundary can be written as 
\begin{equation}
\begin{aligned}
u_{r}^{scat} &= {\sum\limits_{n = 0}^{\infty}{\sum\limits_{m = - \infty}^{\infty}{A_{mn}V_{n}(z){H^{'}}_{m}\left( {k_{n}r} \right)\text{e}^{\text{i}m\theta}}}} + {\sum\limits_{n = 0,2,\cdots}^{\infty}{\sum\limits_{m = - \infty}^{\infty}{\text{i}mB_{mn}{\cos\left( \frac{n\pi z}{2h} \right)}\frac{H_{m}\left( {l_{n}r} \right)}{l_{n}r}\text{e}^{\text{i}m\theta}}}}, \\
u_{\theta}^{scat} &= {\sum\limits_{n = 0}^{\infty}{\sum\limits_{m = - \infty}^{\infty}{\text{i}mA_{mn}V_{n}(z)\frac{H_{m}\left( {k_{n}r} \right)}{k_{n}r}\text{e}^{\text{i}m\theta}}}} - {\sum\limits_{n = 0,2,\cdots}^{\infty}{\sum\limits_{m = - \infty}^{\infty}{B_{mn}{\cos\left( \frac{n\pi z}{2h} \right)}{H^{'}}_{m}\left( {l_{n}r} \right)\text{e}^{\text{i}m\theta}}}}, \\
u_{z}^{scat} &= {\sum\limits_{n = 0}^{\infty}{\sum\limits_{m = - \infty}^{\infty}{A_{mn}W_{n}(z)H_{m}\left( {k_{n}r} \right)\text{e}^{\text{i}m\theta}}}}.
\end{aligned}
\label{eq2.12}
\end{equation}
\begin{equation}
\begin{aligned}
\sigma_{r}^{scat} &= {\sum\limits_{n = 0}^{\infty}{\sum\limits_{m = - \infty}^{\infty}{A_{mn}\left\lbrack {\textstyle \sum_{rr}^{n}(z)H_{m}\left( {k_{n}r} \right) - {\overset{\sim}{\textstyle \sum}}_{rr}^{n}(z)\left( {\frac{1}{r}{H^{'}}_{m}\left( {k_{n}r} \right) - \frac{m^{2}}{k_{n}r^{2}}H_{m}\left( {k_{n}r} \right)} \right)} \right\rbrack\text{e}^{\text{i}m\theta}}}} \\
&+ {\sum\limits_{n = 0,2,\cdots}^{\infty}{\sum\limits_{m = - \infty}^{\infty}{\text{i}mB_{mn}\mu{\cos\left( \frac{n\pi z}{2h} \right)}\left( {\frac{2}{r}{H'}_{m}\left( {l_{n}r} \right) - \frac{2}{l_{n}r^{2}}H_{m}\left( {l_{n}r} \right)} \right)\text{e}^{\text{i}m\theta}}}}, \\
\tau_{r\theta}^{scat} &= {\sum\limits_{n = 0}^{\infty}{\sum\limits_{m = - \infty}^{\infty}{\text{i}mA_{mn}\textstyle \sum_{r\theta}^{n}(z)\left( {\frac{1}{r}{H^{'}}_{m}\left( {k_{n}r} \right) - \frac{1}{k_{n}r^{2}}H_{m}\left( {k_{n}r} \right)} \right)\text{e}^{\text{i}m\theta}}}} \\
&+ {\sum\limits_{n = 0,2,\cdots}^{\infty}{\sum\limits_{m = - \infty}^{\infty}{B_{mn}\mu{\cos\left( \frac{n\pi z}{2h} \right)}\left\lbrack {\frac{2}{r}{H'}_{m}\left( {l_{n}r} \right) + \left( {l_{n} - \frac{2m^{2}}{l_{n}r^{2}}} \right)H_{m}\left( {l_{n}r} \right)} \right\rbrack\text{e}^{\text{i}m\theta}}}} \\
\tau_{rz}^{scat} &= - {\sum\limits_{n = 0}^{\infty}{\sum\limits_{m = - \infty}^{\infty}{A_{mn}\textstyle \sum_{rz}^{n}(z){H^{'}}_{m}\left( {k_{n}r} \right)\text{e}^{\text{i}m\theta}}}} \\
&- {\sum\limits_{n = 2,4,\cdots}^{\infty}{\sum\limits_{m = - \infty}^{\infty}{\text{i}mB_{mn}\mu\frac{n\pi}{2h}{\sin\left( \frac{n\pi z}{2h} \right)}\frac{1}{l_{n}r}H_{m}\left( {l_{n}r} \right)\text{e}^{\text{i}m\theta}}}}.
\end{aligned}
\label{eq2.13}
\end{equation}

\noindent
respectively, where $A_{mn}$ and $B_{mn}$ account for the expansion coefficients which have to be found in order to obtain the scattered field, $\mu$ is the second Lamé constant and 
\begin{equation}
\begin{aligned}
\textstyle \sum_{rr}^{n}(z) &= \mu\left\lbrack {s_{5}{\cos\left( {p_{n}z} \right)} + s_{6}{\cos\left( {q_{n}z} \right)}} \right\rbrack,\\
{\overset{\sim}{\textstyle \sum}}_{rr}^{n}(z) &= \textstyle \sum_{r\theta}^{n}(z) = \mu\left\lbrack {s_{7}{\cos\left( {p_{n}z} \right)} + s_{8}{\cos\left( {q_{n}z} \right)}} \right\rbrack, \\
\textstyle \sum_{rz}^{n}(z) &= \mu\left\lbrack {s_{9}{\sin\left( {p_{n}z} \right)} + s_{10}{\sin\left( {q_{n}z} \right)}} \right\rbrack, \\
\end{aligned}
\label{eq2.14}
\end{equation}
\begin{equation}
s_{5} = \frac{2\left( {2p_{n}^{2} - k_{n}^{2} - q_{n}^{2}} \right)}{k_{n}}{\cos\left( {q_{n}h} \right)}, \quad 
s_{6} = \frac{2\left( {k_{n}^{2} - q_{n}^{2}} \right)}{k_{n}}{\cos\left( {p_{n}h} \right)}, \quad s_{7} = 4{\cos\left( {q_{n}h} \right)},
\label{eq2.15}
\end{equation}
\begin{equation}
s_{8} = - \frac{2\left( {k_{n}^{2} - q_{n}^{2}} \right)}{k_{n}^{2}}{\cos\left( {q_{n}h} \right)}, \quad
s_{9} = 4p_{n}{\cos\left( {q_{n}h} \right)},\quad s_{10} = \frac{\left( {k_{n}^{2} - q_{n}^{2}} \right)^{2}}{q_{n}k_{n}^{2}}{\cos\left( {p_{n}h} \right)}.
\label{eq2.16}
\end{equation}

In actual numerical computations, the summation for the angular dependence $m$ and mode $n$ are truncated at $\left| m \right| = M$ and $n = N$ respectively.

Next, the general equation for the displacement in the incident field is considered. Since the incident field is given by a symmetric plane Lamb wave and only the case of zero mode is considered, the expression for the scalar potential is expressed by the equation below.
\begin{equation}
\phi^{\text{inc}} = \text{e}^{\text{i}k_{0}x} = \text{e}^{\text{i}k_{0}r{\cos\theta}} = {\sum\limits_{m = - \infty}^{\infty}{\text{i}^{m}J_{m}\left( {k_{0}r} \right)\text{e}^{\text{i}m\theta}}}
\label{eq2.17}
\end{equation}
\noindent where $J_m\left(\right)$ is a Bessel function of the first kind.
The incident displacements and stresses can be written as
\begin{equation}
\begin{aligned}
u_{r}^{\text{inc}} &= {\sum\limits_{m = - \infty}^{\infty}{\text{i}^{m}V_{0}(z){J^{'}}_{m}\left( {k_{0}r} \right)\text{e}^{\text{i}m\theta}}}, \\
u_{\theta}^{\text{inc}} &= {\sum\limits_{m = - \infty}^{\infty}{m\text{i}^{m + 1}V_{0}(z)\frac{J_{m}\left( {k_{0}r} \right)}{k_{0}r}\text{e}^{\text{i}m\theta}}}, \\
u_{z}^{\text{inc}} &= {\sum\limits_{m = - \infty}^{\infty}{\text{i}^{m}W_{0}(z)J_{m}\left( {k_{0}r} \right)\text{e}^{\text{i}m\theta}}}.
\end{aligned}
\label{eq2.18}
\end{equation}
\begin{equation}
\begin{aligned}
\sigma_{r}^{\text{inc}} &= {\sum\limits_{m = - \infty}^{\infty}{\text{i}^{m}\left\lbrack {\textstyle \sum_{rr}^{0}(z)J_{m}\left( {k_{0}r} \right) - {\overset{\sim}{\textstyle \sum}}_{rr}^{0}(z)\left( {\frac{1}{r}{J'}_{m}\left( {k_{0}r} \right) - \frac{m^{2}}{k_{0}r^{2}}J_{m}\left( {k_{0}r} \right)} \right)} \right\rbrack\text{e}^{\text{i}m\theta}}}, \\
\tau_{r\theta}^{\text{inc}} &= {\sum\limits_{m = - \infty}^{\infty}{\text{i}^{m + 1}m\textstyle \sum_{r\theta}^{0}(z)\left( {\frac{1}{r}{J'}_{m}\left( {k_{0}r} \right) - \frac{1}{k_{0}r^{2}}J_{m}\left( {k_{0}r} \right)} \right)\text{e}^{\text{i}m\theta}}}, \\
\tau_{rz}^{\text{inc}} &= - {\sum\limits_{m = - \infty}^{\infty}{\textstyle \sum_{rz}^{0}(z)J_{m}\left( {k_{0}r} \right)\text{e}^{\text{i}m\theta}}}.
\end{aligned}
\label{eq2.19}
\end{equation}

\subsection{Far-field DtN conditions at the virtual boundaries}

In order to get the relationship between displacements and tractions at the virtual boundary, the mode orthogonality should be utilized firstly. In simplifying the equations, this study focuses on the orthogonality of each term in the right-hand side of Eq.\ref{eq2.12} and Eq.\ref{eq2.13} and proceeds by multiplying each projection function, and the procedure is described below. Firstly, the projection function in the direction of $\theta$ is $\text{e}^{- \text{i}m^{'}\theta}\left( {- M \ll m^{'} \ll M} \right)$ and the orthogonality of the projection function is
\begin{equation}
{\int_{- \pi}^{\pi}{\text{e}^{\text{i}m\theta}\text{e}^{- \text{i}m^{'}\theta}}}\text{d}\theta = \left\{ \begin{matrix}
{2\pi \quad \left( {m = m^{'}} \right)} \\
{0 \quad \left( {m \neq m^{'}} \right)} \\
\end{matrix} \right.
\label{eq2.20}
\end{equation}

Next, when considering the projection function in the $z$-direction, there is no clear function for the term of the Lamb mode. However, for the SH mode terms, only $\cos{\left (\sfrac {n\pi z}/{2h}\right)}$ for $u_r^{scat}$ and $u_\theta^{scat}$, and $\sin {\left  ( \sfrac {n\pi z}/{2h}\right )}$ for $u_z^{scat}$ can be found to be $z$-dependent functions. Therefore, the projection function in the $z$ direction is

\begin{equation}
{\int_{- h}^{h}{{\cos\left( \frac{n\pi z}{2h} \right)} \cdot {\cos\left( \frac{n^{'}\pi z}{2h} \right)}}}\text{d}\theta = \left\{ \begin{matrix}
{2h \quad \left( {n = n^{'} = 0} \right)} \\
{h \quad \left( {n = n^{'} \neq 0} \right)} \\
{0 \quad \left( {n \neq n^{'}} \right)} \\
\end{matrix} \right.,
\label{eq2.21}
\end{equation}

\begin{equation}
{\int_{- h}^{h}{{\sin \left ( \frac{n\pi z}{2h} \right )} \cdot {\sin \left ( \frac{n^{'}\pi z}{2h} \right )}}}\text{d}\theta = \left\{ \begin{matrix}
{0 \quad \left( {n = n^{'} = 0} \right)} \\
{h \quad \left( {n = n^{'} \neq 0} \right)} \\
{0 \quad \left( {n \neq n^{'}} \right)} \\
\end{matrix} \right.
\label{eq2.22}
\end{equation}

Therefore, the $z$-dependent orthogonal functions for the expressions $u_r^{{\rm{sca}}}$ and $u_\theta ^{{\rm{sca}}}$ are $\cos \left( {n'\pi z/2h} \right)$ where $n' = 0,2,4, \cdots ,N$ and the $z$-dependent orthogonal functions for the expression $u_z^{{\rm{sca}}}$ are $\sin \left( {n'\pi z/2h} \right)$ where $n' = 2,4, \cdots ,N$. Thus, the orthogonal projection functions for $u_r^{{\rm{sca}}}$ and $u_\theta ^{{\rm{sca}}}$ can be defined as
\begin{equation}
\frac{1}{2\pi\mu}{\cos\left( \frac{\acute{n}\pi z}{2h} \right)}\text{e}^{- \text{i}\acute{m}\theta}~~~~~\acute{n} = 0,2,\cdots,N
\label{eq2.23}
\end{equation}

\noindent
and the orthogonal projection functions for $u_z^{{\rm{sca}}}$ are selected as
\begin{equation}
\frac{1}{2\pi\mu}{\sin\left( \frac{n^{'}\pi z}{2h} \right)}\text{e}^{- \text{i}m^{'}\theta}~~~~~\acute{n} = 2,4,\cdots,N
\label{eq2.24}
\end{equation}

Firstly, the equation $u_r^{{\rm{sca}}}$ is transformed. Multiplying both sides of Eq.\ref{eq2.12} by Eq.\ref{eq2.23} yields
\begin{equation}
\begin{aligned}
\frac{1}{2\pi\mu}&{\cos\left( \frac{n^{'}\pi z}{2h} \right)}\text{e}^{- \text{i}m^{'}\theta}u_{r}^{scat} \\ 
= &{\sum\limits_{n = 0}^{N}{\sum\limits_{m = - M}^{M}{\frac{1}{2\pi\mu}{\cos\left( \frac{n^{'}\pi z}{2h} \right)}\text{e}^{- \text{i}m^{'}\theta}A_{mn}V_{n}(z){H^{'}}_{m}\left( {k_{n}r} \right)\text{e}^{\text{i}m\theta}}}} \\
&+ {\sum\limits_{n = 0,2,\cdots}^{N}{\sum\limits_{m = - M}^{M}{\frac{1}{2\pi\mu}{\cos\left( \frac{n^{'}\pi z}{2h} \right)}\text{e}^{- \text{i}m^{'}\theta}imB_{mn}{\cos\left( \frac{n\pi z}{2h} \right)}\frac{H_{m}\left( {l_{n}r} \right)}{l_{n}r}\text{e}^{\text{i}m\theta}}}}
\end{aligned}
\label{eq2.25}
\end{equation}

Applying a surface integral about the virtual surface of radius $a$ to both sides of the above equation, we can achieve
\begin{equation}
\begin{aligned}
{\int_{S}^{}{\frac{1}{2\pi\mu}{\cos\left( \frac{\acute{n}\pi z}{2h} \right)}\text{e}^{- \text{i}\acute{m}\theta}}}u_{r}^{\text{sca}}dS \\
= &{\sum\limits_{n = 0}^{N}{\sum\limits_{m = - M}^{M}{\int_{- \pi}^{\pi}{{\int_{- h}^{h}{\frac{1}{2\pi\mu}{\cos\left( \frac{n^{'}\pi z}{2h} \right)}V_{n}(z){H^{'}}_{m}\left( {k_{n}a} \right)\text{e}^{\text{i}{({m{- m}^{'}})}\theta}A_{mn}}}a\text{d}\theta\text{d}z}}}} \\
&+ {\sum\limits_{n = 0,2,\cdots}^{N}{\sum\limits_{m = - M}^{M}{\int_{- \pi}^{\pi}{{\int_{- h}^{h}{\frac{\text{i}m}{2\pi\mu}{\cos\left( \frac{n^{'}\pi z}{2h} \right)}{\cos\left( \frac{n\pi z}{2h} \right)}\frac{H_{m}\left( {l_{n}a} \right)}{l_{n}a}\text{e}^{\text{i}{({m{- m}^{'}})}\theta}B_{mn}}}a\text{d}\theta\text{d}z}}}}
\end{aligned}
\label{eq2.26}
\end{equation}

\begin{equation}
\left. \Longleftrightarrow{\int_{S}^{}{\frac{1}{2\pi\mu}{\cos\left( \frac{n^{'}\pi z}{2h} \right)}\text{e}^{- \text{i}m^{'}\theta}u_{r}^{\text{sca}}\text{d}S}} \right. \\
= {\sum\limits_{n = 0}^{N}{\left\lbrack {\frac{a}{\mu}{H^{'}}_{m^{'}}\left( {k_{n}a} \right){\int_{- h}^{h}{{\cos\left( \frac{n^{'}\pi z}{2h} \right)}V_{n}(z)}}\text{d}z} \right\rbrack A_{m^{'}n}}} \\
+\left\lbrack {\frac{\text{i}m^{'}H_{m^{'}}\left({l_{n^{'}}a}\right)}{\mu l_{n^{'}}}I_{n^{'}n^{'}}^{S}}\right\rbrack B_{m^{'}n^{'}}\left({n^{'}=0,~2,~\cdots,~N}\right).
\label{eq2.27}
\end{equation}

Secondly, the equation $u_{\theta}^{{\rm{sca}}}$ is transformed. Multiplying both sides of Eq.\ref{eq2.12} by Eq.\ref{eq2.23} yields
\begin{equation}
\begin{aligned}
\frac{1}{2\pi\mu}&{\cos\left( \frac{n^{'}\pi z}{2h} \right)}\text{e}^{- \text{i}m^{'}\theta}u_{\theta}^{\text{sca}} \\
= &{\sum\limits_{n = 0}^{N}{\sum\limits_{m = - M}^{M}{\frac{1}{2\pi\mu}{\cos\left( \frac{n^{'}\pi z}{2h} \right)}\text{e}^{- \text{i}m^{'}\theta}\text{i}mA_{mn}V_{n}(z)\frac{H_{m}\left( {k_{n}r} \right)}{k_{n}r}\text{e}^{\text{i}m\theta}}}} \\
&- {\sum\limits_{n = 0,2,\cdots}^{N}{\sum\limits_{m = - M}^{M}{\frac{1}{2\pi\mu}{\cos\left( \frac{n^{'}\pi z}{2h} \right)}\text{e}^{- \text{i}m^{'}\theta}B_{mn}{\cos\left( \frac{n\pi z}{2h} \right)}{H^{'}}_{m}\left( {l_{n}r} \right)\text{e}^{\text{i}m\theta}}}}
\end{aligned}
\label{eq2.28}
\end{equation}

Applying the same surface integral about the virtual surface of radius $a$ to both sides of the above equation, we can achieve
\begin{equation}
\begin{aligned}
\int_{S}^{}\frac{1}{2\pi\mu}&{\cos\left( \frac{n^{'}\pi z}{2h} \right)}\text{e}^{-\text{i}m^{'}\theta}u_{\theta}^{\text{sca}}\text{d}S\\
= &{\sum\limits_{n = 0}^{N}{\sum\limits_{m = - M}^{M}{\int_{- \pi}^{\pi}{{\int_{- h}^{h}{\frac{\text{i}m}{2\pi\mu}{\cos\left( \frac{n^{'}\pi z}{2h} \right)}V_{n}(z)\frac{H_{m}\left( {k_{n}a} \right)}{k_{n}a}\text{e}^{\text{i}{({m{- m}^{'}})}\theta}A_{mn}}}a\text{d}\theta\text{d}z}}}} \\
&- {\sum\limits_{n = 0,2,\cdots}^{N}{\sum\limits_{m = - M}^{M}{\int_{- \pi}^{\pi}{{\int_{- h}^{h}{\frac{1}{2\pi\mu}{\cos\left( \frac{n^{'}\pi z}{2h} \right)}{\cos\left( \frac{n\pi z}{2h} \right)}{H^{'}}_{m}\left( {l_{n}a} \right)\text{e}^{\text{i}{({m{- m}^{'}})}\theta}B_{mn}}}a\text{d}\theta\text{d}z}}}}
\end{aligned}
\label{eq2.29}
\end{equation}

\begin{equation}
\begin{aligned}
\Longleftrightarrow \int_{S}^{}\frac{1}{2\pi\mu}&{\cos\left( \frac{n^{'}\pi z}{2h} \right)}\text{e}^{- \text{i}m^{'}\theta}u_{\theta}^{\text{sca}}\text{d}S \\
= &{\sum\limits_{n = 0}^{N}{\left\lbrack {\frac{\text{i}m^{'}H_{m^{'}}\left( {k_{n}a} \right)}{\mu k_{n}}{\int_{- h}^{h}{{\cos\left( \frac{n^{'}\pi z}{2h} \right)}V_{n}(z)}}\text{d}z} \right\rbrack A_{m^{'}n}}} \\
&-\left\lbrack {\frac{a{H^{'}}_{m^{'}}\left({l_{n^{'}}a}\right)}{\mu}I_{nn^{'}}^{S}} \right\rbrack B_{m^{'}n^{'}}\left({n^{'}=0,~2,~\cdots,~N}\right).
\end{aligned}
\label{eq2.30}
\end{equation}

Finally, the equation $u_{z}^{{\rm{sca}}}$ is transformed. Multiplying both sides of Eq.\ref{eq2.12} by Eq.\ref{eq2.24} yields
\begin{equation}
\begin{aligned}
\frac{1}{2\pi\mu}&{\sin\left( \frac{n^{'}\pi z}{2h} \right)}\text{e}^{- \text{i}m^{'}\theta}u_{z}^{\text{sca}} \\
= &{\sum\limits_{n = 0}^{N}{\sum\limits_{m = - M}^{M}{\frac{1}{2\pi\mu}{\sin\left( \frac{n^{'}\pi z}{2h} \right)}\text{e}^{- \text{i}m^{'}\theta}A_{mn}W_{n}(z)H_{m}\left( {k_{n}r} \right)\text{e}^{\text{i}m\theta}}}}.
\end{aligned}
\label{eq2.31}
\end{equation}

Applying the same surface integral about the virtual surface of radius $a$ to both sides of the above equation, we can achieve
\begin{equation}
\begin{aligned}
\int_{S}^{}\frac{1}{2\pi\mu}&{\sin\left( \frac{n^{'}\pi z}{2h} \right)}\text{e}^{- \text{i}m^{'}\theta}u_{z}^{\text{sca}}\text{d}S \\ 
= &{\sum\limits_{n = 0}^{N}{\sum\limits_{m = - M}^{M}{\int_{- \pi}^{\pi}{{\int_{- h}^{h}{\frac{1}{2\pi\mu}{\sin\left( \frac{n^{'}\pi z}{2h} \right)}W_{n}(z)H_{m}\left( {k_{n}a} \right)\text{e}^{\text{i}{({m{- m}^{'}})}\theta}}}A_{mn}a\text{d}\theta\text{d}z}}}}
\end{aligned}
\label{eq2.32}
\end{equation}

\begin{equation}
\begin{aligned}
\Longleftrightarrow\int_{S}^{}\frac{1}{2\pi\mu}&{\sin\left( \frac{n^{'}\pi z}{2h} \right)}\text{e}^{- \text{i}m^{'}\theta}u_{z}^{\text{sca}}\text{d}S \\
=&{\sum\limits_{n=0}^{N}\left\lbrack{\frac{aH_{m^{'}}\left({k_{n}a}\right)}{\mu}{\int_{-h}^{h}{{\sin\left(\frac{n^{'}\pi z}{2h}\right)}W_{n}(z)}}\text{d}z}\right\rbrack}A_{m^{'}n}\left({n^{'}=2,~4,~\cdots,~N}\right)
\end{aligned}
\label{eq2.33}
\end{equation}

The above orthogonality technique allowed $A_{mn}$ and $B_{mn}$ to be obtained. In the above equations, $m$ can be determined to any value $m^\prime$ and $A_{mn}$, $B_{mn}$ can be considered separately for each $m^\prime$. Eq.\ref{eq2.27}, Eq.\ref{eq2.30} and Eq.\ref{eq2.33} can be rewritten into the matrix form as follows
\begin{equation}
\begin{bmatrix}
{\int_{S}^{}{\frac{1}{2\pi\mu}{\cos\left( \frac{0 \cdot \pi z}{2h} \right)}\text{e}^{- \text{i}m^{'}\theta}u_{r}^{\text{sca}}\text{d}S}} \\
 \vdots \\
{\int_{S}^{}{\frac{1}{2\pi\mu}{\cos\left( \frac{N\pi z}{2h} \right)}\text{e}^{- \text{i}m^{'}\theta}u_{r}^{\text{sca}}\text{d}S}} \\
{\int_{S}^{}{\frac{1}{2\pi\mu}{\cos\left( \frac{0 \cdot \pi z}{2h} \right)}\text{e}^{- \text{i}m^{'}\theta}u_{\theta}^{\text{sca}}\text{d}S}} \\
 \vdots \\
{\int_{S}^{}{\frac{1}{2\pi\mu}{\cos\left( \frac{N\pi z}{2h} \right)}\text{e}^{- \text{i}m^{'}\theta}u_{\theta}^{\text{sca}}\text{d}S}} \\
{\int_{S}^{}{\frac{1}{2\pi\mu}{\sin\left( \frac{2\pi z}{2h} \right)}\text{e}^{- \text{i}m^{'}\theta}u_{z}^{\text{sca}}\text{d}S}} \\
 \vdots \\
{\int_{S}^{}{\frac{1}{2\pi\mu}{\sin\left( \frac{N\pi z}{2h} \right)}\text{e}^{- \text{i}m^{'}\theta}u_{z}^{\text{sca}}\text{d}S}} \\
\end{bmatrix} = \left\lbrack {\overset{-}{AB}}_{m^{'}} \right\rbrack_{{({\frac{3}{2}N + 3})} \times {({N + 2})}}\begin{bmatrix}
A_{m^{'}0} \\
A_{m^{'}1} \\
 \vdots \\
A_{m^{'}N} \\
B_{m^{'}0} \\
B_{m^{'}2} \\
 \vdots \\
B_{m^{'}N} \\
\end{bmatrix}.
\label{eq2.34}
\end{equation}

\noindent
where
\begin{equation}
\left\lbrack {\overset{-}{AB}}_{m^{'}} \right\rbrack = \begin{bmatrix}
{\overset{-}{A}}_{m^{'}} & {\overset{-}{B}}_{m^{'}} \\
\end{bmatrix},
\label{eq2.35}
\end{equation}

\begin{equation}
\begin{aligned}
&\left\lbrack {\overset{-}{A}}_{m^{'}} \right\rbrack_{{({\frac{3}{2}N + 3})} \times {({\frac{1}{2}N + 1})}} = \\ 
&\begin{bmatrix}
{\frac{a}{\mu}{H^{'}}_{m^{'}}\left( {k_{0}a} \right){\int_{- h}^{h}{{\cos\left( \frac{0 \cdot \pi z}{2h} \right)}V_{0}(z)}}\text{d}z} & \cdots & {\frac{a}{\mu}{H^{'}}_{m^{'}}\left( {k_{N}a} \right){\int_{- h}^{h}{{\cos\left( \frac{0 \cdot \pi z}{2h} \right)}V_{N}(z)}}\text{d}z} \\
 \vdots & \ddots & \vdots \\
{\frac{a}{\mu}{H^{'}}_{m^{'}}\left( {k_{0}a} \right){\int_{- h}^{h}{{\cos\left( \frac{N\pi z}{2h} \right)}V_{0}(z)}}\text{d}z} & \cdots & {\frac{a}{\mu}{H^{'}}_{m^{'}}\left( {k_{N}a} \right){\int_{- h}^{h}{{\cos\left( \frac{N\pi z}{2h} \right)}V_{N}(z)}}\text{d}z} \\
{\frac{\text{i}m^{'}H_{m^{'}}\left( {k_{0}a} \right)}{\mu k_{0}}{\int_{- h}^{h}{{\cos\left( \frac{0 \cdot \pi z}{2h} \right)}V_{0}(z)}}\text{d}z} & \cdots & {\frac{\text{i}m^{'}H_{m^{'}}\left( {k_{N}a} \right)}{\mu k_{n}}{\int_{- h}^{h}{{\cos\left( \frac{0 \cdot \pi z}{2h} \right)}V_{N}(z)}}\text{d}z} \\
 \vdots & \ddots & \vdots \\
{\frac{\text{i}m^{'}H_{m^{'}}\left( {k_{0}a} \right)}{\mu k_{0}}{\int_{- h}^{h}{{\cos\left( \frac{N\pi z}{2h} \right)}V_{0}(z)}}\text{d}z} & \cdots & {\frac{\text{i}m^{'}H_{m^{'}}\left( {k_{N}a} \right)}{\mu k_{N}}{\int_{- h}^{h}{{\cos\left( \frac{N\pi z}{2h} \right)}V_{N}(z)}}\text{d}z} \\
{\frac{aH_{m^{'}}\left( {k_{0}a} \right)}{\mu}{\int_{- h}^{h}{{\sin\left( \frac{2\pi z}{2h} \right)}W_{0}(z)}}\text{d}z} & \cdots & {\frac{aH_{m^{'}}\left( {k_{N}a} \right)}{\mu}{\int_{- h}^{h}{{\sin\left( \frac{2\pi z}{2h} \right)}W_{N}(z)}}\text{d}z} \\
 \vdots & \ddots & \vdots \\
{\frac{aH_{m^{'}}\left( {k_{0}a} \right)}{\mu}{\int_{- h}^{h}{{\sin\left( \frac{N\pi z}{2h} \right)}W_{0}(z)}}\text{d}z} & \cdots & {\frac{aH_{m^{'}}\left( {k_{N}a} \right)}{\mu}{\int_{- h}^{h}{{\sin\left( \frac{N\pi z}{2h} \right)}W_{N}(z)}}\text{d}z} \\
\end{bmatrix},
\end{aligned}
\label{eq2.36}
\end{equation}

\begin{equation}
\left\lbrack {\overset{-}{B}}_{m^{'}} \right\rbrack_{{({\frac{3}{2}N + 3})} \times {({\frac{1}{2}N + 1})}} = \begin{bmatrix}
{\frac{\text{i}m^{'}H_{m^{'}}\left( {l_{0}a} \right)}{\mu l_{0}}I_{00}^{S}} & 0 & \cdots & 0 \\
0 & {\frac{\text{i}m^{'}H_{m^{'}}\left( {l_{2}a} \right)}{\mu l_{2}}I_{22}^{S}} & \cdots & 0 \\
 \vdots & \vdots & \ddots & \vdots \\
0 & 0 & \cdots & {\frac{\text{i}m^{'}H_{m^{'}}\left( {l_{N}a} \right)}{\mu l_{n}}I_{NN}^{S}} \\
{- \frac{a{H^{'}}_{m^{'}}\left( {l_{0}a} \right)}{\mu}I_{00}^{S}} & 0 & \cdots & 0 \\
0 & {- \frac{a{H^{'}}_{m^{'}}\left( {l_{2}a} \right)}{\mu}I_{22}^{S}} & \cdots & 0 \\
 \vdots & \vdots & \ddots & \vdots \\
0 & 0 & \cdots & {- \frac{a{H^{'}}_{m^{'}}\left( {l_{n}a} \right)}{\mu}I_{NN}^{S}} \\
0 & 0 & \cdots & 0 \\
 \vdots & \vdots & \ddots & \vdots \\
0 & 0 & \cdots & 0 \\
\end{bmatrix}.
\label{eq2.37}
\end{equation}

The actual calculation requires the left-hand side of Eq.\ref{eq2.34} to be calculated in discretised form, and the process is explained here. The integral on the left-hand side of Eq.\ref{eq2.34} is a surface integral, so it needs to be discretised using a two-dimensional shape function. Here, a four-node quadrilateral element (Q4) is adopted. Then, we can get
\begin{equation}
\begin{aligned}
\int_{S}^{}\frac{1}{2\pi\mu}&{\cos\left( \frac{n^{'}\pi z}{2h} \right)}\text{e}^{- \text{i}m^{'}\theta}u_{r}^{\text{sca}}\text{d}S \\
= &{\int_{S}^{}{\frac{1}{2\pi\mu}{\cos\left( \frac{n^{'}\pi z}{2h} \right)}\text{e}^{- \text{i}m^{'}\theta}N_{J}{u_{r}^{\text{sca}}}_{J}\text{d}S}} \\
= &\left\{ {\int_{S}^{}{\frac{1}{2\pi\mu}{\cos\left( \frac{n^{'}\pi z}{2h} \right)}\text{e}^{- \text{i}m^{'}\theta}~N_{J}\text{d}S}} \right\}{u_{r}^{\text{sca}}}_{J} \\
= &\left\{ {\int_{S}^{}{\frac{1}{2\pi\mu}{\cos\left( \frac{n^{'}\pi z}{2h} \right)}\text{e}^{- \text{i}m^{'}\theta}~N_{J}\text{d}S}} \right\}\left( {{u_{x}^{\text{sca}}}_{J}{\cos\theta} + {u_{y}^{\text{sca}}}_{J}{\sin\theta}} \right),
\end{aligned}
\label{eq2.38}
\end{equation}

\begin{equation}
\begin{aligned}
\int_{S}^{}\frac{1}{2\pi\mu}&{\cos\left( \frac{n^{'}\pi z}{2h} \right)}\text{e}^{- \text{i}m^{'}\theta}u_{\theta}^{\text{sca}}\text{d}S \\
= &{\int_{S}^{}{\frac{1}{2\pi\mu}{\cos\left( \frac{n^{'}\pi z}{2h} \right)}\text{e}^{- \text{i}m^{'}\theta}N_{J}{u_{\theta}^{\text{sca}}}_{J}\text{d}S}} \\
= &\left\{ {\int_{S}^{}{\frac{1}{2\pi\mu}{\cos\left( \frac{n^{'}\pi z}{2h} \right)}\text{e}^{- \text{i}m^{'}\theta}~N_{J}\text{d}S}} \right\}{u_{\theta}^{\text{sca}}}_{J} \\
= &\left\{ {\int_{S}^{}{\frac{1}{2\pi\mu}{\cos\left( \frac{n^{'}\pi z}{2h} \right)}\text{e}^{- \text{i}m^{'}\theta}~N_{J}\text{d}S}} \right\}\left( {- {u_{x}^{\text{sca}}}_{J}{\cos\theta} + {u_{y}^{\text{sca}}}_{J}{\sin\theta}} \right),
\end{aligned}
\label{eq2.39}
\end{equation}

\begin{equation}
\begin{aligned}
\int_{S}^{}\frac{1}{2\pi\mu}&{\sin\left( \frac{n^{'}\pi z}{2h} \right)}\text{e}^{- \text{i}m^{'}\theta}u_{z}^{\text{sca}}\text{d}S \\
= &{\int_{S}^{}{\frac{1}{2\pi\mu}{\sin\left( \frac{n^{'}\pi z}{2h} \right)}\text{e}^{- \text{i}m^{'}\theta}N_{J}{u_{z}^{\text{sca}}}_{J}\text{d}S}} \\
= &\left\{ {\int_{S}^{}{\frac{1}{2\pi\mu}{\sin\left( \frac{n^{'}\pi z}{2h} \right)}\text{e}^{- \text{i}m^{'}\theta}N_{J}\text{d}S}} \right\}{u_{z}^{\text{sca}}}_{J}.
\end{aligned}
\label{eq2.40}
\end{equation}

\noindent
where $J$ is the local node number. Projecting the integral value on the local node $J$ in Eq.\ref{eq2.38}, Eq.\ref{eq2.39} and Eq.\ref{eq2.40} onto the global node $p$, we get
\begin{equation}
\begin{aligned}
\int_{S}^{}\frac{1}{2\pi\mu}{\cos\left( \frac{n^{'}\pi z}{2h} \right)}\text{e}^{- \text{i}m^{'}\theta}~N_{J}\text{d}S \cdot {\cos\theta}&\rightarrow D_{m^{'}n^{'}}^{r~~x_{p}} = D_{m^{'}n^{'}}^{r~~p}{\cos{\theta_{p}~}},  \\
\int_{S}^{}\frac{1}{2\pi\mu}{\cos\left( \frac{n^{'}\pi z}{2h} \right)}\text{e}^{- \text{i}m^{'}\theta}~N_{J}\text{d}S \cdot {\sin\theta}&\rightarrow D_{m^{'}n^{'}}^{r~~y_{p}} = D_{m^{'}n^{'}}^{r~~p}{\sin\theta_{p}}~, \\
\int_{S}^{}\frac{1}{2\pi\mu}{\cos\left( \frac{n^{'}\pi z}{2h} \right)}\text{e}^{- \text{i}m^{'}\theta}~N_{J}\text{d}S \cdot \left( {{- \cos}\theta} \right)&\rightarrow D_{m^{'}n^{'}}^{\theta~~x_{p}} = - D_{m^{'}n^{'}}^{\theta~~p}{\sin{\theta_{p}~}}, \\
\int_{S}^{}\frac{1}{2\pi\mu}{\cos\left( \frac{n^{'}\pi z}{2h} \right)}\text{e}^{- \text{i}m^{'}\theta}~N_{J}\text{d}S \cdot {\sin\theta}&\rightarrow D_{m^{'}n^{'}}^{\theta~~y_{p}} = D_{m^{'}n^{'}}^{\theta~~p}{\cos\theta_{p}}~, \\
\int_{S}^{}\frac{1}{2\pi\mu}{\sin\left( \frac{n^{'}\pi z}{2h} \right)}\text{e}^{- \text{i}m^{'}\theta}~N_{J}\text{d}S&\rightarrow D_{m^{'}n^{'}}^{z~~z_{p}} = D_{m^{'}n^{'}}^{z~~p}~~.
\end{aligned}
\label{eq2.41}
\end{equation}

Then the left-hand side of Eq.\ref{eq2.34} can be rewritten as
\begin{equation}
\left\lbrack D_{m^{'}} \right\rbrack_{{({\frac{3}{2}N + 2})} \times 3p} \cdot \left\lbrack u_{xyz}^{\text{sca}} \right\rbrack_{3p \times 1} = \left\lbrack {\overset{-}{AB}}_{m^{'}} \right\rbrack_{{({\frac{3}{2}N + 2})} \times {({\frac{3}{2}N + 2})}} \cdot \begin{bmatrix}
A_{m^{'}} \\
B_{m^{'}} \\
\end{bmatrix}_{{({\frac{3}{2}N + 2})} \times 1}
\label{eq2.42}
\end{equation}

\noindent
where
\begin{equation}
\left\lbrack D_{m^{'}} \right\rbrack_{{({\frac{3}{2}N + 2})} \times 3p} = \begin{bmatrix}
D_{m^{'}0}^{r~~x_{1}} & D_{m^{'}0}^{r~~y_{1}} & 0 & \cdots & D_{m^{'}0}^{r~~x_{p}} & D_{m^{'}0}^{r~~y_{p}} & 0 \\
 \vdots & \vdots & \vdots & \ddots & \vdots & \vdots & \vdots \\
D_{m^{'}N}^{r~~x_{1}} & D_{m^{'}N}^{r~~y_{1}} & 0 & \cdots & D_{m^{'}N}^{r~~x_{p}} & D_{m^{'}N}^{r~~y_{p}} & 0 \\
D_{m^{'}0}^{\theta~~x_{1}} & D_{m^{'}0}^{\theta~~y_{1}} & 0 & \cdots & D_{m^{'}0}^{\theta~~x_{p}} & D_{m^{'}0}^{\theta~~y_{p}} & 0 \\
 \vdots & \vdots & \vdots & \ddots & \vdots & \vdots & \vdots \\
D_{m^{'}N}^{\theta~~x_{1}} & D_{m^{'}N}^{\theta~~y_{1}} & 0 & \cdots & D_{m^{'}N}^{\theta~~x_{p}} & D_{m^{'}N}^{\theta~~y_{p}} & 0 \\
0 & 0 & D_{m^{'}2}^{z~~z_{1}} & \cdots & 0 & 0 & D_{m^{'}2}^{z~~z_{p}} \\
 \vdots & \vdots & \vdots & \ddots & \vdots & \vdots & \vdots \\
0 & 0 & D_{m^{'}N}^{z~~z_{1}} & \cdots & 0 & 0 & D_{m^{'}N}^{z~~z_{p}} \\
\end{bmatrix}~,
\label{eq2.43}
\end{equation}

\begin{equation}
\left\lbrack u_{xyz}^{\text{sca}} \right\rbrack_{3p \times 1} = \begin{bmatrix}
u_{x_{1}}^{\text{sca}} \\
u_{y_{1}}^{\text{sca}} \\
u_{z_{1}}^{\text{sca}} \\
u_{x_{2}}^{\text{sca}} \\
 \vdots \\
u_{z_{p}}^{\text{sca}} \\
\end{bmatrix}~.
\label{eq2.44}
\end{equation}

According to Eq.\ref{eq2.42}, the unknown scattered coefficients for each $m'$ can be expressed as
\begin{equation}
\begin{bmatrix}
A_{m^{'}} \\
B_{m^{'}} \\
\end{bmatrix}_{{({\frac{3}{2}N + 2})} \times 1} = \left\lbrack {\overset{-}{AB}}_{m^{'}} \right\rbrack_{{({\frac{3}{2}N + 2})} \times {({\frac{3}{2}N + 2})}}^{- 1} \cdot \left\lbrack D_{m^{'}} \right\rbrack_{{({\frac{3}{2}N + 2})} \times 3p} \cdot \left\lbrack u_{xyz}^{\text{sca}} \right\rbrack_{3p \times 1}
\label{eq2.45}
\end{equation}

Next, the relationship between the unknown scattered coefficients and the node stresses at the virtual boundary will be established. Due to Eq.\ref{eq2.13}, the stress expressions at the virtual boundary can be rewritten as
\begin{equation}
\begin{aligned}
\sigma_{rr}^{scat} &= {\sum\limits_{m^{'} = - M}^{M}{\left\lbrack {E_{rr}^{m^{'}}\left( {\theta,z} \right)} \right\rbrack\begin{bmatrix}
A_{m^{'}} \\
B_{m^{'}} \\
\end{bmatrix}}} = {\sum\limits_{m^{'} = - M}^{M}{\left\lbrack {E_{rr}^{m^{'}}\left( {\theta,z} \right)} \right\rbrack\left\lbrack {\overset{-}{AB}}_{m^{'}} \right\rbrack^{- 1}\left\lbrack D_{m^{'}} \right\rbrack\left\lbrack u_{xyz}^{\text{sca}} \right\rbrack}}\mathbf{~}, \\
\sigma_{r\theta}^{scat} &= {\sum\limits_{m^{'} = - M}^{\mathbf{M}}{\left\lbrack {E_{r\theta}^{m^{'}}\left( {\theta,z} \right)} \right\rbrack\begin{bmatrix}
A_{m^{'}} \\
B_{m^{'}} \\
\end{bmatrix}}} = {\sum\limits_{m^{'} = - M}^{M}{\left\lbrack {E_{r\theta}^{m^{'}}\left( {\theta,z} \right)} \right\rbrack\left\lbrack {\overset{-}{AB}}_{m^{'}} \right\rbrack^{- 1}\left\lbrack D_{m^{'}} \right\rbrack\left\lbrack u_{xyz}^{\text{sca}} \right\rbrack}}\mathbf{~}, \\
\sigma_{rz}^{scat} &= {\sum\limits_{m^{'} = - M}^{M}{\left\lbrack {E_{rz}^{m^{'}}\left( {\theta,z} \right)} \right\rbrack\begin{bmatrix}
A_{m^{'}} \\
B_{m^{'}} \\
\end{bmatrix}}} = {\sum\limits_{m^{'} = - M}^{M}{\left\lbrack {E_{rz}^{m^{'}}\left( {\theta,z} \right)} \right\rbrack\left\lbrack {\overset{-}{AB}}_{m^{'}} \right\rbrack^{- 1}\left\lbrack D_{m^{'}} \right\rbrack\left\lbrack u_{xyz}^{\text{sca}} \right\rbrack}}\mathbf{~}.
\end{aligned}
\label{eq2.46}
\end{equation}

\noindent
where
\begin{equation}
\begin{aligned}
\left\lbrack E_{rr}^{m^{'}} \right\rbrack &= \begin{bmatrix}
E_{rrA}^{m^{'}0}&E_{rrA}^{m^{'}1}&\cdots&E_{rrA}^{m^{'}n}&E_{rrB}^{m^{'}0}&E_{rrB}^{m^{'}2}&\cdots&E_{rrB}^{m^{'}n}\\
\end{bmatrix}\mathbf{~}, \\
\left\lbrack E_{r\theta}^{m^{'}} \right\rbrack &= \begin{bmatrix}
E_{r\theta A}^{m^{'}0} & E_{r\theta A}^{m^{'}1} & \cdots & E_{r\theta A}^{m^{'}n} & E_{r\theta B}^{m^{'}0} & E_{r\theta B}^{m^{'}2} & \cdots & E_{r\theta B}^{m^{'}n} \\
\end{bmatrix}\mathbf{~}, \\
\left\lbrack E_{rz}^{m^{'}} \right\rbrack &= \begin{bmatrix}
E_{rz A}^{m^{'}0} & E_{rz A}^{m^{'}1} & \cdots & E_{rzA}^{m^{'}n}&0&E_{rzB}^{m^{'}2}&\cdots&E_{rzB}^{m^{'}n}\\
\end{bmatrix}\mathbf{~},
\end{aligned}
\label{eq2.47}
\end{equation}

\begin{equation}
\begin{aligned}
E_{rrA}^{mn}&=\left\lbrack{\Sigma_{rr}^{n}(z)H_{m}\left({k_{n}a}\right)-{\overset{\sim}{\Sigma}}_{rr}^{n}(z)\left({\frac{1}{a}{H^{'}}_{m}\left({k_{n}a}\right)-\frac{m^{2}}{k_{n}a^{2}}H_{m}\left({k_{n}a}\right)}\right)}\right\rbrack\text{e}^{\text{i}m\theta},\\
E_{rrB}^{mn}&=\text{i}m\mu{\cos\left(\frac{n\pi z}{2h}\right)}\left({\frac{2}{a}{H^{'}}_{m}\left({l_{n}a}\right)-\frac{2}{l_{n}a^{2}}H_{m}\left({l_{n}a}\right)}\right)\text{e}^{\text{i}m\theta},\\
E_{r\theta A}^{mn}&=\text{i}m\Sigma_{r\theta}^{n}(z)\left( {\frac{1}{a}{H^{'}}_{m}\left( {k_{n}a} \right) - \frac{1}{k_{n}a^{2}}H_{m}\left( {k_{n}a} \right)} \right)\text{e}^{\text{i}m\theta}, \\
E_{r\theta B}^{mn} &= \mu{\cos\left( \frac{n\pi z}{2h} \right)}\left\lbrack {\frac{2}{a}{H^{'}}_{m}\left( {l_{n}a} \right) + \left( {l_{n} - \frac{2m^{2}}{l_{n}a^{2}}} \right)H_{m}\left( {l_{n}a} \right)} \right\rbrack\text{e}^{\text{i}m\theta}, \\
E_{rzA}^{mn}&=-\Sigma_{rz}^{n}(z){H^{'}}_{m}\left({k_{n}a}\right)\text{e}^{\text{i}m\theta},\\
E_{rzB}^{mn}&= - \text{i}m\mu\frac{n\pi}{2h}{\sin\left( \frac{n\pi z}{2h}\right)}\frac{1}{l_{n}a}H_{m}\left( {l_{n}a} \right)\text{e}^{\text{i}m\theta}.
\end{aligned}
\label{eq2.48}
\end{equation}

Therefore, the equivalent nodal forces at the virtual boundary can be denoted as
\begin{equation}
\begin{aligned}
F_{x}^{\text{sca}} &= {\int_{S}^{}{N_{J}T_{x}^{\text{sca}}}}\text{d}S = {\int_{S}^{}{N_{J}\left( {\sigma_{r}^{\text{sca}}{\cos\theta} - \tau_{r\theta}^{\text{sca}}{\sin\theta}} \right)}}\text{d}S \\
&= {\int_{S}^{}{N_{J}\left\{ {{\sum\limits_{m^{'} = - M}^{M}\left( {\left\lbrack {E_{rr}^{m^{'}}\left( {\theta,z} \right)} \right\rbrack\begin{bmatrix}
A_{m^{'}} \\
B_{m^{'}} \\
\end{bmatrix}} \right)}{\cos\theta} - {\sum\limits_{m^{'} = - M}^{M}\left( {\left\lbrack {E_{r\theta}^{m^{'}}\left( {\theta,z} \right)} \right\rbrack\begin{bmatrix}
A_{m^{'}} \\
B_{m^{'}} \\
\end{bmatrix}} \right)}{\sin\theta}} \right\}}}\text{d}S \\
&= {\int_{S}^{}{N_{J}\begin{Bmatrix}
{{\sum\limits_{m^{'} = - M}^{M}\left( {\left\lbrack {E_{rr}^{m^{'}}\left( {\theta,z} \right)} \right\rbrack\left\lbrack {\overset{-}{AB}}_{m^{'}} \right\rbrack^{- 1}\left\lbrack D_{m^{'}} \right\rbrack\left\lbrack u_{xyz}^{\text{sca}} \right\rbrack} \right)}{\cos\theta}} \\
{- {\sum\limits_{m^{'} = - M}^{M}\left( {\left\lbrack {E_{r\theta}^{m^{'}}\left( {\theta,z} \right)} \right\rbrack\left\lbrack {\overset{-}{AB}}_{m^{'}} \right\rbrack^{- 1}\left\lbrack D_{m^{'}} \right\rbrack\left\lbrack u_{xyz}^{\text{sca}} \right\rbrack} \right)}{\sin\theta}} \\
\end{Bmatrix}}}\text{d}S,
\end{aligned}
\label{eq2.49}
\end{equation}

\begin{equation}
\begin{aligned}
F_{y}^{\text{sca}} &= {\int_{S}^{}{N_{J}T_{y}^{\text{sca}}}}\text{d}S = {\int_{S}^{}{N_{J}\left( {\sigma_{r}^{\text{sca}}{\sin\theta} + \tau_{r\theta}^{\text{sca}}{\cos\theta}} \right)}}\text{d}S \\
&= {\int_{S}^{}{N_{J}\left\{ {{\sum\limits_{m^{'} = - M}^{M}\left( {\left\lbrack {E_{rr}^{m^{'}}\left( {\theta,z} \right)} \right\rbrack\begin{bmatrix}
A_{m^{'}} \\
B_{m^{'}} \\
\end{bmatrix}} \right)}{\sin\theta} + {\sum\limits_{m^{'} = - M}^{M}\left( {\left\lbrack {E_{r\theta}^{m^{'}}\left( {\theta,z} \right)} \right\rbrack\begin{bmatrix}
A_{m^{'}} \\
B_{m^{'}} \\
\end{bmatrix}} \right)}{\cos\theta}} \right\}}}\text{d}S \\
&= {\int_{S}^{}{N_{J}\begin{Bmatrix}
{{\sum\limits_{m^{'} = - M}^{M}\left( {\left\lbrack {E_{rr}^{m^{'}}\left( {\theta,z} \right)} \right\rbrack\left\lbrack {\overset{-}{AB}}_{m^{'}} \right\rbrack^{- 1}\left\lbrack D_{m^{'}} \right\rbrack\left\lbrack u_{xyz}^{\text{sca}} \right\rbrack} \right)}{\sin\theta}} \\
{+ {\sum\limits_{m^{'} = - M}^{M}\left( {\left\lbrack {E_{r\theta}^{m^{'}}\left( {\theta,z} \right)} \right\rbrack\left\lbrack {\overset{-}{AB}}_{m^{'}} \right\rbrack^{- 1}\left\lbrack D_{m^{'}} \right\rbrack\left\lbrack u_{xyz}^{\text{sca}} \right\rbrack} \right)}{\cos\theta}} \\
\end{Bmatrix}}}\text{d}S,
\end{aligned}
\label{eq2.50}
\end{equation}

\begin{equation}
\begin{aligned}
F_{z}^{\text{sca}} &= {\int_{S}^{}{N_{J}T_{z}^{\text{sca}}}}\text{d}S = {\int_{S}^{}{N_{J}\tau_{rz}^{\text{sca}}}}\text{d}S \\
&= {\int_{S}^{}{N_{J}{\sum\limits_{m^{'} = - M}^{M}\left( {\left\lbrack {E_{rz}^{m^{'}}\left( {\theta,z} \right)} \right\rbrack\begin{bmatrix}
A_{m^{'}} \\
B_{m^{'}} \\
\end{bmatrix}} \right)}}}\text{d}S \\
&= {\int_{S}^{}{N_{J}{\sum\limits_{m^{'} = - M}^{M}\left( {\left\lbrack {E_{rz}^{m^{'}}\left( {\theta,z} \right)} \right\rbrack\left\lbrack {\overset{-}{AB}}_{m^{'}} \right\rbrack^{- 1}\left\lbrack D_{m^{'}} \right\rbrack\left\lbrack u_{xyz}^{\text{sca}} \right\rbrack} \right)}}}\text{d}S.
\end{aligned}
\label{eq2.51}
\end{equation}

Projecting the integral value on the local node $J$ in Eq.\ref{eq2.49}, Eq.\ref{eq2.50} and Eq.\ref{eq2.51} onto the global node $p$, we get
\begin{equation}
\begin{aligned}
\int_{S}^{}N_{J}\left( {\left\lbrack {E_{rr}^{m^{'}}\left( {\theta,z} \right)} \right\rbrack{\cos{\theta - \left\lbrack {E_{r\theta}^{m^{'}}\left( {\theta,z} \right)} \right\rbrack{\sin\theta}}}} \right)\text{d}S&\rightarrow\left\lbrack G_{m^{'}}^{x_{p}} \right\rbrack \\
\int_{S}^{}N_{J}\left( {\left\lbrack {E_{rr}^{m^{'}}\left( {\theta,z} \right)} \right\rbrack{\sin{\theta + \left\lbrack {E_{r\theta}^{m^{'}}\left( {\theta,z} \right)} \right\rbrack{\cos\theta}}}} \right)\text{d}S&\rightarrow\left\lbrack G_{m^{'}}^{y_{p}} \right\rbrack \\
\int_{S}^{}N_{J}\left\lbrack {E_{rz}^{m^{'}}\left( {\theta,z} \right)} \right\rbrack \text{d}S&\rightarrow\left\lbrack G_{m^{'}}^{z_{p}} \right\rbrack
\end{aligned}
\label{eq2.52}
\end{equation}

Finally, the equivalent nodal force matrix can be expressed as
\begin{equation}
\left\lbrack F_{xyz}^{\text{sca}} \right\rbrack_{3p \times 1} = \left\lbrack {\overset{-}{F}}^{sca} \right\rbrack_{3p \times 3p} \cdot \left\lbrack u_{xyz}^{\text{sca}} \right\rbrack_{3p \times 1}
\label{eq2.53}
\end{equation}

\noindent
where
\begin{equation}
\left\lbrack {\overset{-}{F}}^{sca} \right\rbrack_{3p \times 3p} = {\sum\limits_{m^{'} = - M}^{M}{\left\lbrack G_{m^{'}} \right\rbrack_{3p \times {({\frac{3}{2}N + 2})}}\left\lbrack {\overset{-}{AB}}_{m^{'}} \right\rbrack_{{({\frac{3}{2}N + 2})} \times {({\frac{3}{2}N + 2})}}^{- 1} \left\lbrack D_{m^{'}} \right\rbrack_{{({\frac{3}{2}N + 2})} \times 3p}}}
\label{eq2.54}
\end{equation}

\begin{equation}
\left\lbrack G_{m^{'}} \right\rbrack_{3p \times {({\frac{3}{2}N + 2})}} = \begin{bmatrix}
\begin{matrix}
G_{m^{'}}^{x_{1}} \\
G_{m^{'}}^{y_{1}} \\
G_{m^{'}}^{z_{1}} \\
\end{matrix} \\
 \vdots \\
\begin{matrix}
G_{m^{'}}^{x_{p}} \\
G_{m^{'}}^{y_{p}} \\
G_{m^{'}}^{z_{p}} \\
\end{matrix} \\
\end{bmatrix}
\label{eq2.55}
\end{equation}

\subsection{Three-dimensional DtN-FEM formulation}

The elasto-dynamic finite element formulation is 
\begin{equation}
\begin{gathered}
\left( {\left[ K \right] - {\omega ^2}\left[ M \right]} \right)\left[ {{U^{{\rm{tot}}}}} \right] = \left[ {{F^{{\rm{tot}}}}} \right]
\end{gathered}
\label{eq2.56}
\end{equation}

\noindent
where $\omega $ is the angular frequency, and

\begin{equation}
\begin{gathered}
\left[ K \right] = \int\limits_S {{{\left[ B \right]}^T}\left[ D \right]\left[ B \right]{\rm{d}}S} ,\quad \left[ M \right] = \rho \int\limits_S {{{\left[ N \right]}^T}\left[ N \right]{\rm{d}}S} , \\
\left[ {{U^{{\rm{tot}}}}} \right] = \left[ {{U^{{\rm{inc}}}}} \right] + \left[ {{U^{{\rm{sca}}}}} \right],\quad \left[ {{F^{{\rm{tot}}}}} \right] = \left[ {{F^{{\rm{inc}}}}} \right] + \left[ {{F^{{\rm{sca}}}}} \right].
\end{gathered}
\label{eq2.57}
\end{equation}

\noindent
and ${\left[ N \right]}$ is the matrix of shape function, ${\left[ B \right]}$ is the geometric matrix, ${\left[ D \right]}$ is the elastic matrix, ${\left[ K \right]}$ is the stiffness matrix, ${\left[ M \right]}$ is the mass matrix, and $\rho $ is the material density.

Utilizing Eq.\ref{eq2.53}, Eq.\ref{eq2.56} can be finally simplified as
\begin{equation}
\begin{gathered}
\left[ {{U^{{\rm{sca}}}}} \right] = {\left\{ {\left( {\left[ K \right] - {\omega ^2}\left[ M \right]} \right) - \left[ {{{\bar F}^{{\rm{sca}}}}} \right]} \right\}^{ - 1}}\left\{ {\left[ {{F^{{\rm{inc}}}}} \right] - \left( {\left[ K \right] - {\omega ^2}\left[ M \right]} \right)\left[ {{U^{{\rm{inc}}}}} \right]} \right\}
\end{gathered}
\label{eq2.58}
\end{equation}

\noindent
which is the basic DtN-FEM formulation.

\section{Numerical examples}

In this section, the validity will be presented in order to verify the correctness and efficiency of our proposed DtN-FEM. Here, the DtN-FEM results will be compared with the BEM results. Furthermore, parametric analysis about the effect on different defect shapes will be discussed detailedly.

\begin{figure}[h]
    \centering
    \subfigure[;]{
    \begin{minipage}{0.32\linewidth}
        \centering
        \includegraphics[width=5cm]{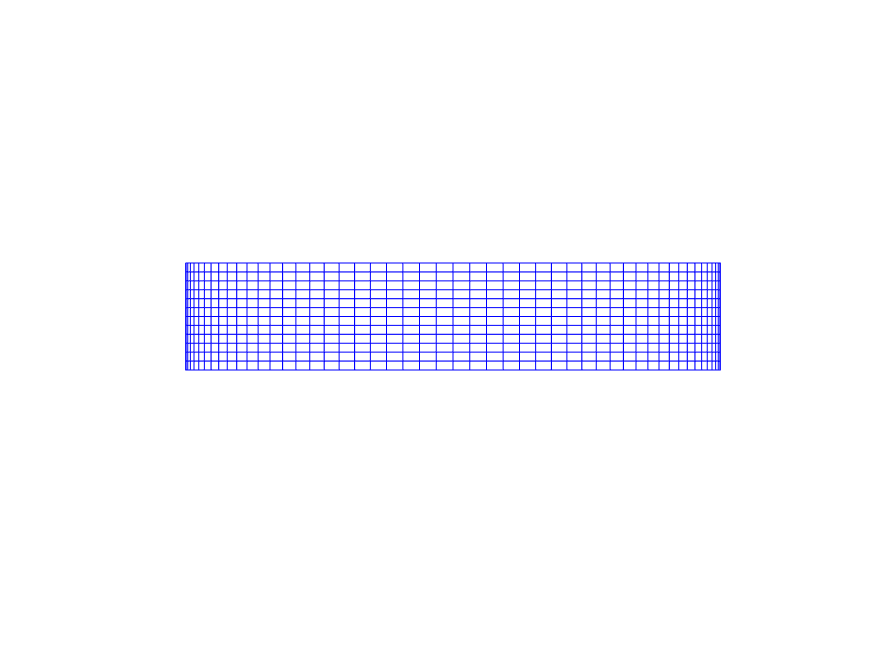}
    \end{minipage}}
    \hfill
    \subfigure[.]{
    \begin{minipage}{0.32\linewidth}
        \centerline{\includegraphics[width=5cm]{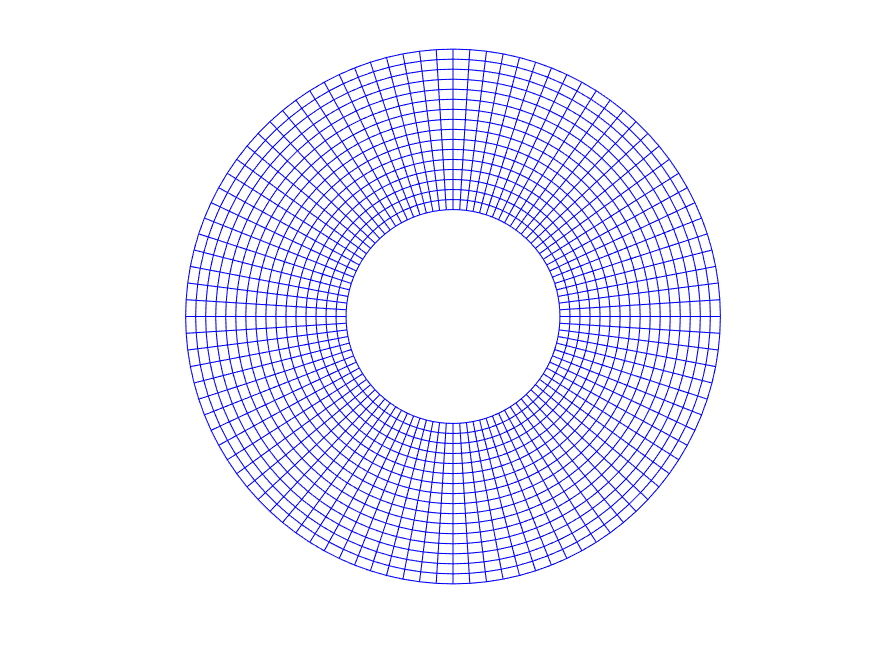}}
    \end{minipage}}
    \hfill
    \subfigure[.]{
    \begin{minipage}{0.32\linewidth}
        \centerline{\includegraphics[width=5cm]{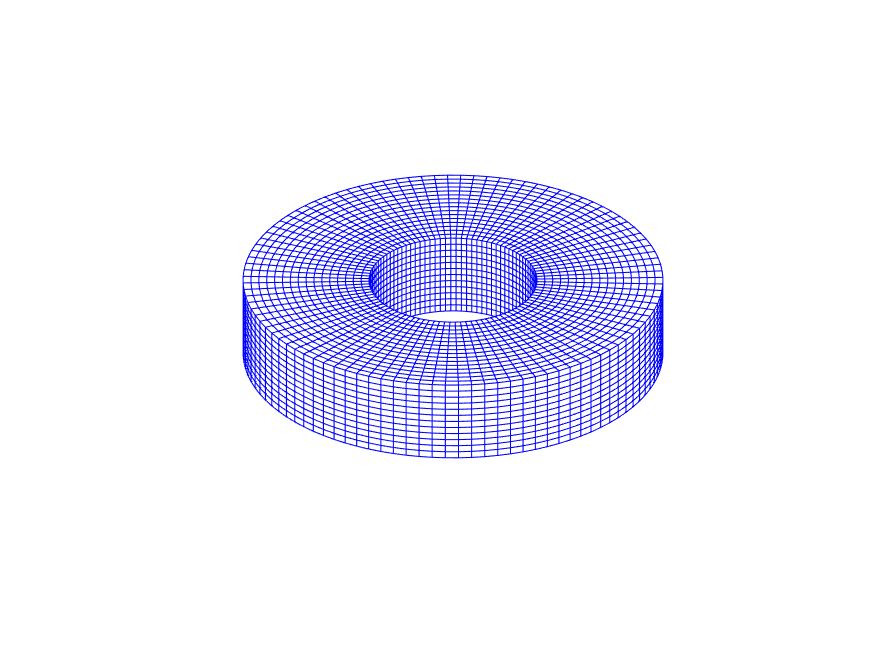}}
    \end{minipage}}
    \caption{aa}
    \label{fig3.1}
\end{figure}

\begin{figure}[h]
\centering
\subfigure[.]{\includegraphics[scale=0.56]{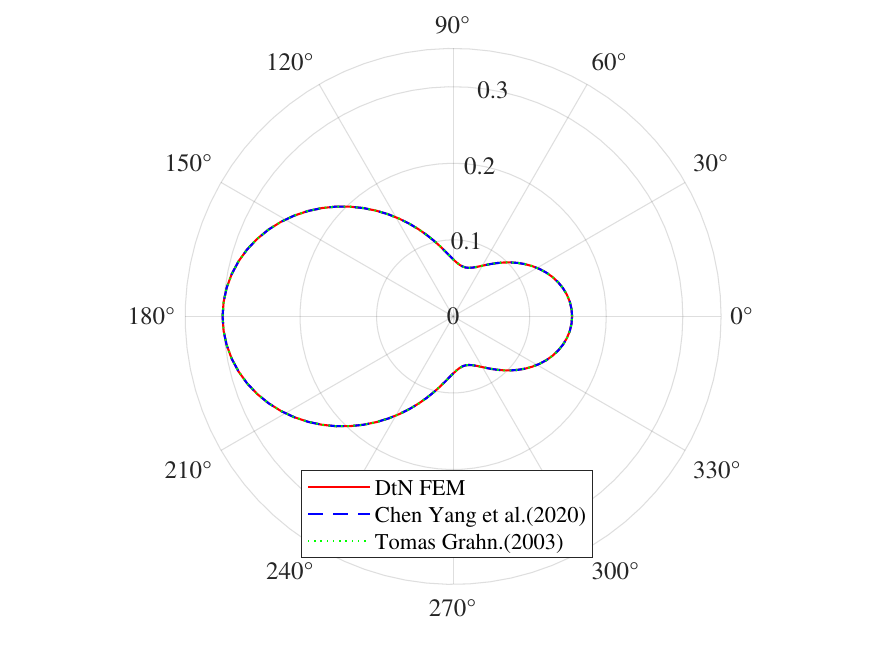}}
\subfigure[.]{\includegraphics[scale=0.56]{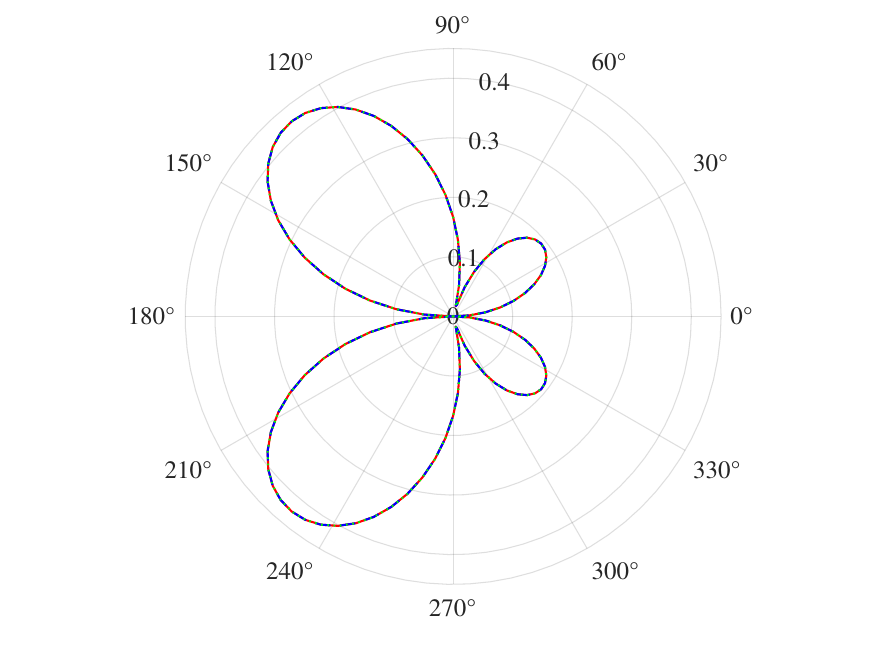}}
\caption{Error analysis of energy balance}
\label{fig2}
\end{figure}

\begin{figure}[hp]
    \centering
    \subfigure[;]{
    \begin{minipage}{0.32\linewidth}
        \centerline{\includegraphics[width=6cm]{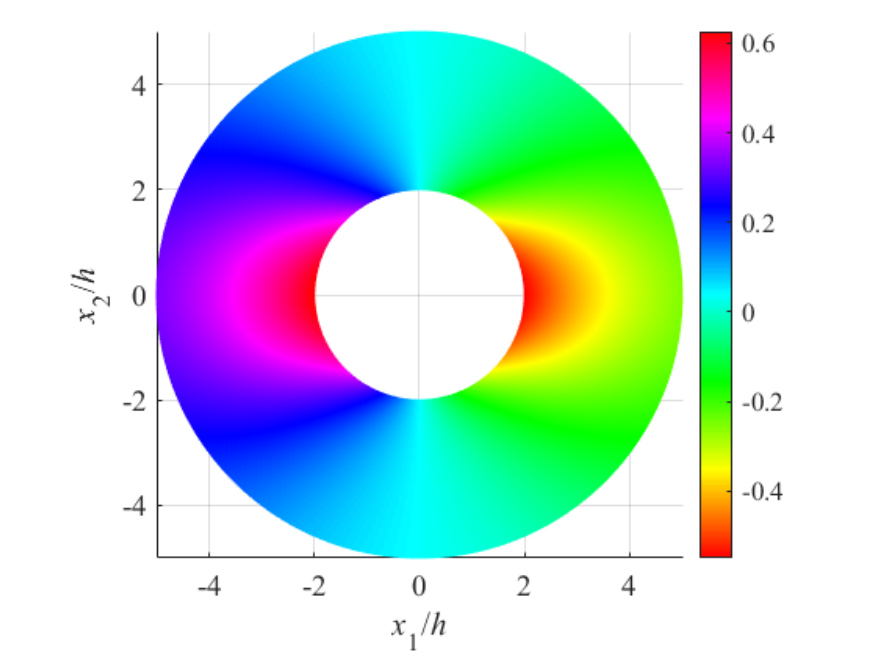}}
    \end{minipage}}
    \hfill
    \subfigure[.]{
    \begin{minipage}{0.32\linewidth}
        \centerline{\includegraphics[width=6cm]{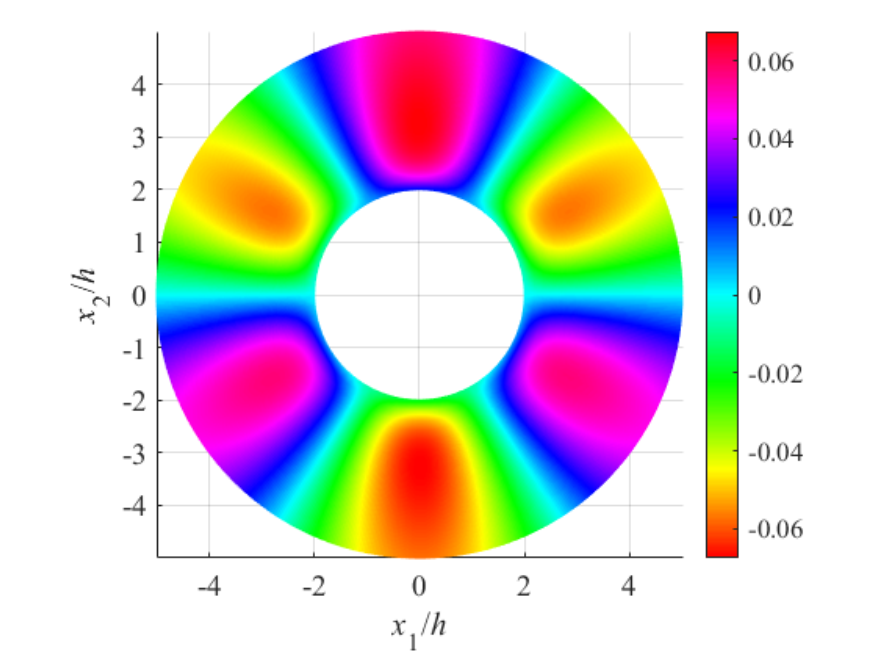}}
    \end{minipage}}
    \hfill
    \subfigure[.]{
    \begin{minipage}{0.32\linewidth}
        \centerline{\includegraphics[width=6cm]{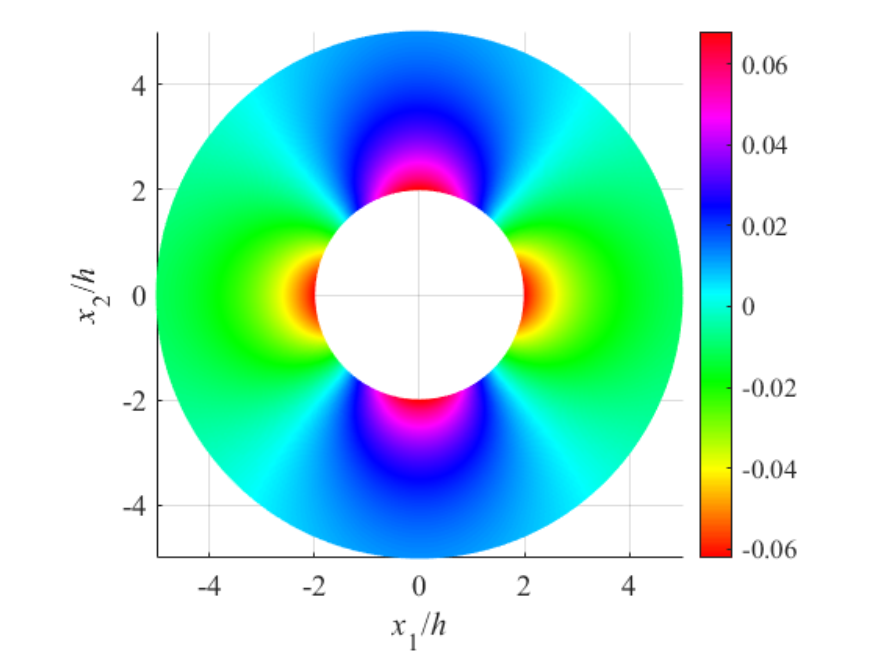}}
    \end{minipage}}
    \hfill
    \subfigure[.]{
    \begin{minipage}{0.32\linewidth}
        \centerline{\includegraphics[width=6cm]{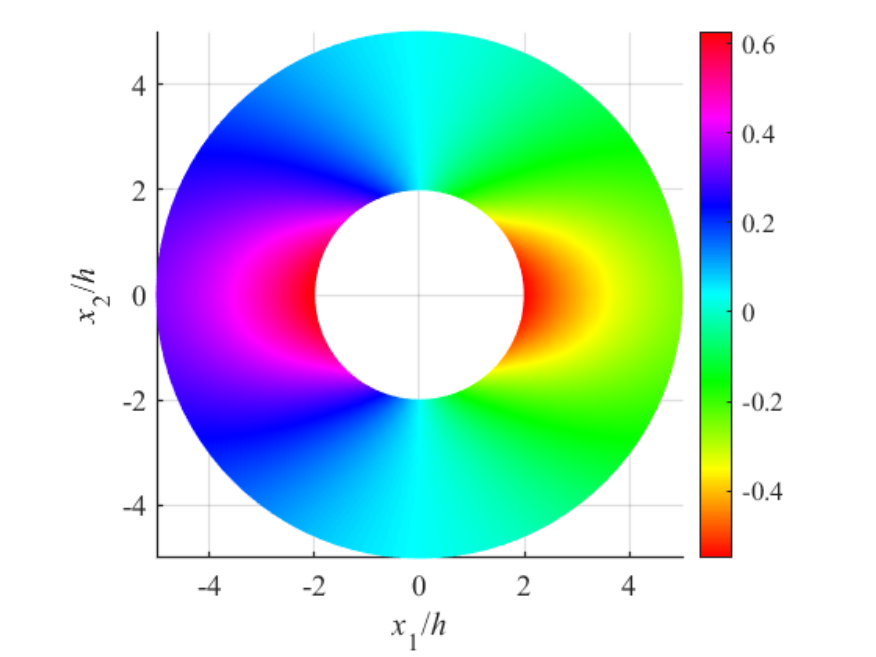}}
    \end{minipage}}
    \hfill
    \subfigure[.]{
    \begin{minipage}{0.32\linewidth}
        \centerline{\includegraphics[width=6cm]{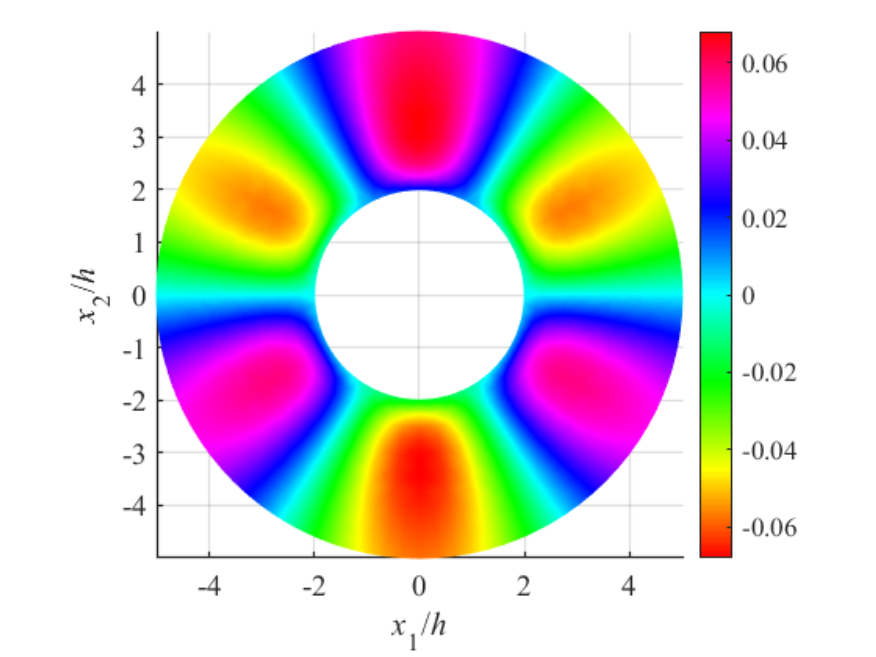}}
    \end{minipage}}
    \hfill
    \subfigure[.]{
    \begin{minipage}{0.32\linewidth}
        \centerline{\includegraphics[width=6cm]{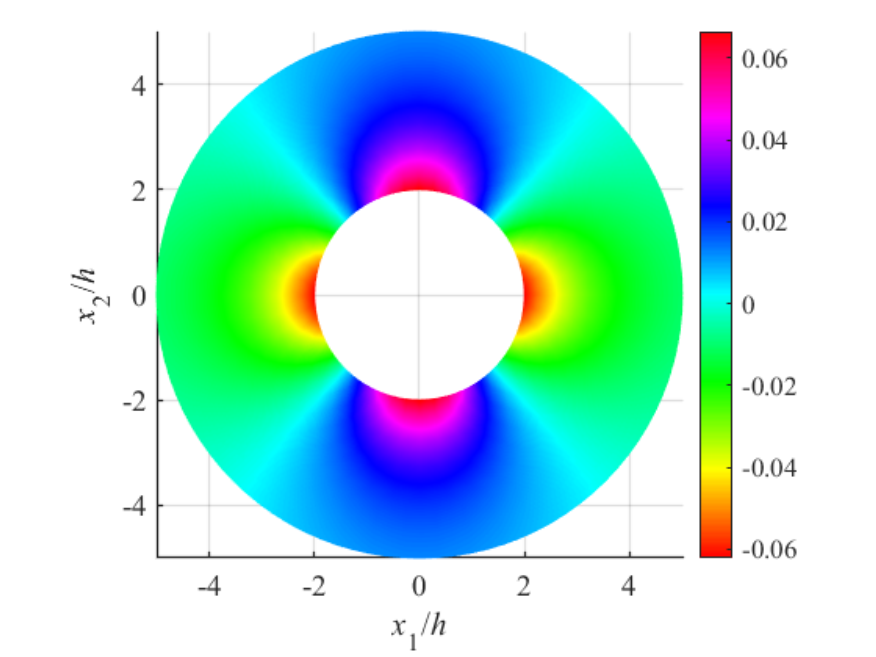}}
    \end{minipage}}
    \hfill
    \subfigure[.]{
    \begin{minipage}{0.32\linewidth}
        \centerline{\includegraphics[width=6cm]{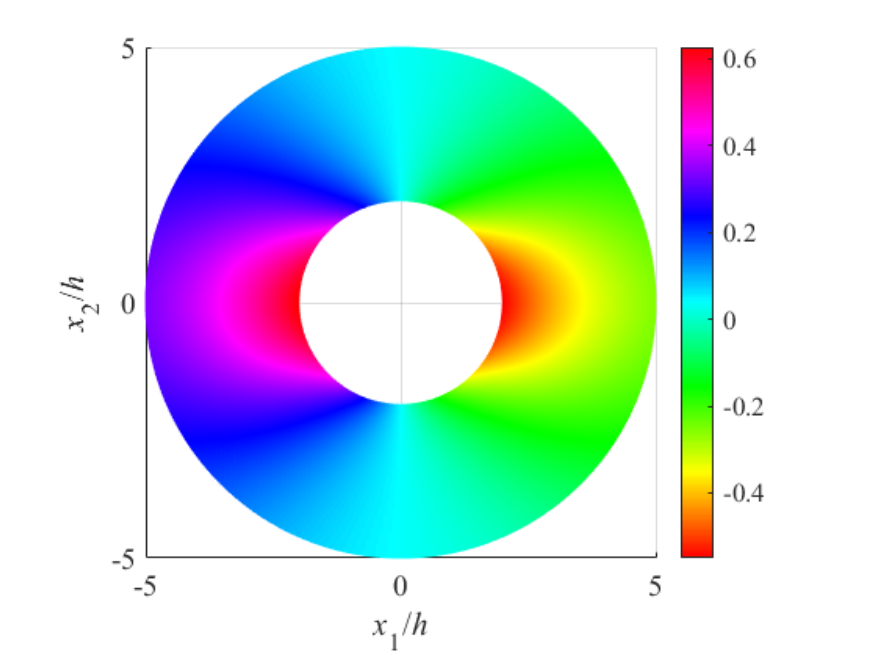}}
    \end{minipage}}
    \hfill
    \subfigure[.]{
    \begin{minipage}{0.32\linewidth}
        \centerline{\includegraphics[width=6cm]{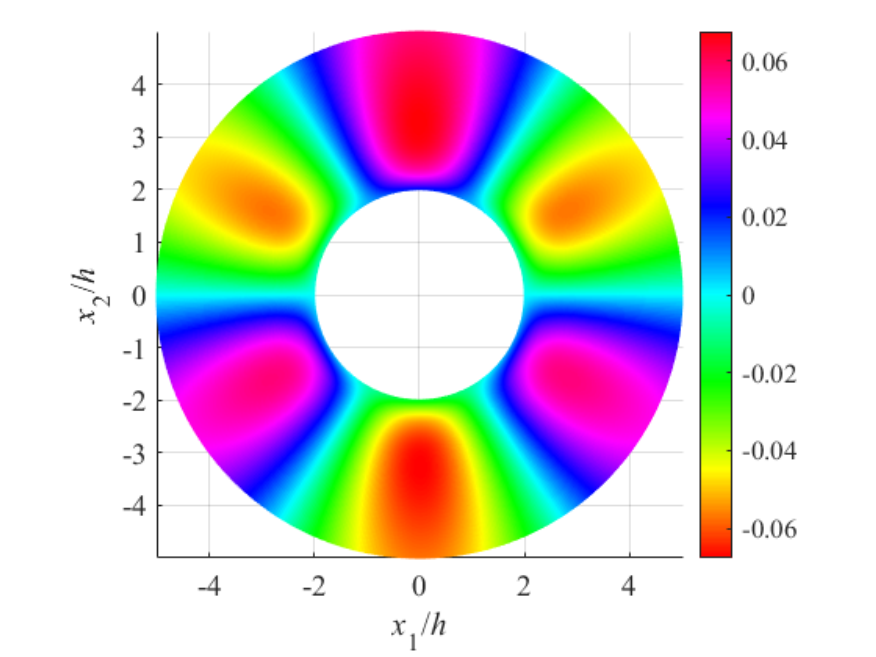}}
    \end{minipage}}
    \hfill
    \subfigure[.]{
    \begin{minipage}{0.32\linewidth}
        \centerline{\includegraphics[width=6cm]{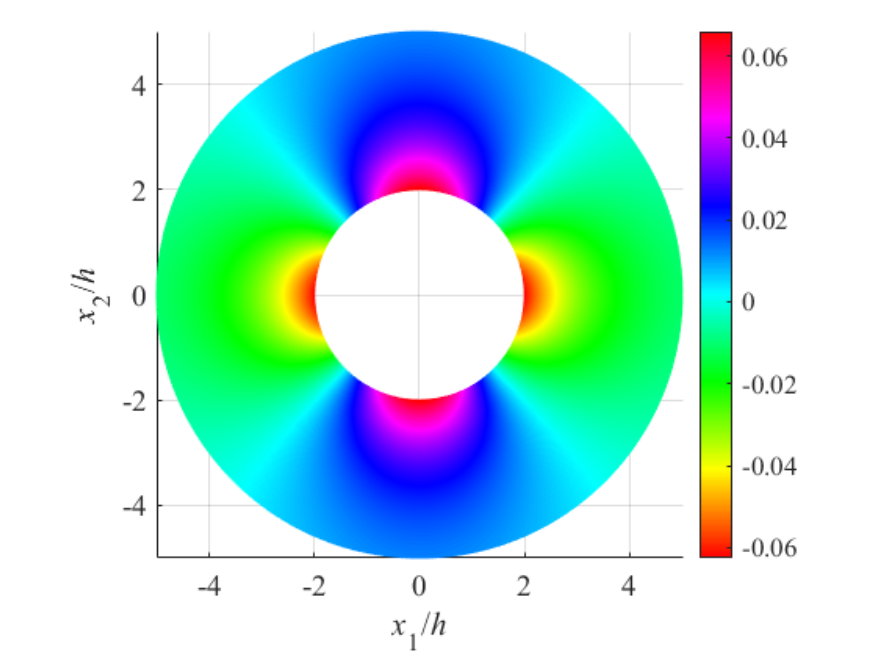}}
    \end{minipage}}
    \hfill
    \subfigure[.]{
    \begin{minipage}{0.32\linewidth}
        \centerline{\includegraphics[width=6cm]{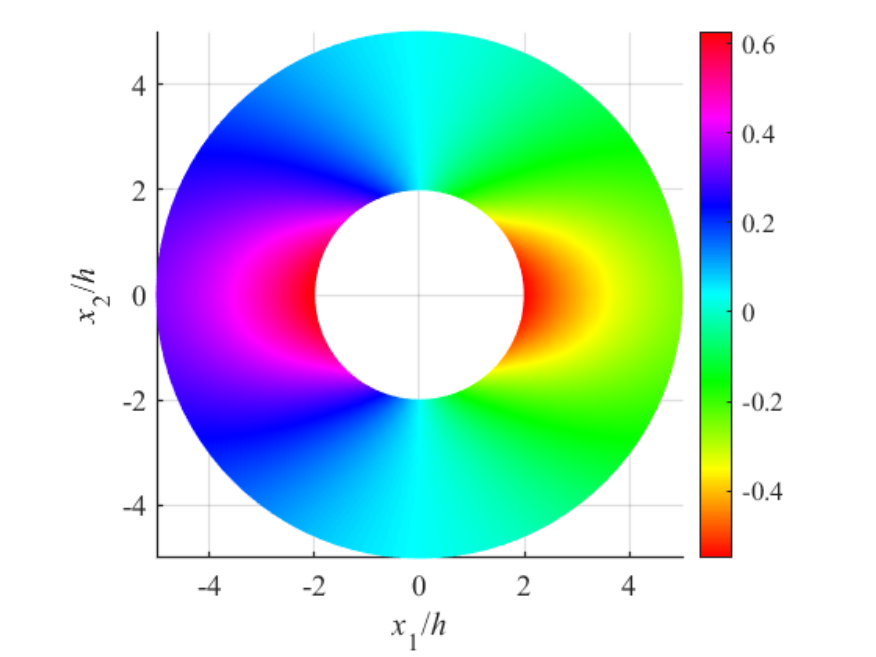}}
    \end{minipage}}
    \hfill
    \subfigure[.]{
    \begin{minipage}{0.32\linewidth}
        \centerline{\includegraphics[width=6cm]{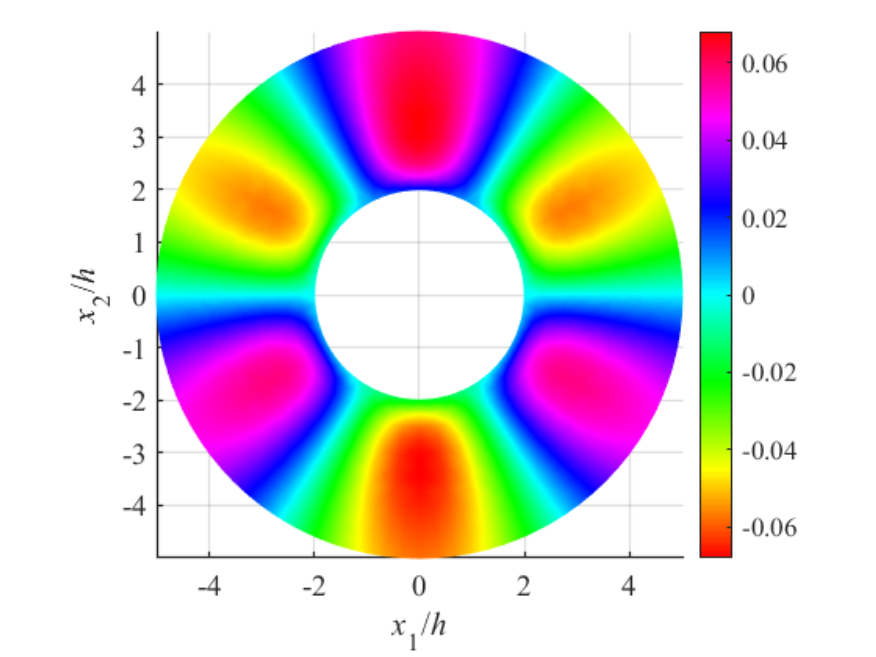}}
    \end{minipage}}
    \hfill
    \subfigure[.]{
    \begin{minipage}{0.32\linewidth}
        \centerline{\includegraphics[width=6cm]{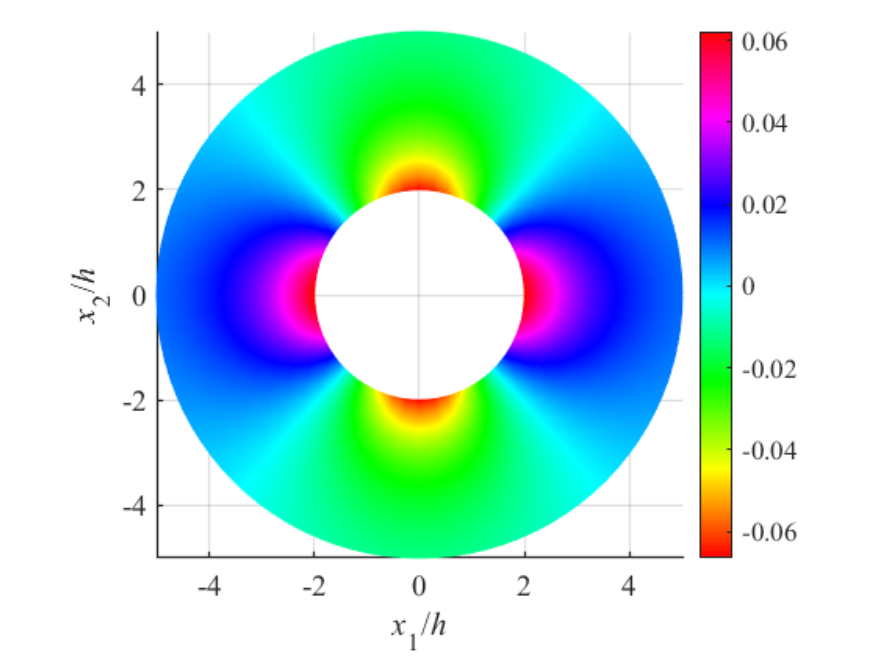}}
    \end{minipage}}
    \caption{aa}
    \label{fig3.3}
\end{figure}

\begin{figure}[h]
\centering
\subfigure[.]{\includegraphics[scale=0.36]{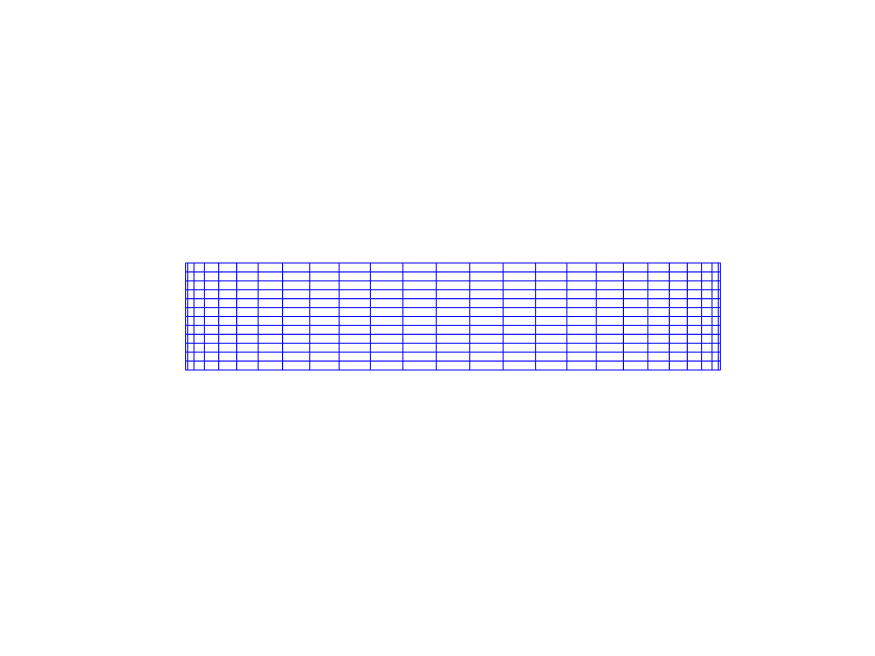}}
\subfigure[.]{\includegraphics[scale=0.36]{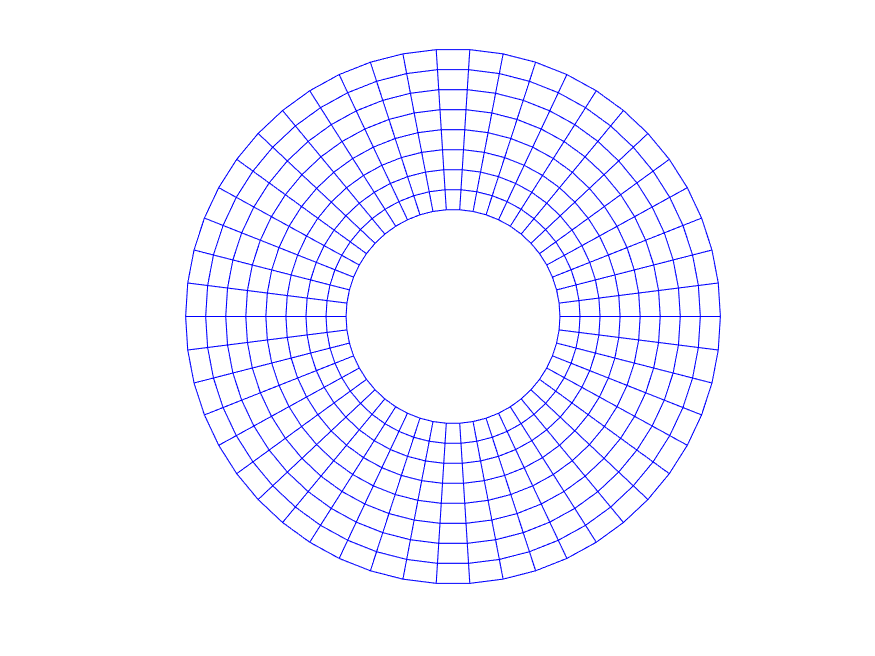}}
\subfigure[.]{\includegraphics[scale=0.36]{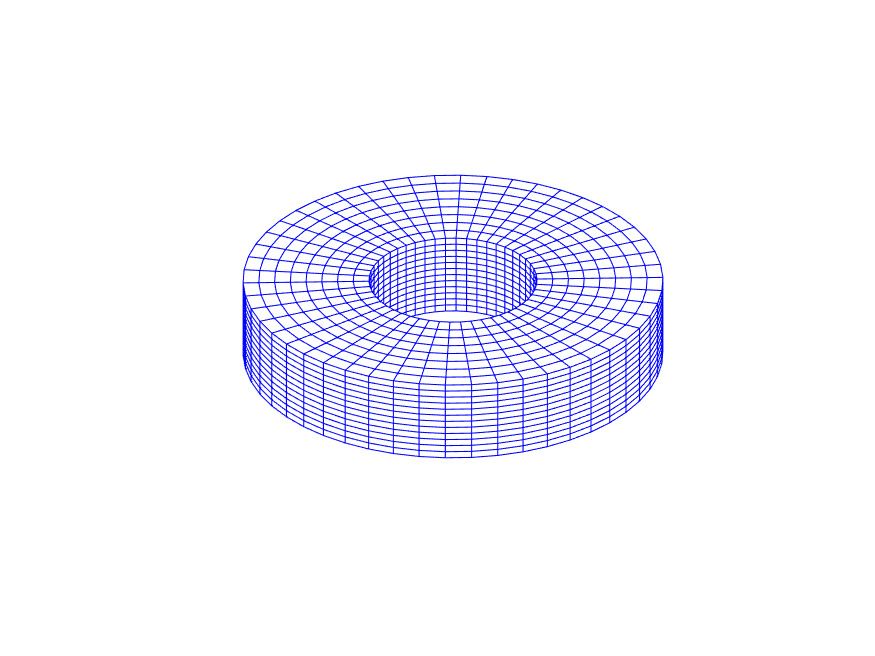}}
\subfigure[.]{\includegraphics[scale=0.36]{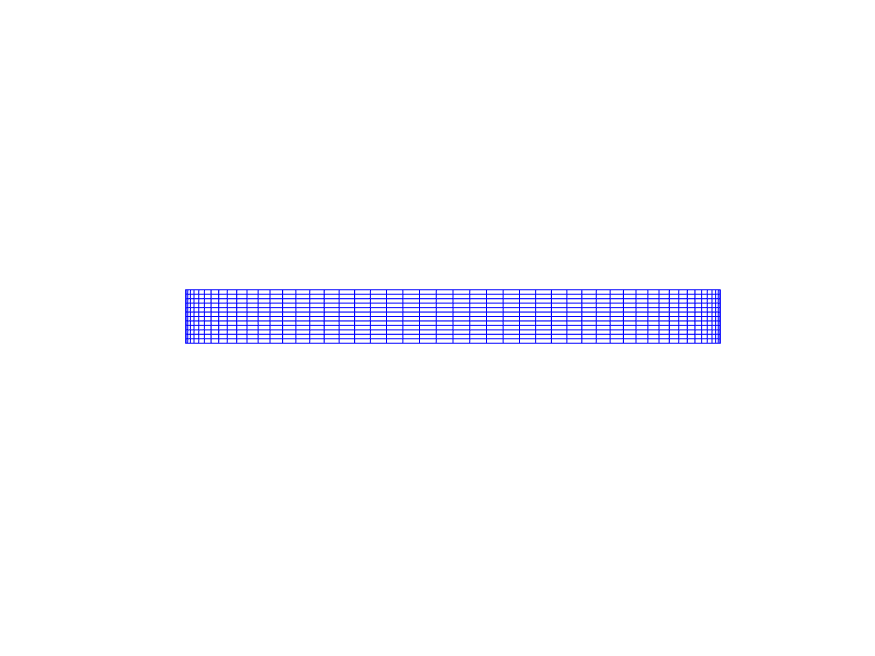}}
\subfigure[.]{\includegraphics[scale=0.36]{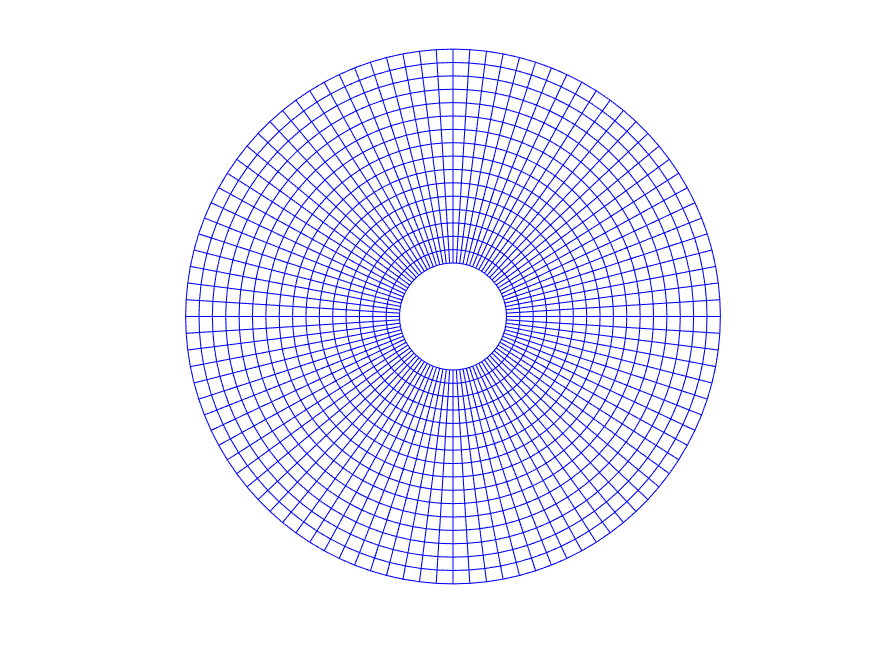}}
\subfigure[.]{\includegraphics[scale=0.36]{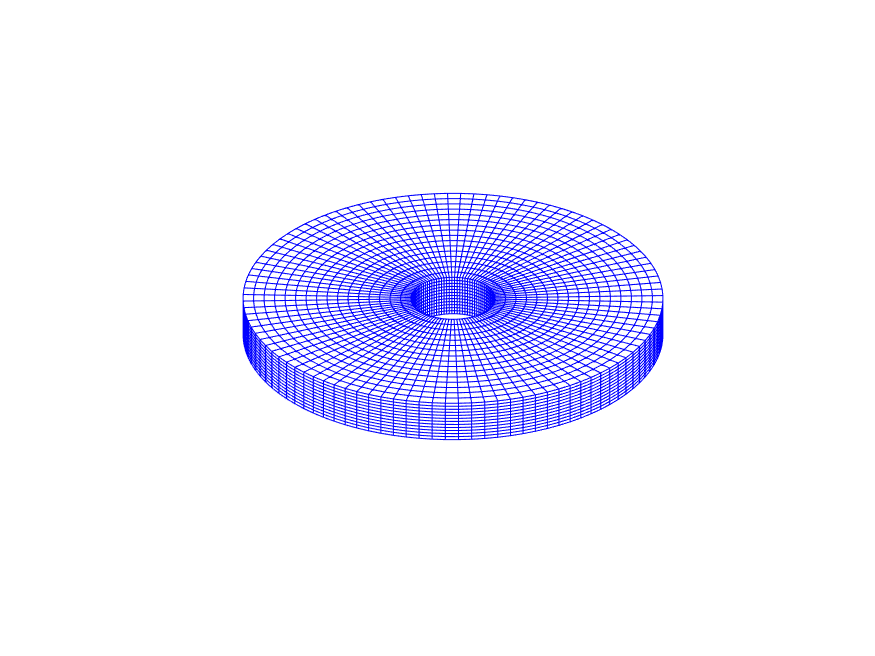}}
\subfigure[.]{\includegraphics[scale=0.36]{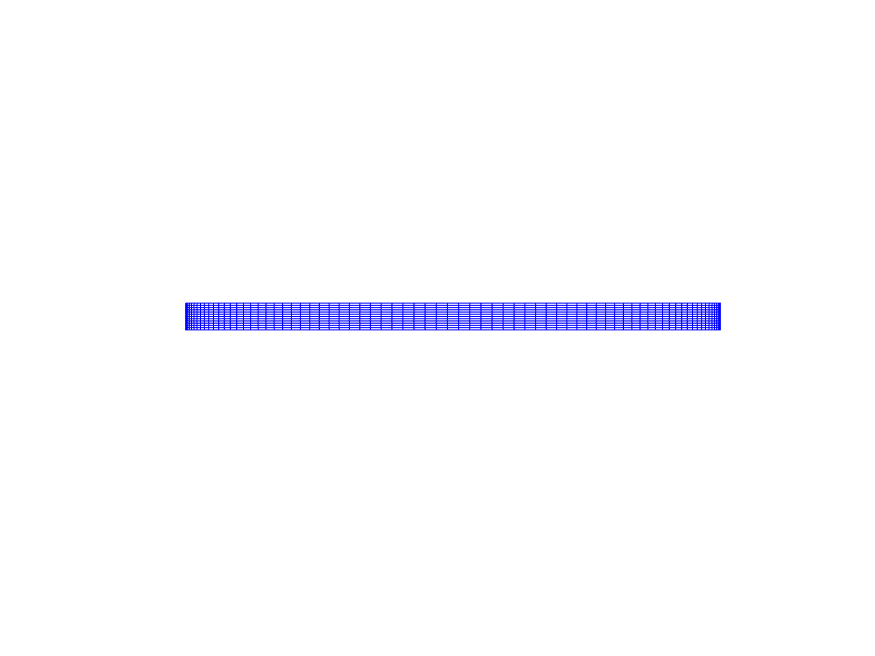}}
\subfigure[.]{\includegraphics[scale=0.36]{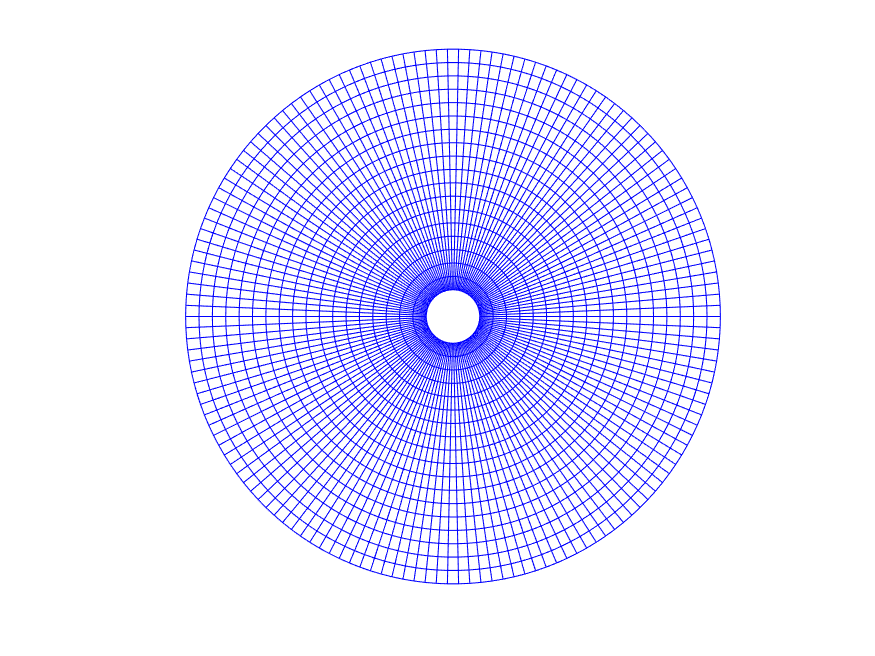}}
\subfigure[.]{\includegraphics[scale=0.36]{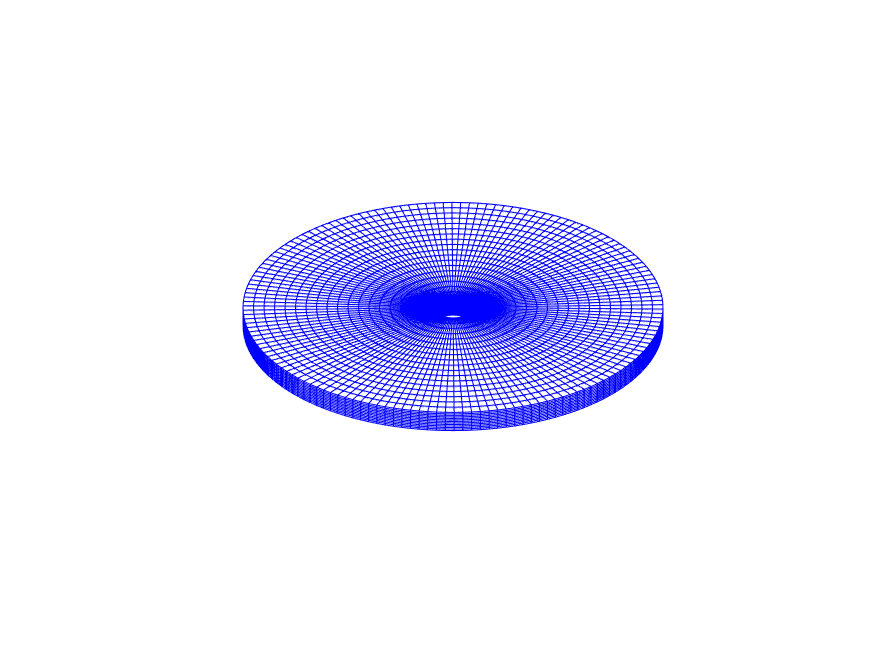}}
\caption{Error analysis of energy balance}
\label{fig4}
\end{figure}

\begin{figure}[h]
\centering
 \subfigure[.]{\includegraphics[scale=0.56]{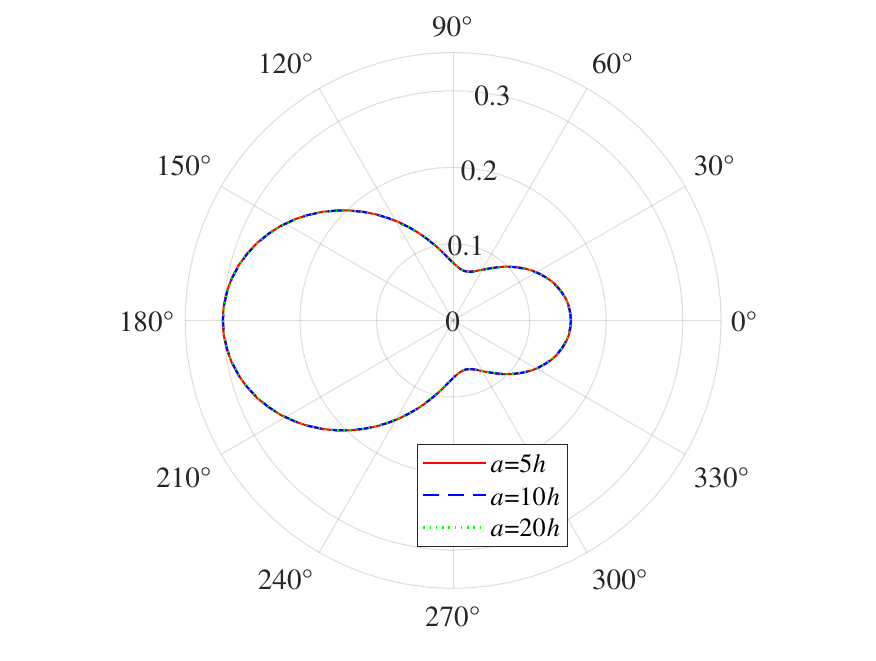}}
\subfigure[.]{\includegraphics[scale=0.56]{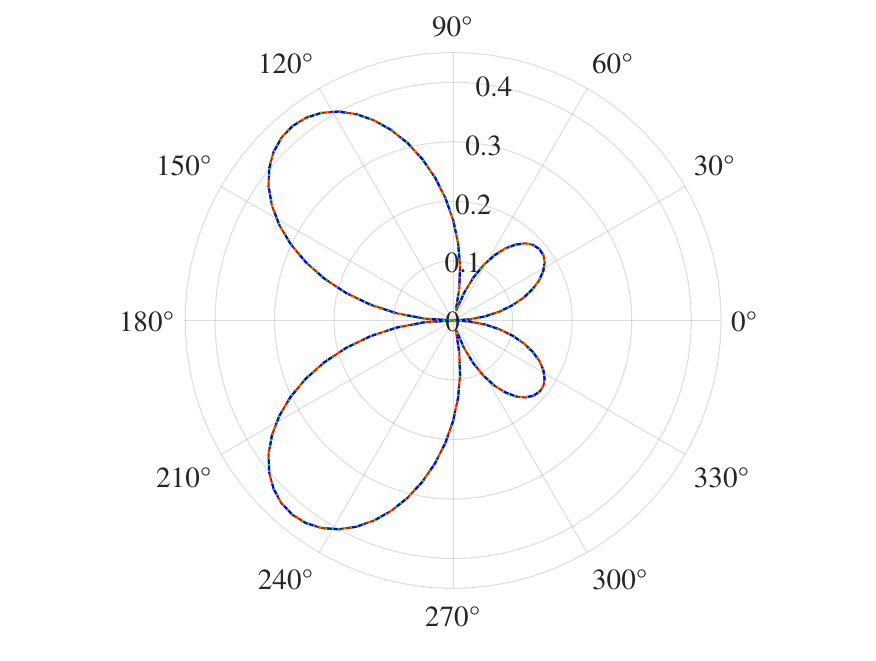}}
\caption{Error analysis of energy balance}
\label{fig5}
\end{figure}

\begin{figure}[h]
\centering
 \subfigure[.]{\includegraphics[scale=0.36]{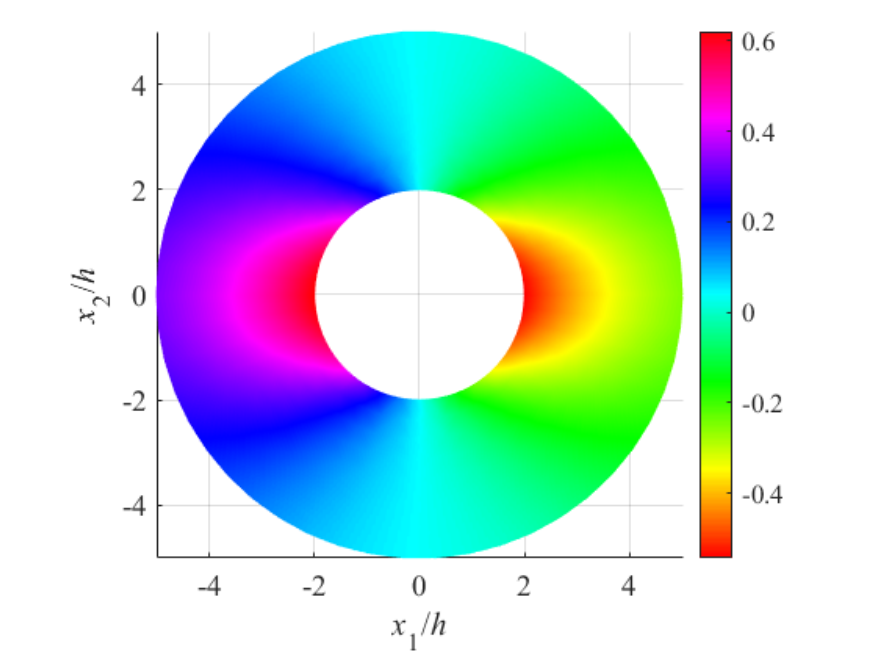}}
\subfigure[.]{\includegraphics[scale=0.36]{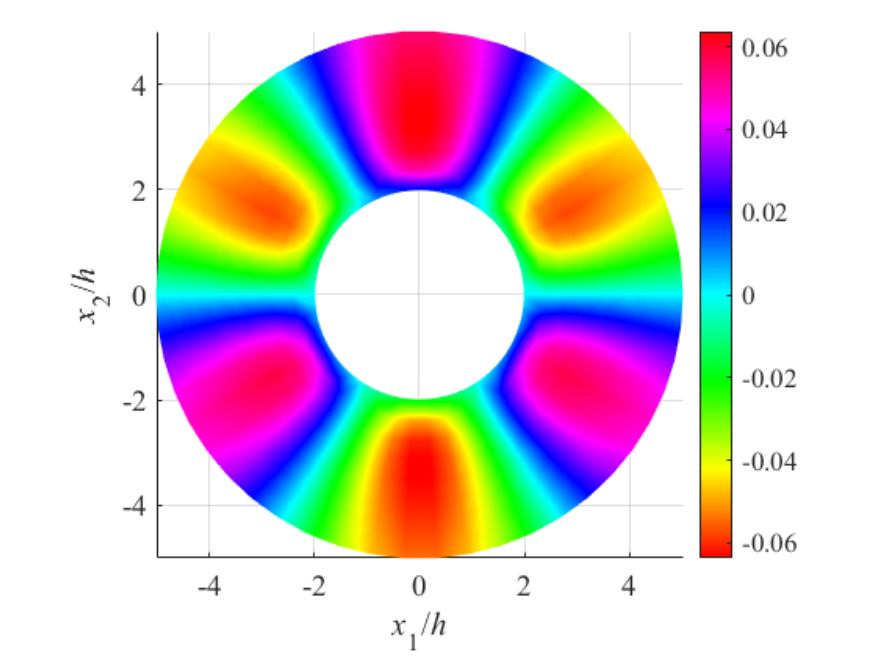}}
\subfigure[.]{\includegraphics[scale=0.36]{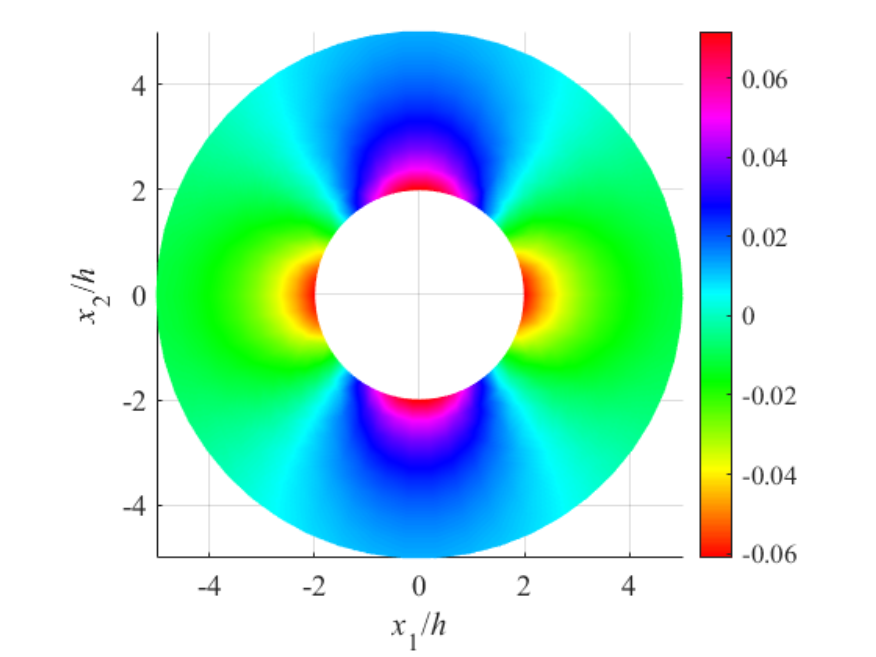}}
\subfigure[.]{\includegraphics[scale=0.36]{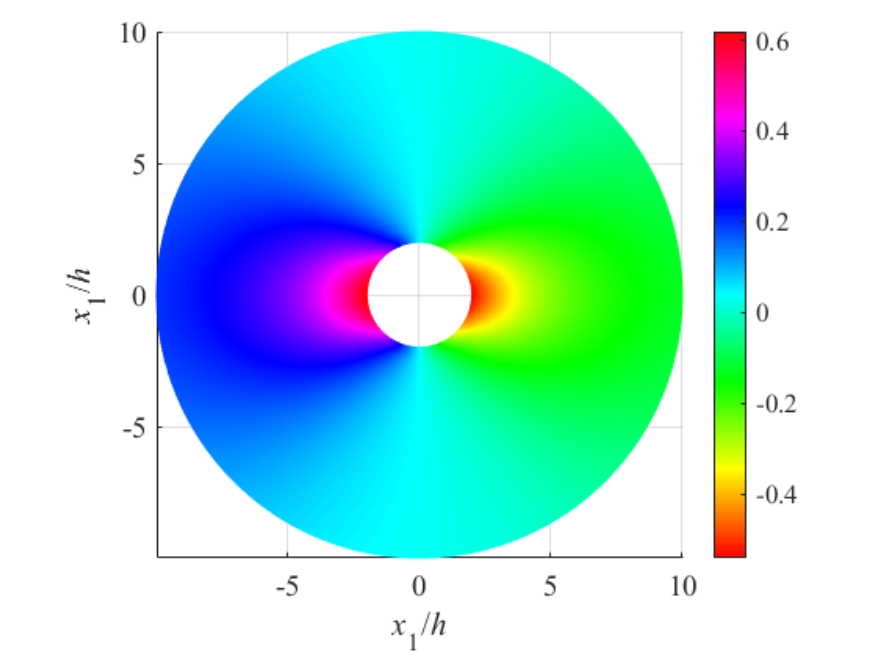}}
 \subfigure[.]{\includegraphics[scale=0.36]{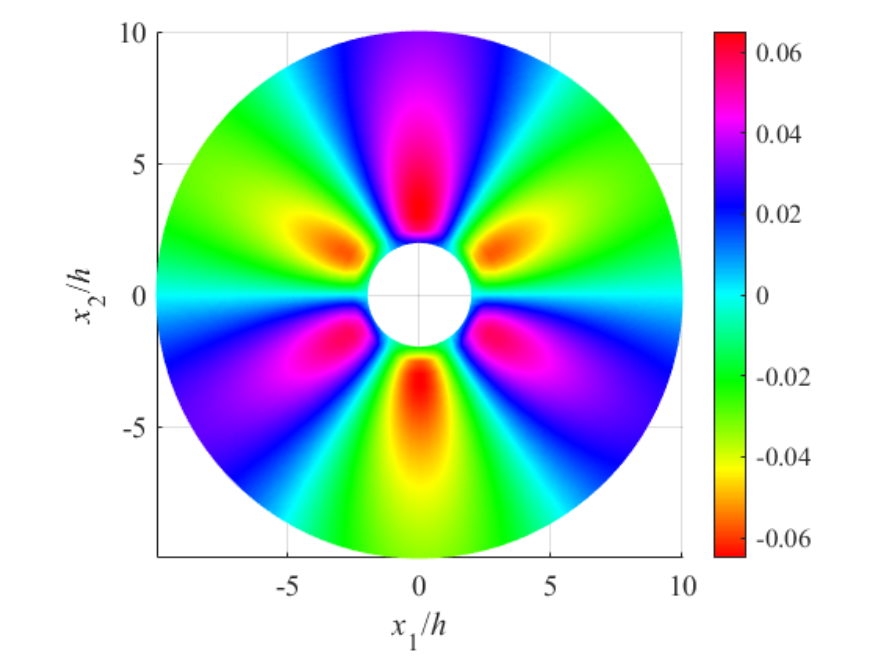}}
\subfigure[.]{\includegraphics[scale=0.36]{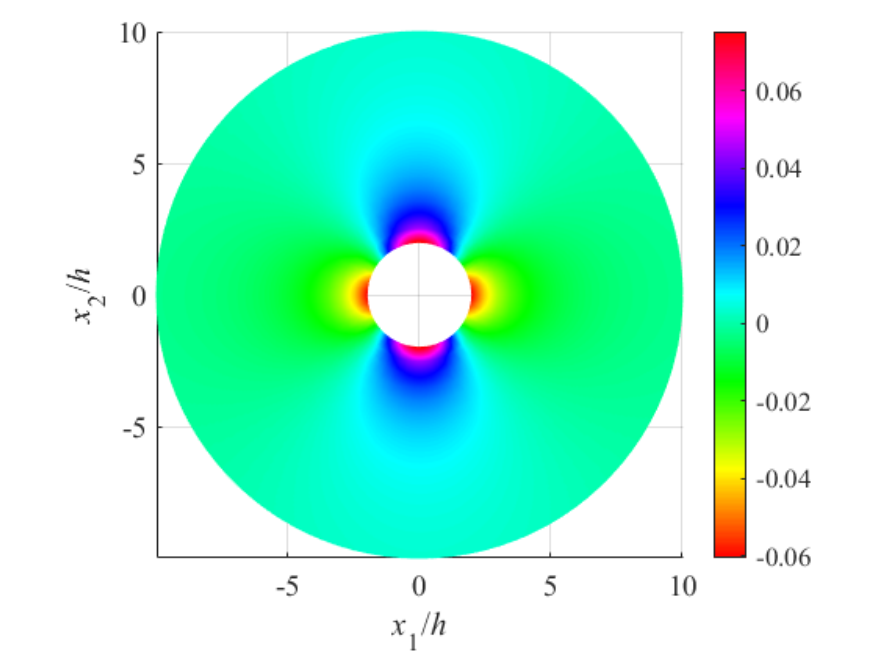}}
\subfigure[.]{\includegraphics[scale=0.36]{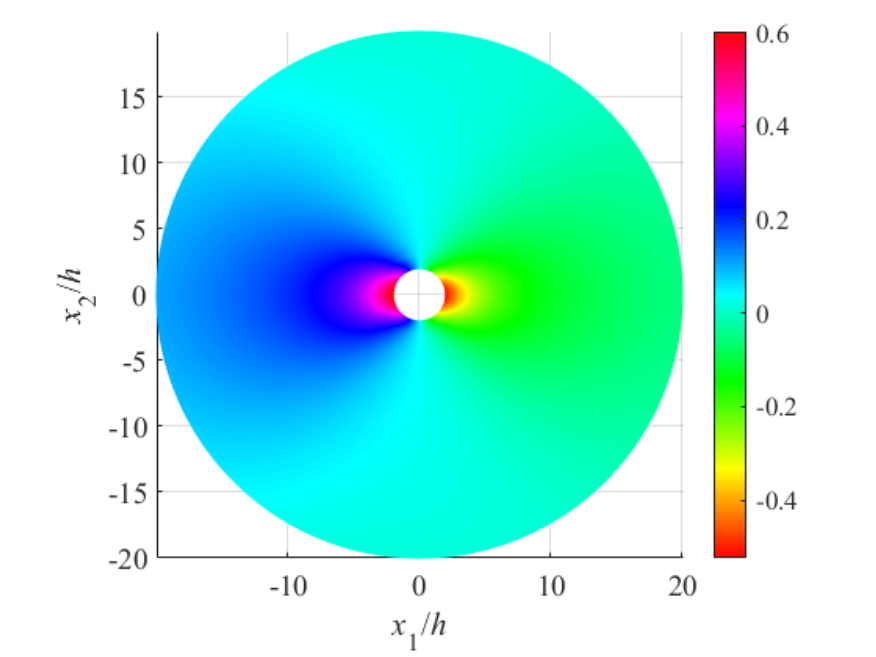}}
\subfigure[.]{\includegraphics[scale=0.36]{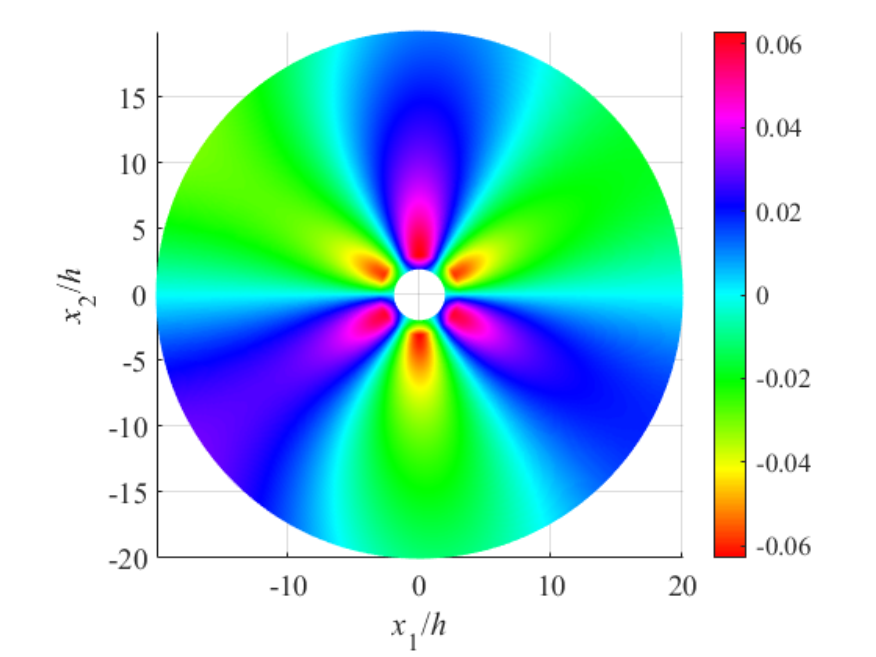}}
 \subfigure[.]{\includegraphics[scale=0.36]{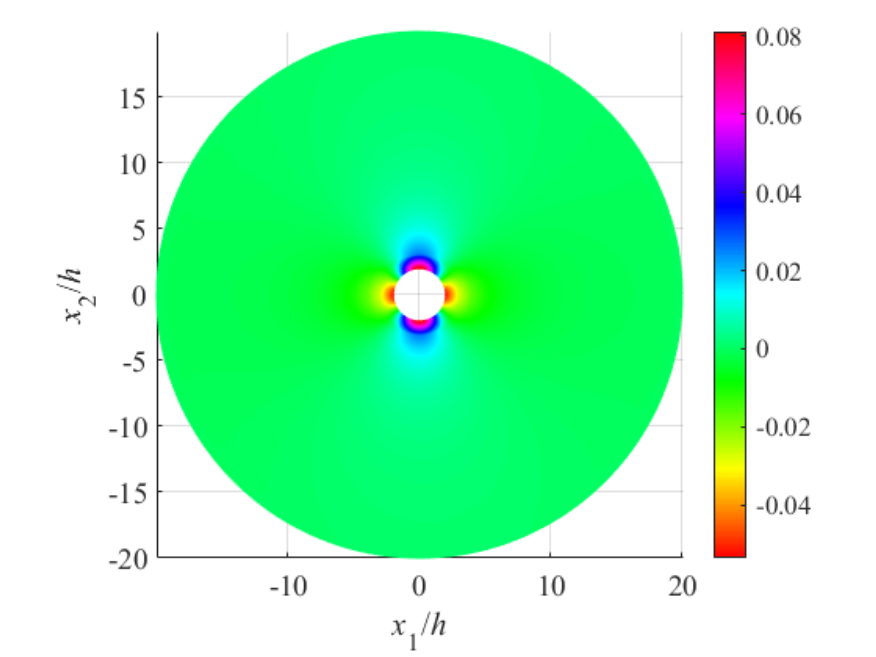}}
\caption{Error analysis of energy balance}
\label{fig6}
\end{figure}

\begin{figure}[h]
\centering
\subfigure[.]{\includegraphics[scale=0.36]{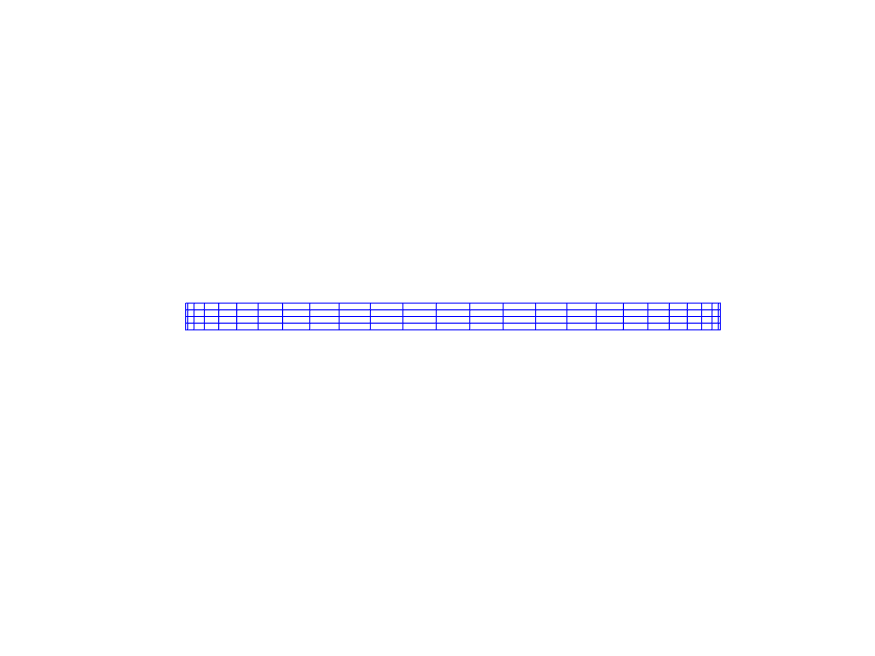}}
\subfigure[.]{\includegraphics[scale=0.36]{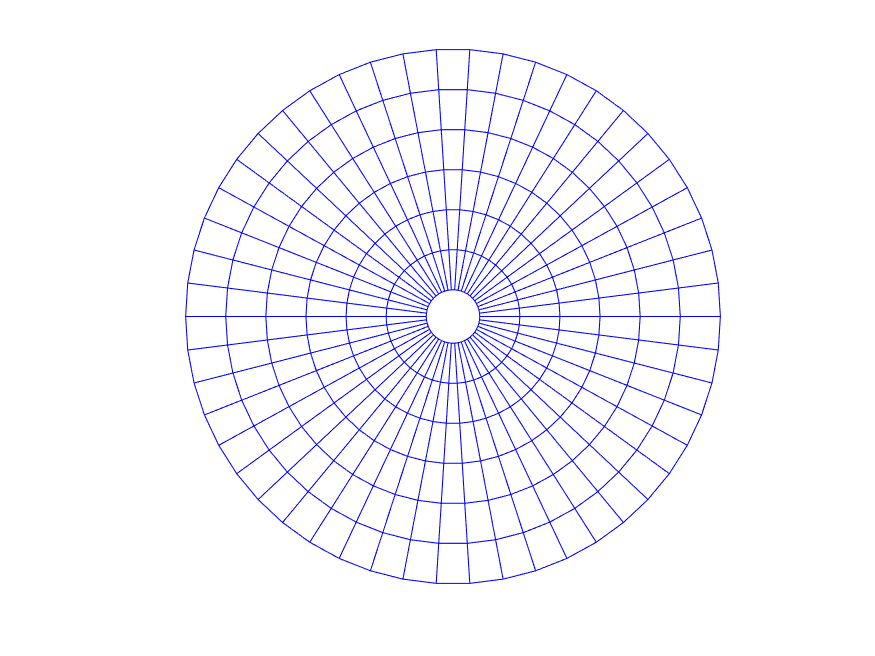}}
\subfigure[.]{\includegraphics[scale=0.36]{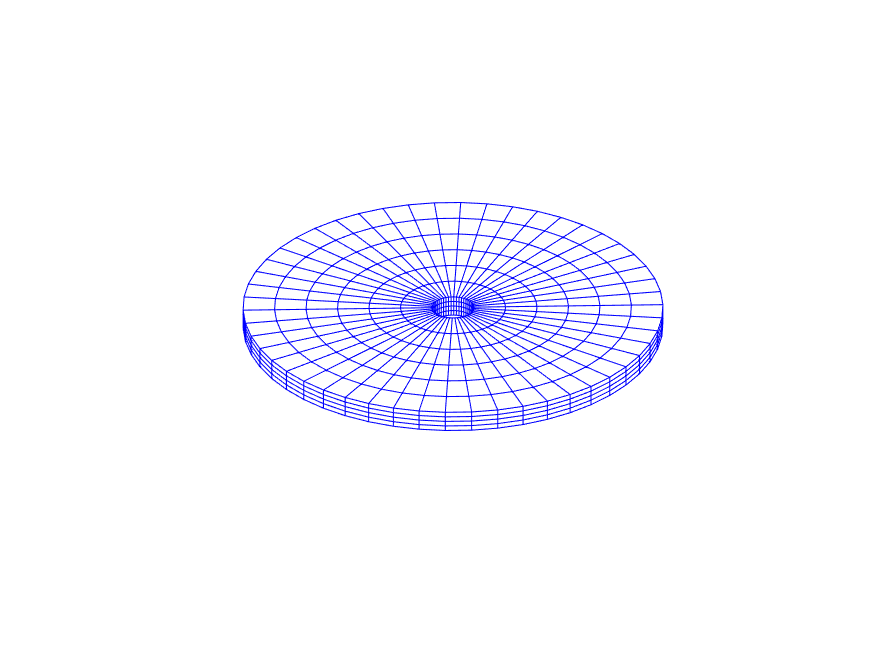}}
\subfigure[.]{\includegraphics[scale=0.36]{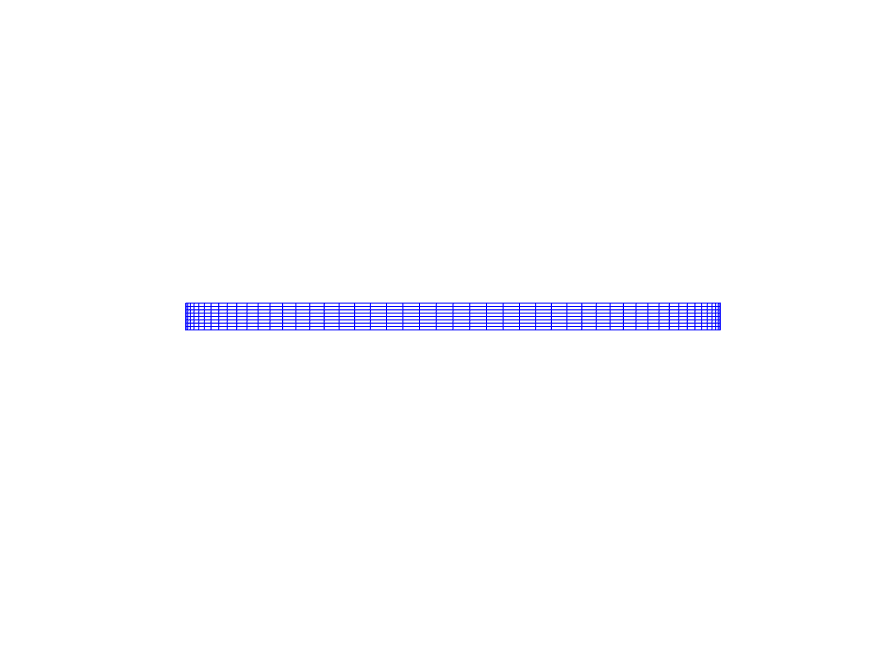}}
\subfigure[.]{\includegraphics[scale=0.36]{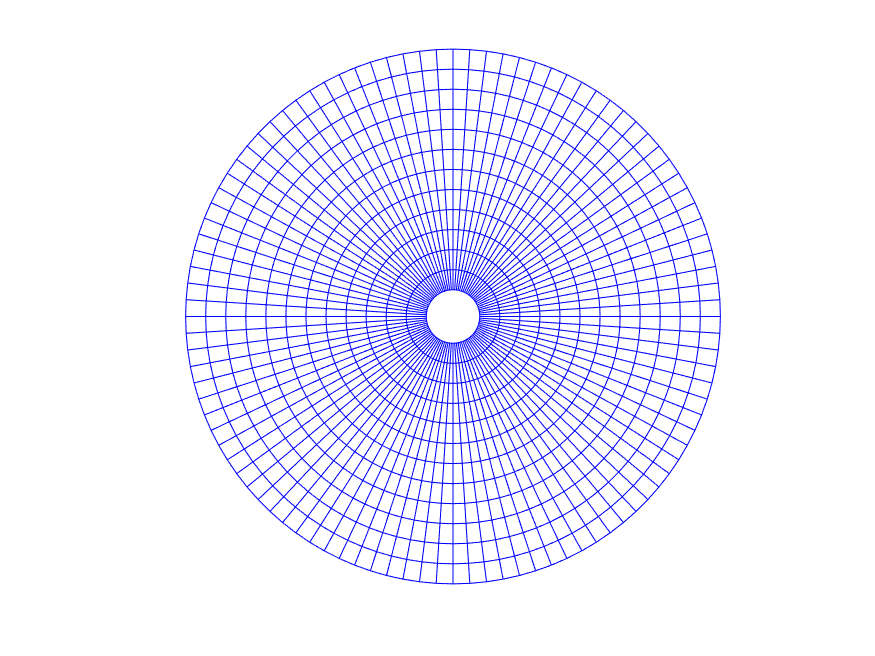}}
\subfigure[.]{\includegraphics[scale=0.36]{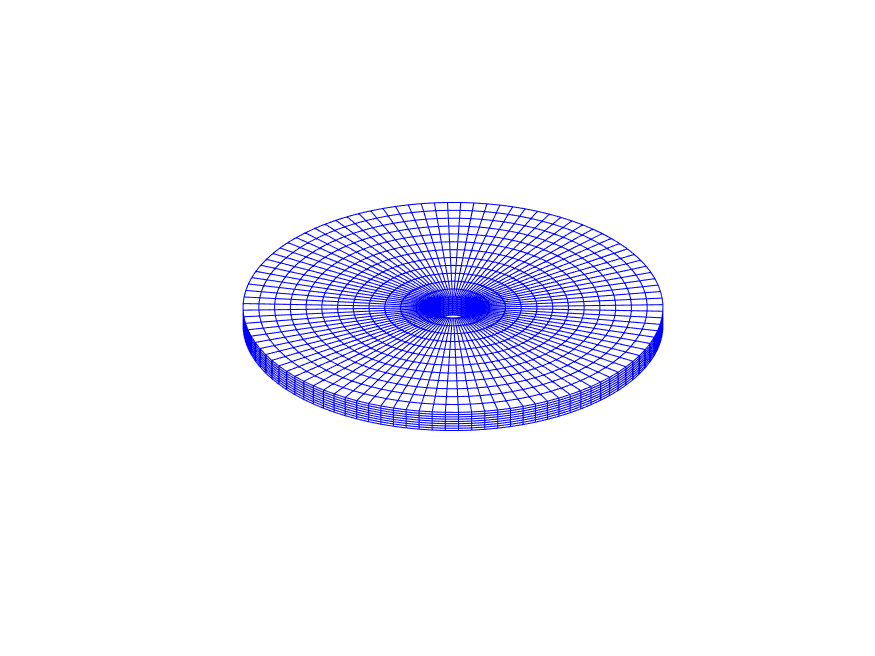}}
\subfigure[.]{\includegraphics[scale=0.36]{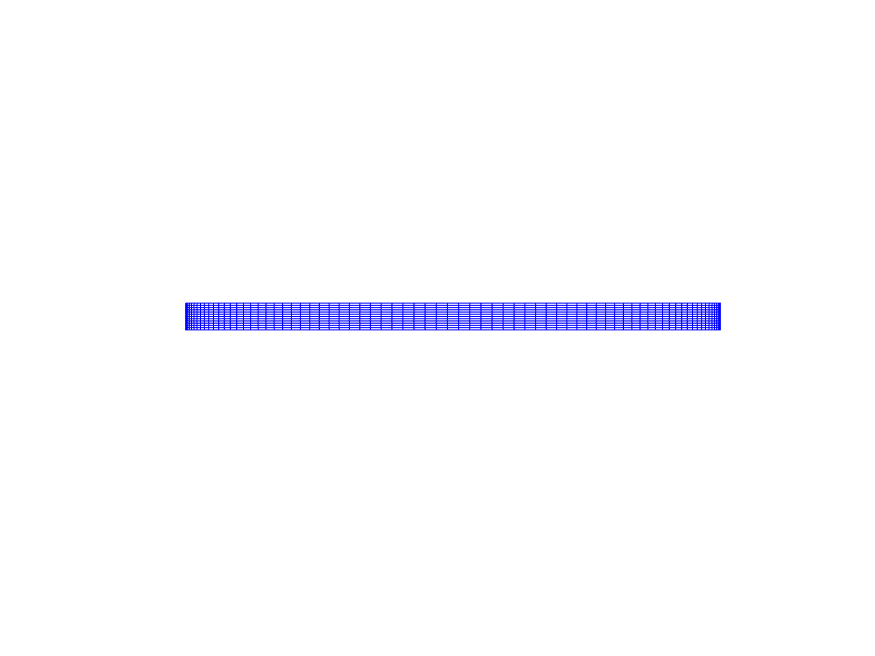}}
\subfigure[.]{\includegraphics[scale=0.36]{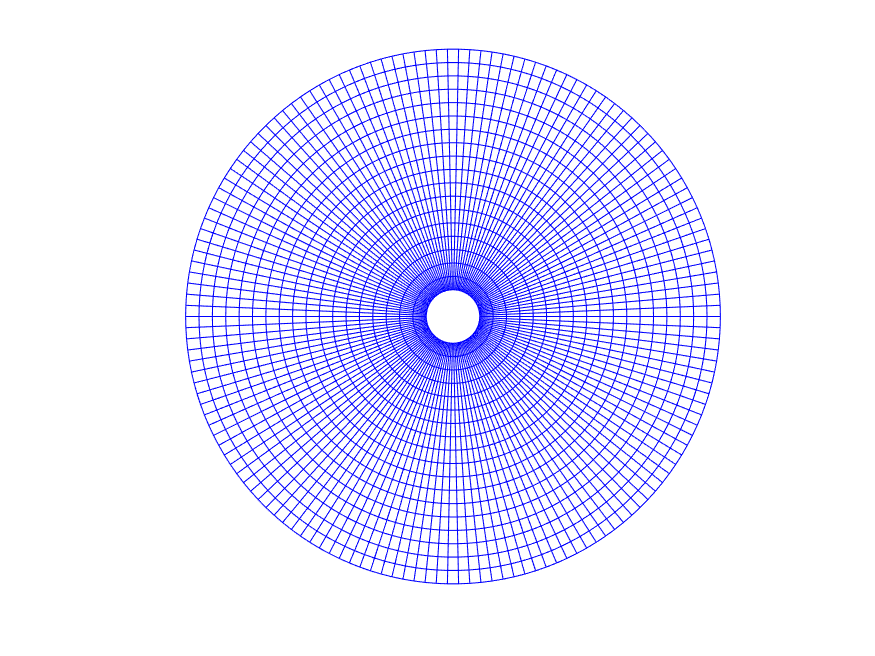}}
\subfigure[.]{\includegraphics[scale=0.36]{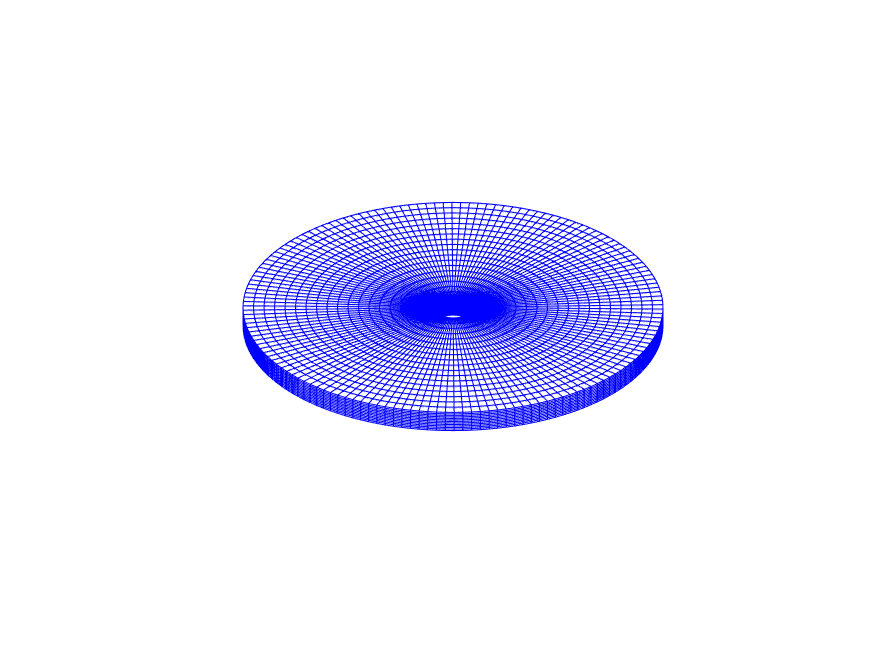}}
\caption{Error analysis of energy balance}
\label{fig7}
\end{figure}

\begin{figure}[h]
\centering
 \subfigure[.]{\includegraphics[scale=0.56]{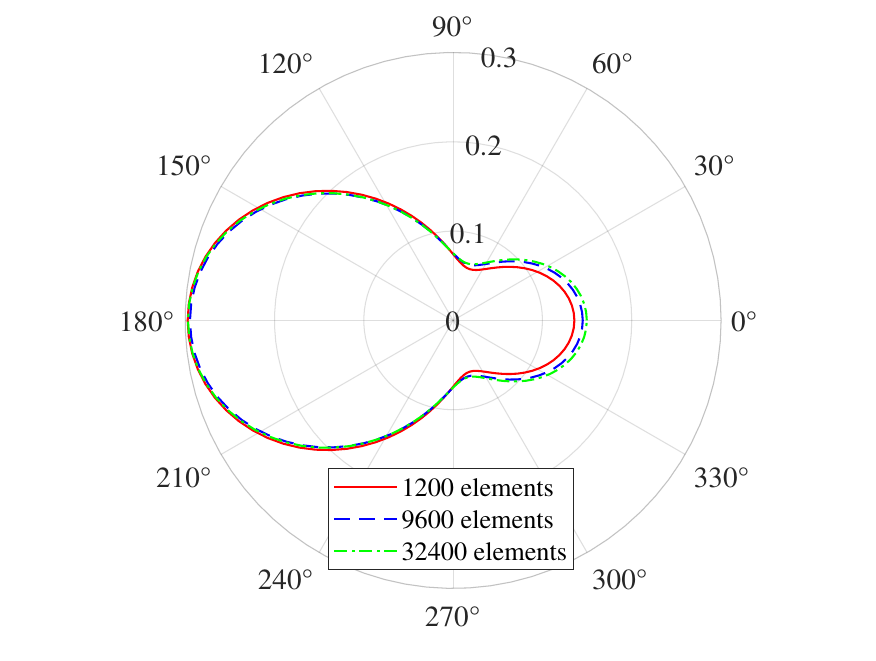}}
\subfigure[.]{\includegraphics[scale=0.56]{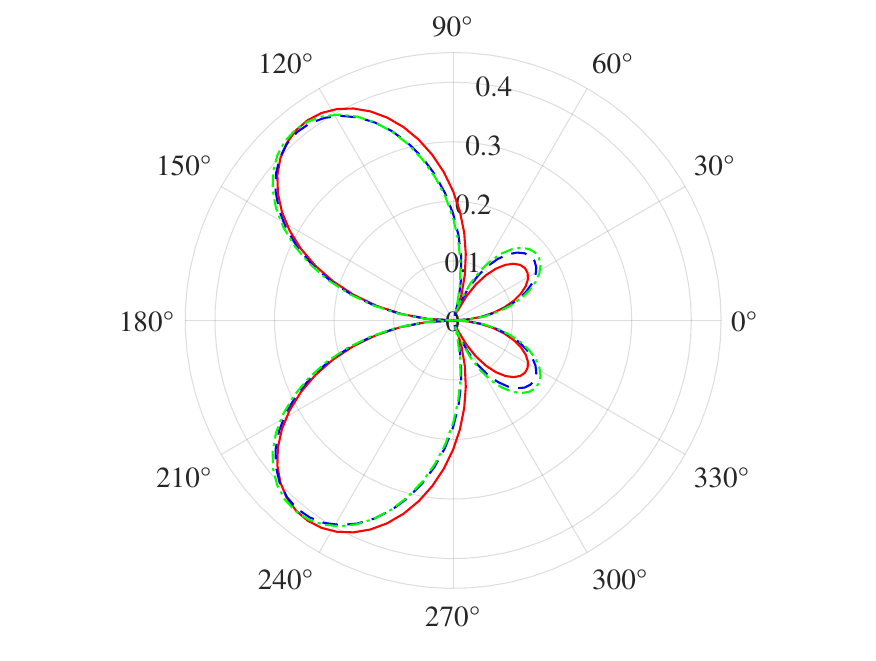}}
\caption{Error analysis of energy balance}
\label{fig8}
\end{figure}

\begin{figure}[h]
\centering
 \subfigure[.]{\includegraphics[scale=0.36]{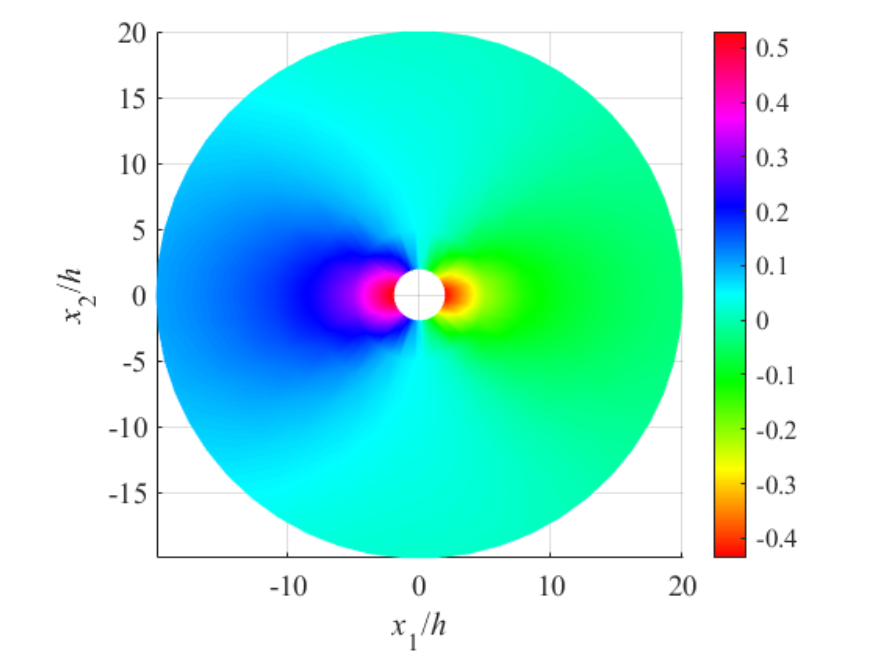}}
\subfigure[.]{\includegraphics[scale=0.36]{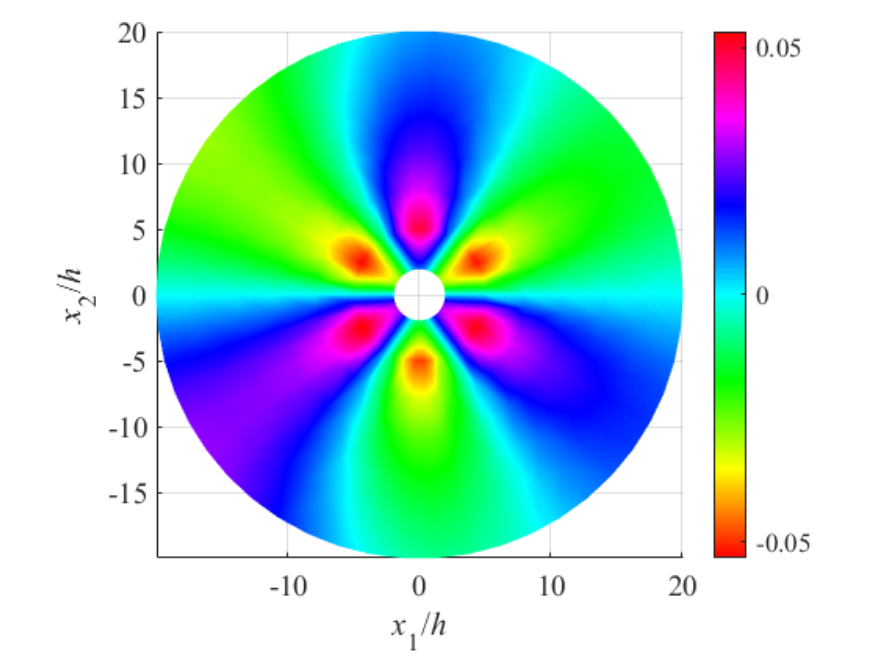}}
\subfigure[.]{\includegraphics[scale=0.36]{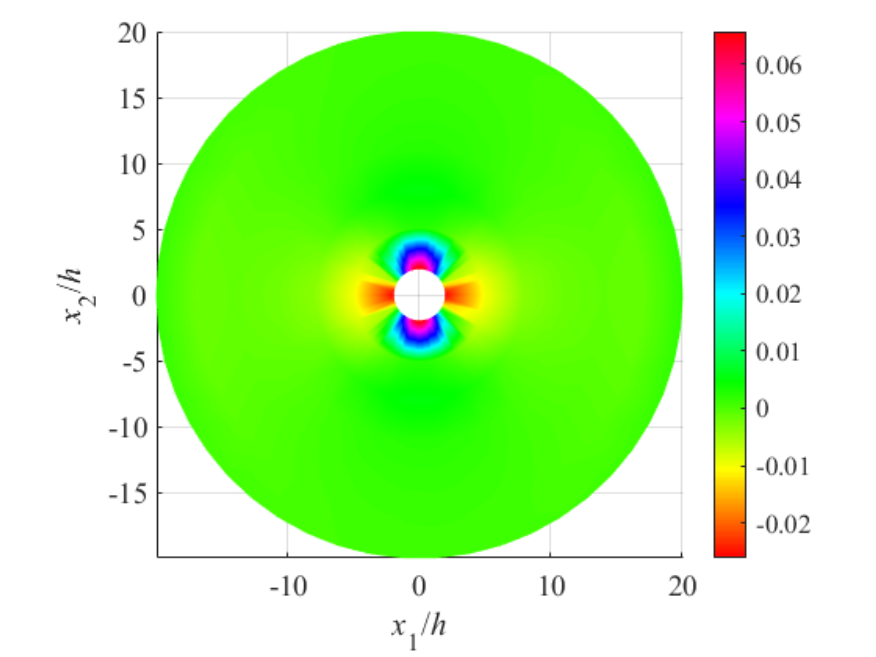}}
\subfigure[.]{\includegraphics[scale=0.36]{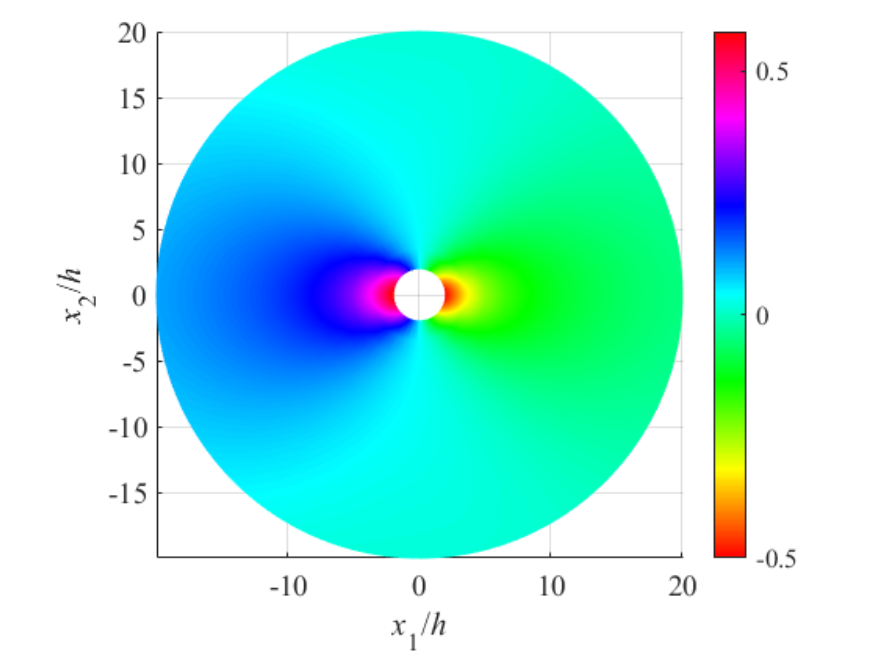}}
 \subfigure[.]{\includegraphics[scale=0.36]{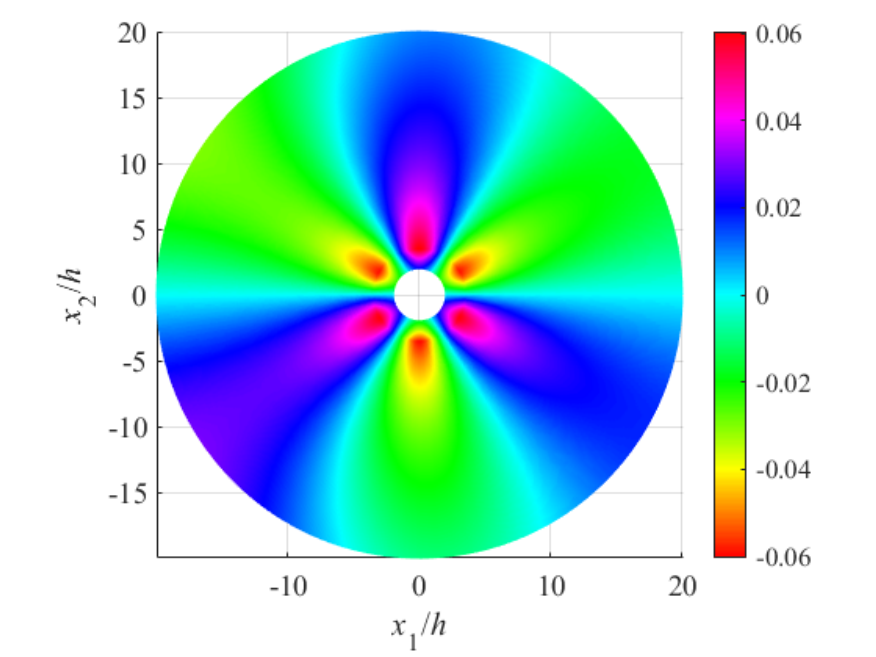}}
\subfigure[.]{\includegraphics[scale=0.36]{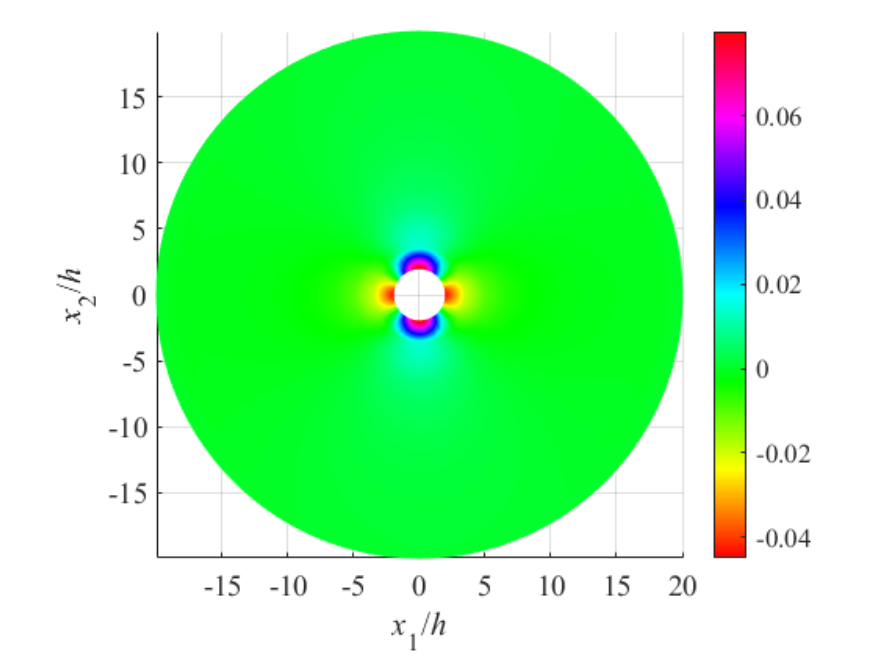}}
\subfigure[.]{\includegraphics[scale=0.36]{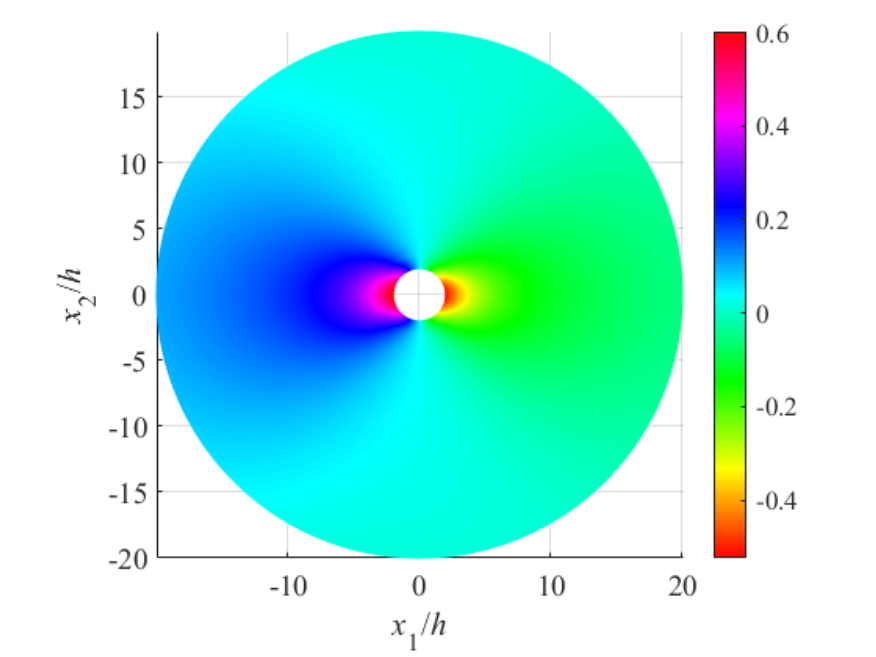}}
\subfigure[.]{\includegraphics[scale=0.36]{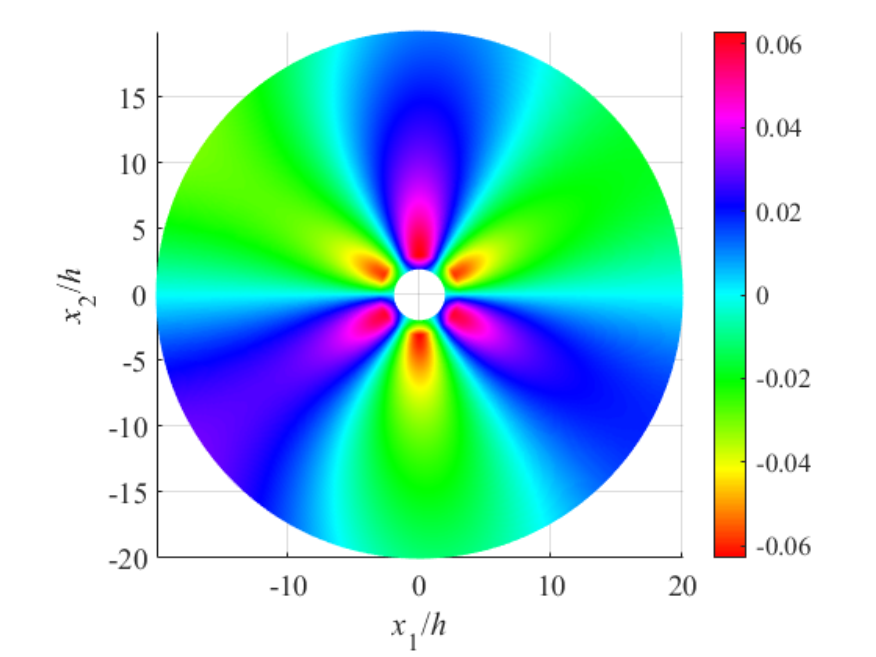}}
 \subfigure[.]{\includegraphics[scale=0.36]{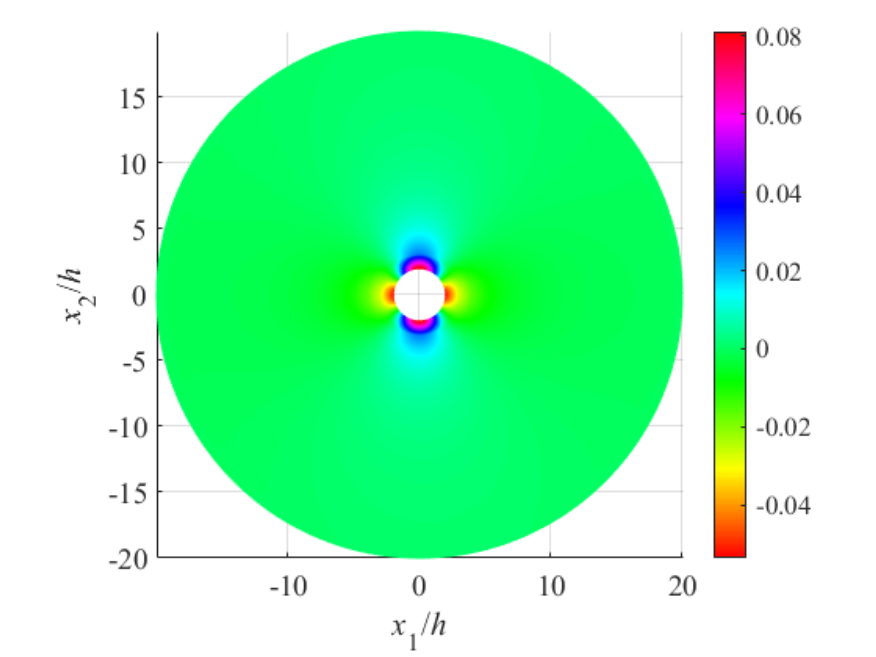}}
\caption{Error analysis of energy balance}
\label{fig9}
\end{figure}

\begin{figure}[h]
\centering
\subfigure[.]{\includegraphics[scale=0.36]{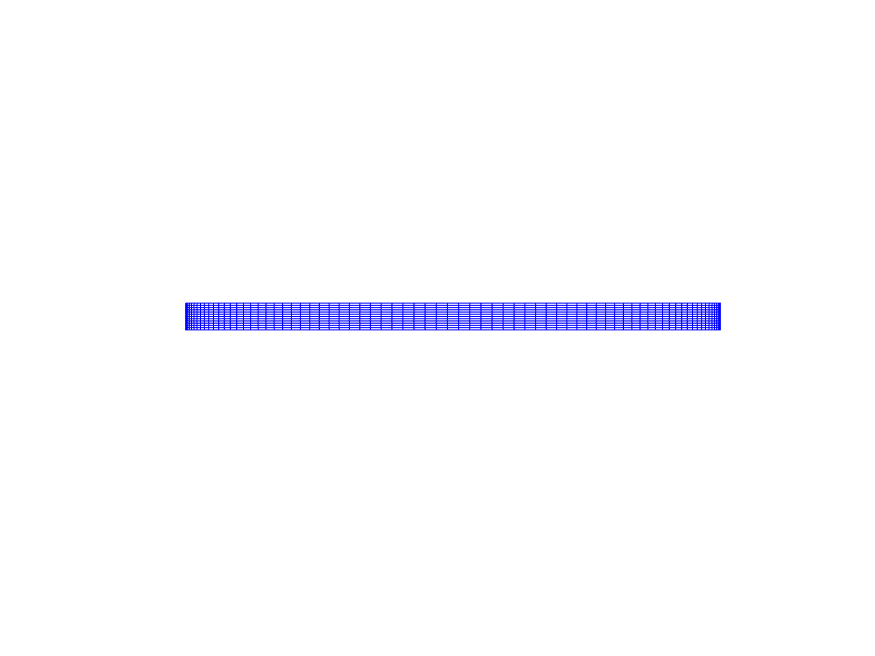}}
\subfigure[.]{\includegraphics[scale=0.36]{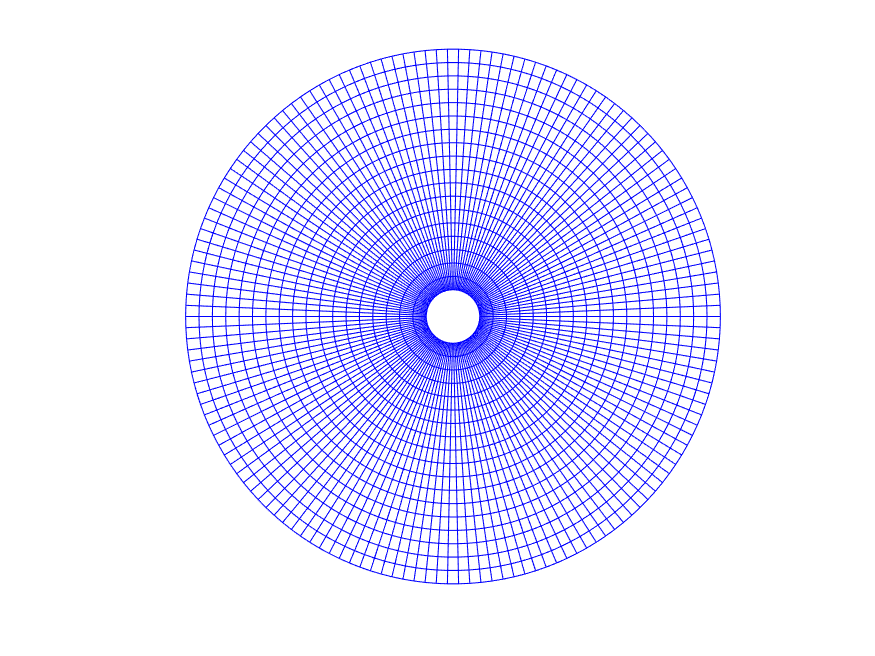}}
\subfigure[.]{\includegraphics[scale=0.36]{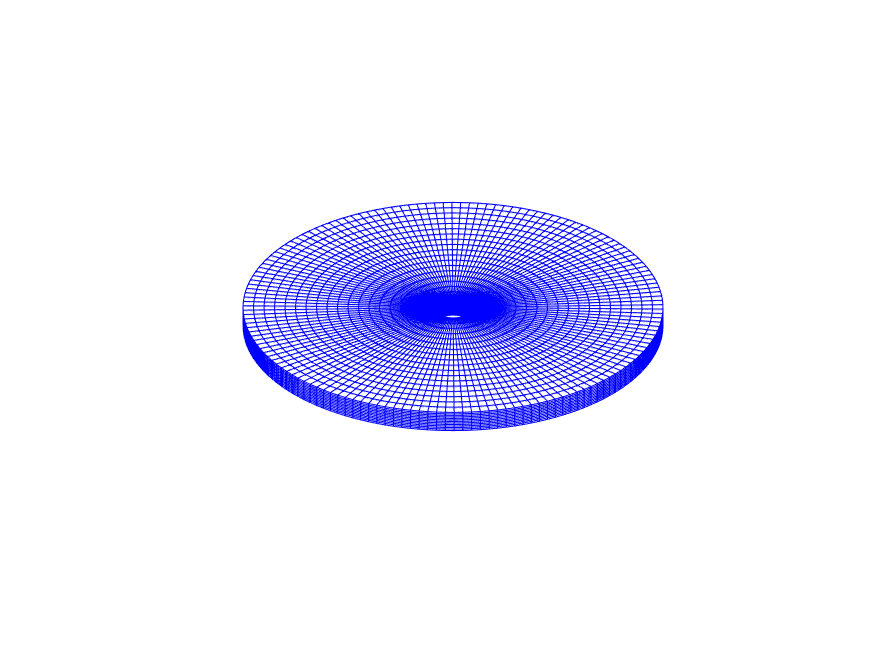}}
\subfigure[.]{\includegraphics[scale=0.36]{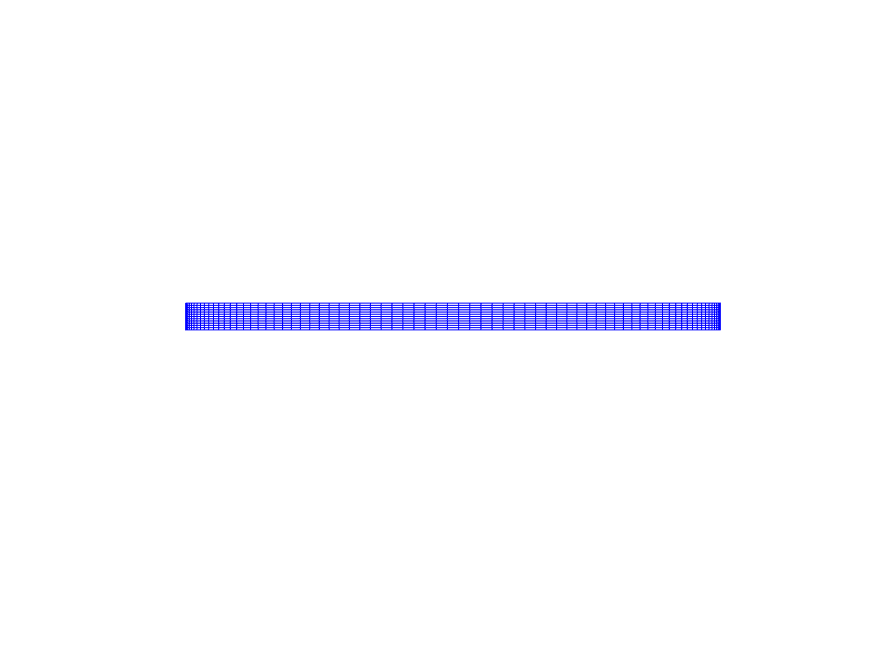}}
\subfigure[.]{\includegraphics[scale=0.36]{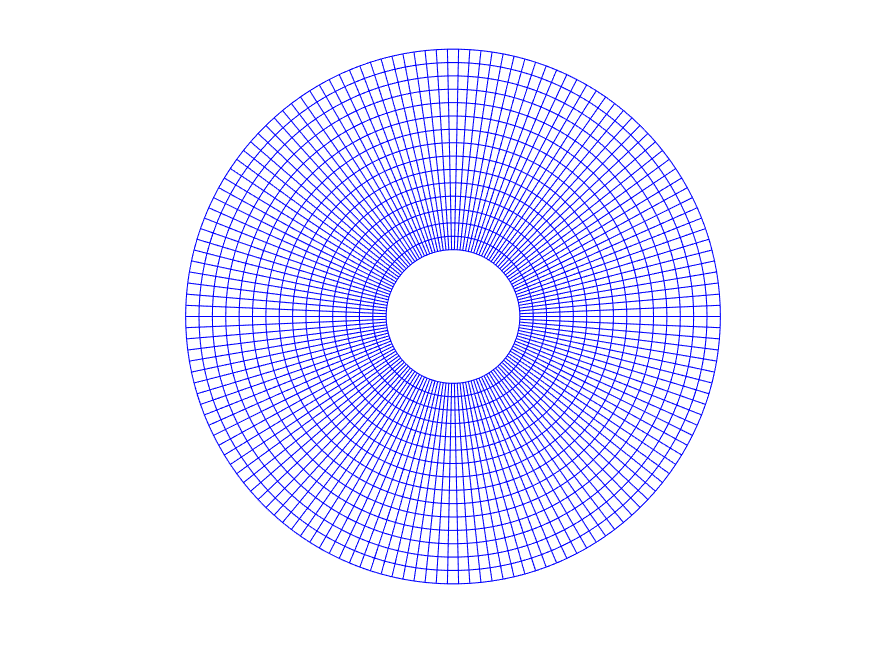}}
\subfigure[.]{\includegraphics[scale=0.36]{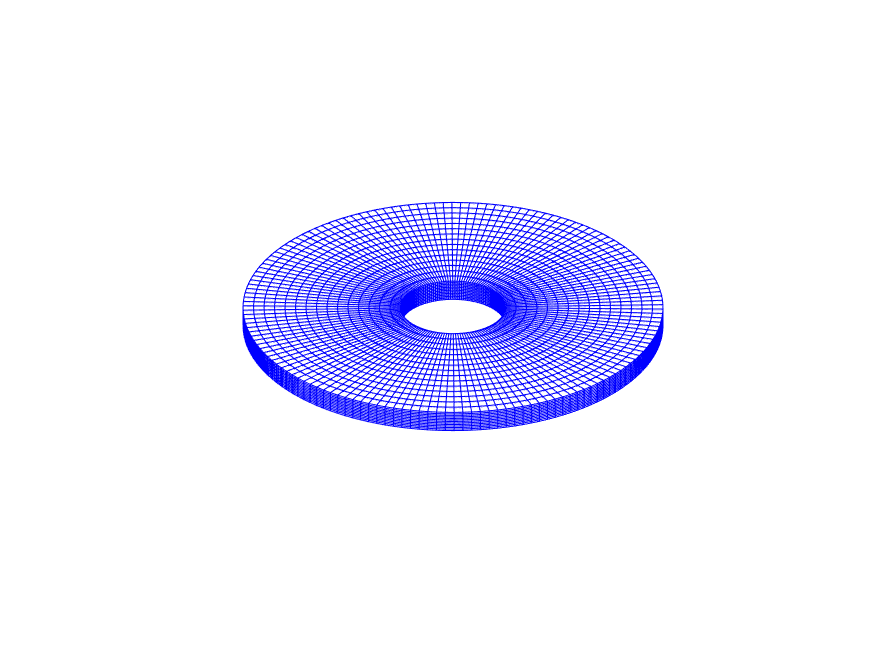}}
\subfigure[.]{\includegraphics[scale=0.36]{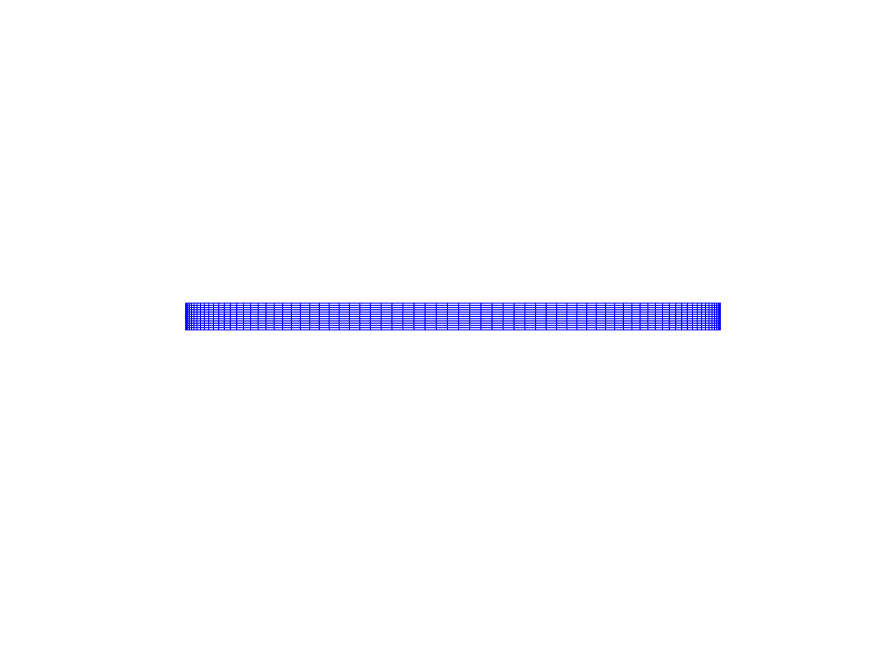}}
\subfigure[.]{\includegraphics[scale=0.36]{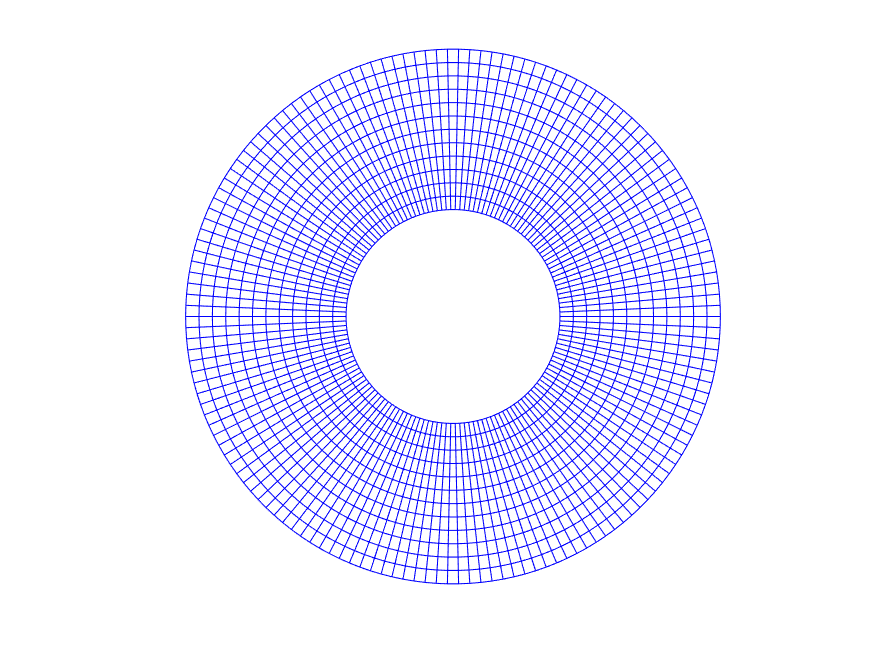}}
\subfigure[.]{\includegraphics[scale=0.36]{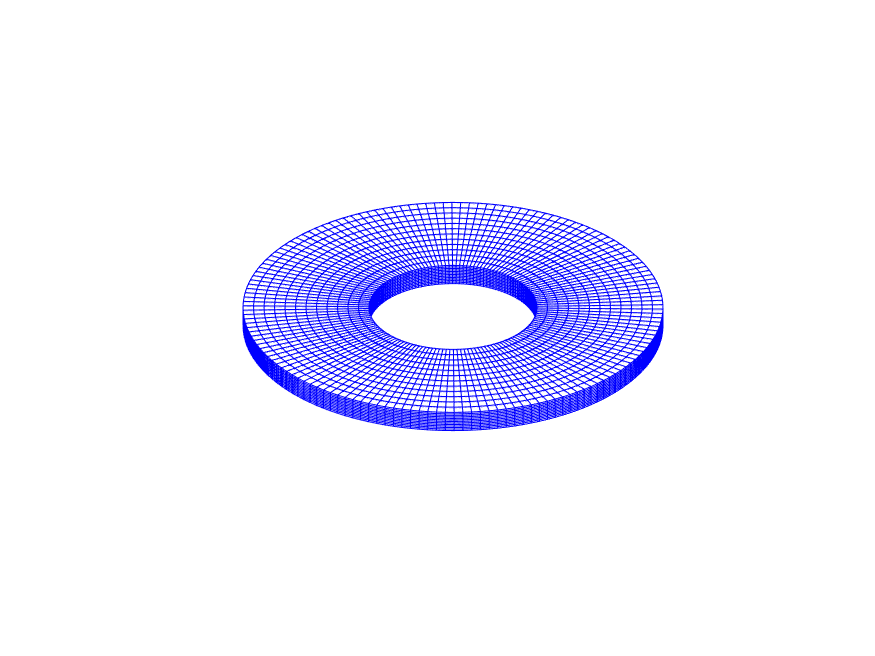}}
\caption{Error analysis of energy balance}
\label{fig10}
\end{figure}

\begin{figure}[h]
\centering
 \subfigure[.]{\includegraphics[scale=0.56]{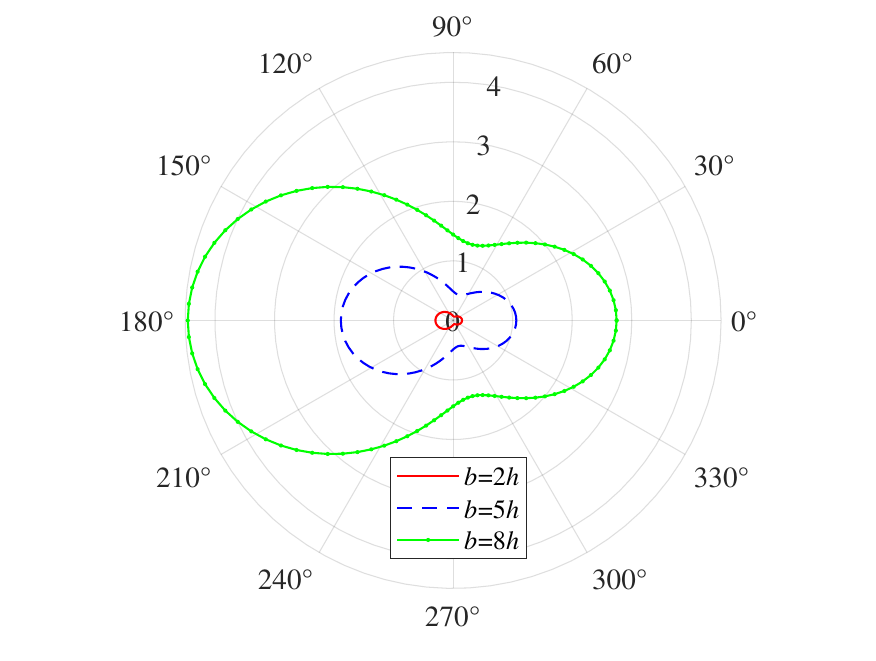}}
\subfigure[.]{\includegraphics[scale=0.56]{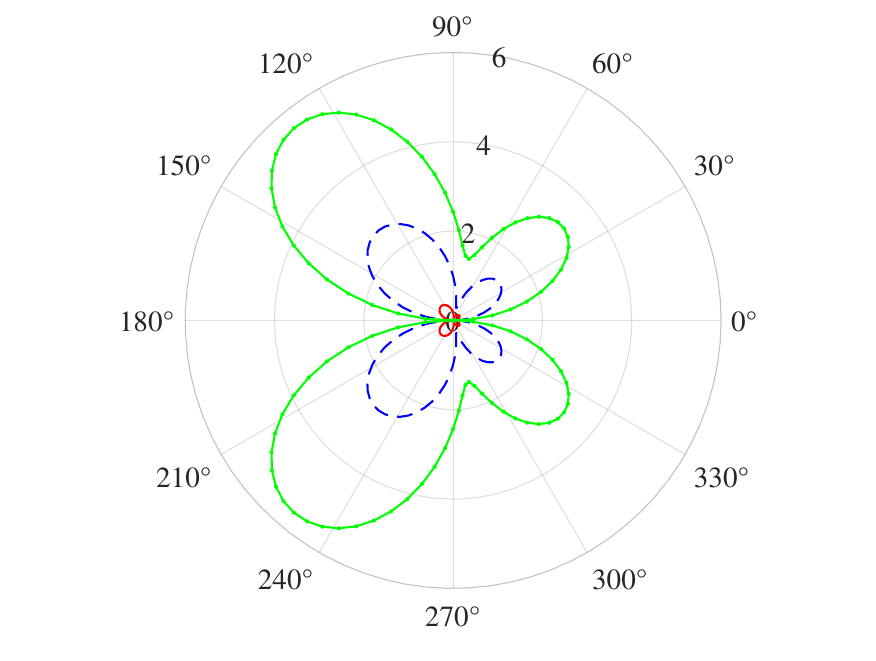}}
\caption{Error analysis of energy balance}
\label{fig11}
\end{figure}

\begin{figure}[h]
\centering
 \subfigure[.]{\includegraphics[scale=0.36]{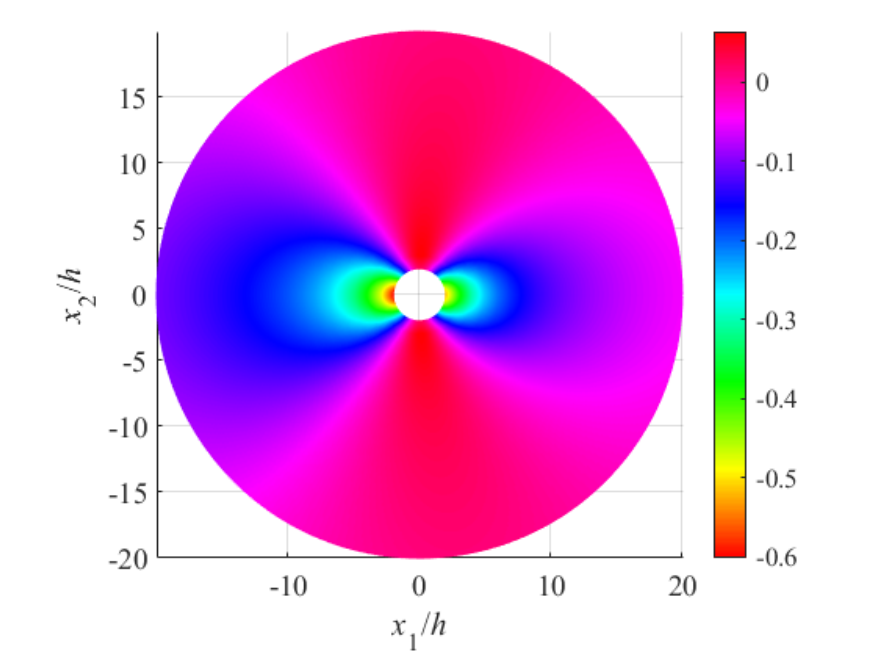}}
\subfigure[.]{\includegraphics[scale=0.36]{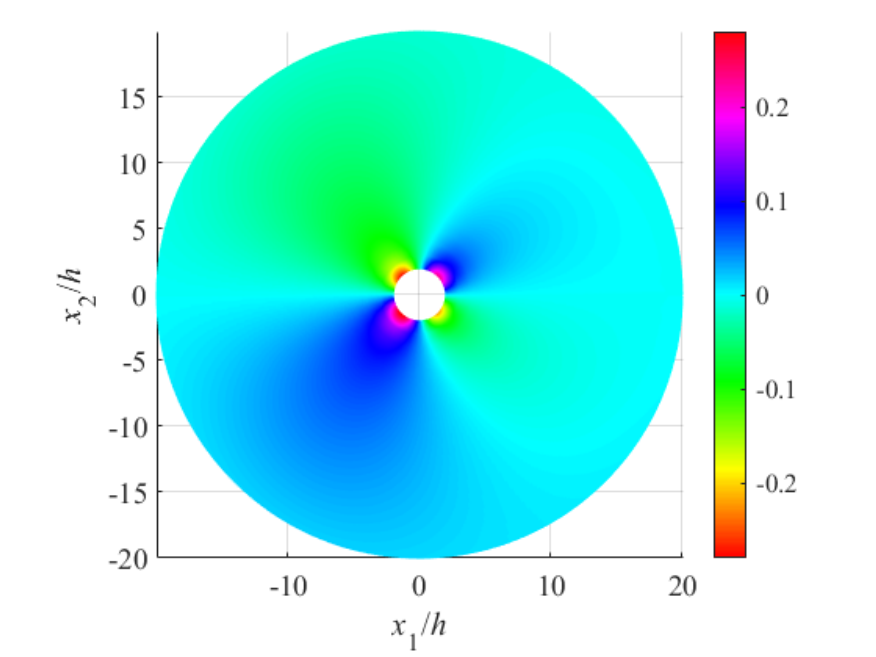}}
\subfigure[.]{\includegraphics[scale=0.36]{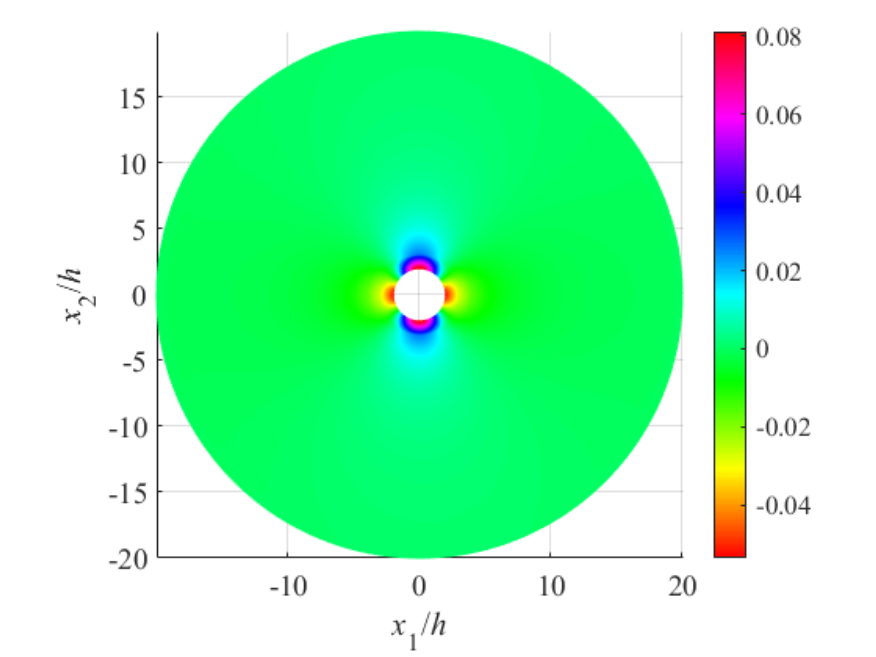}}
\subfigure[.]{\includegraphics[scale=0.36]{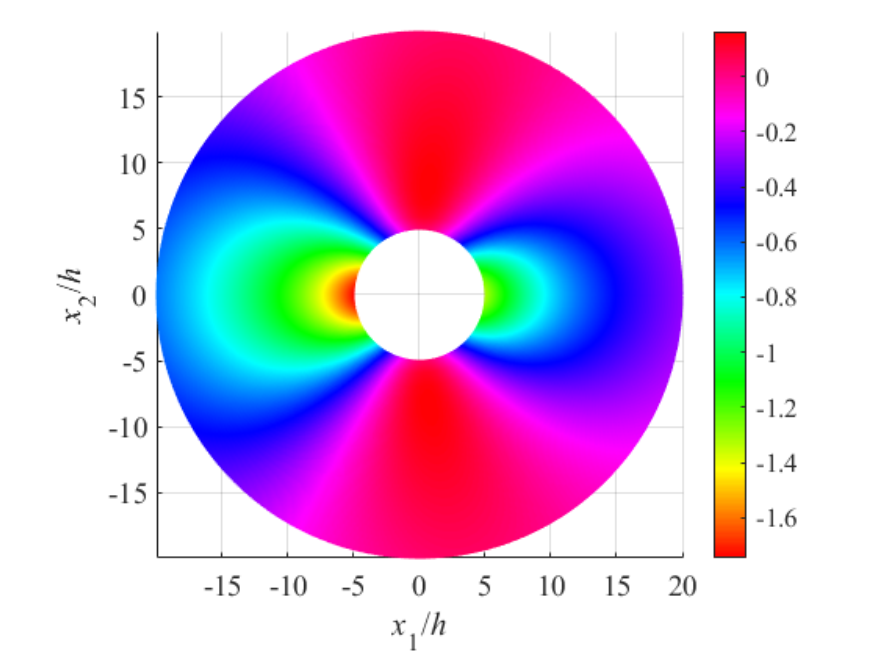}}
 \subfigure[.]{\includegraphics[scale=0.36]{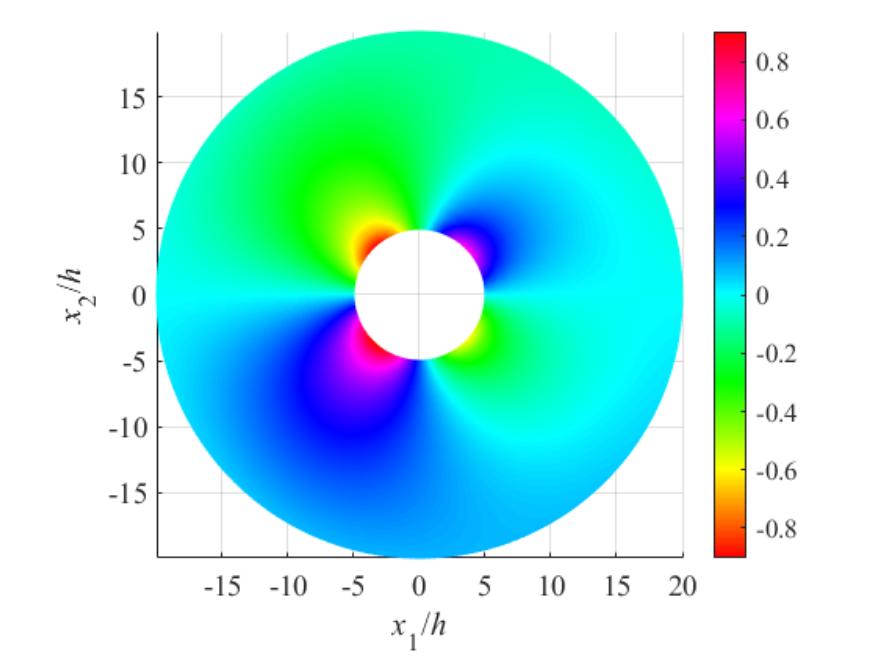}}
\subfigure[.]{\includegraphics[scale=0.36]{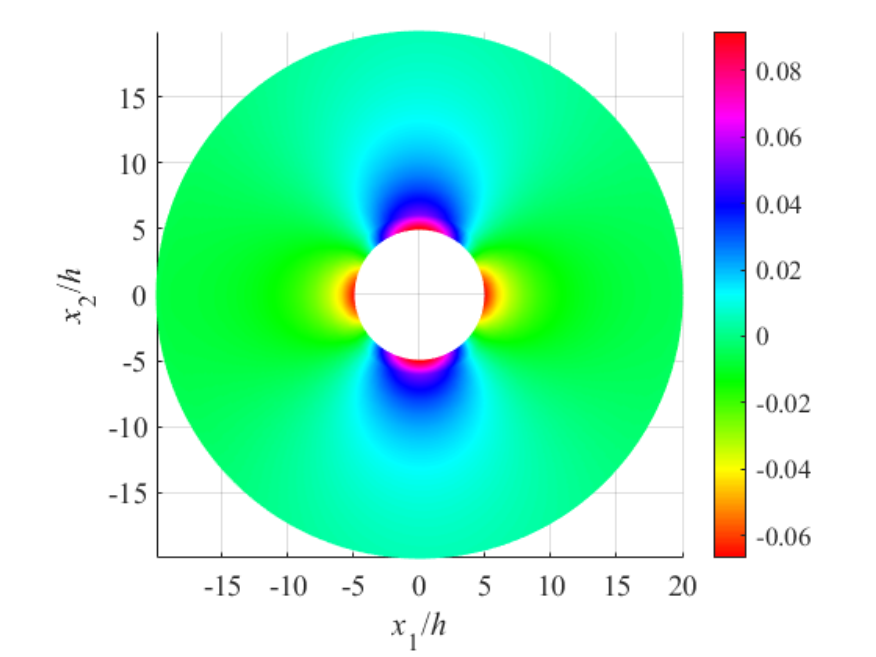}}
\subfigure[.]{\includegraphics[scale=0.36]{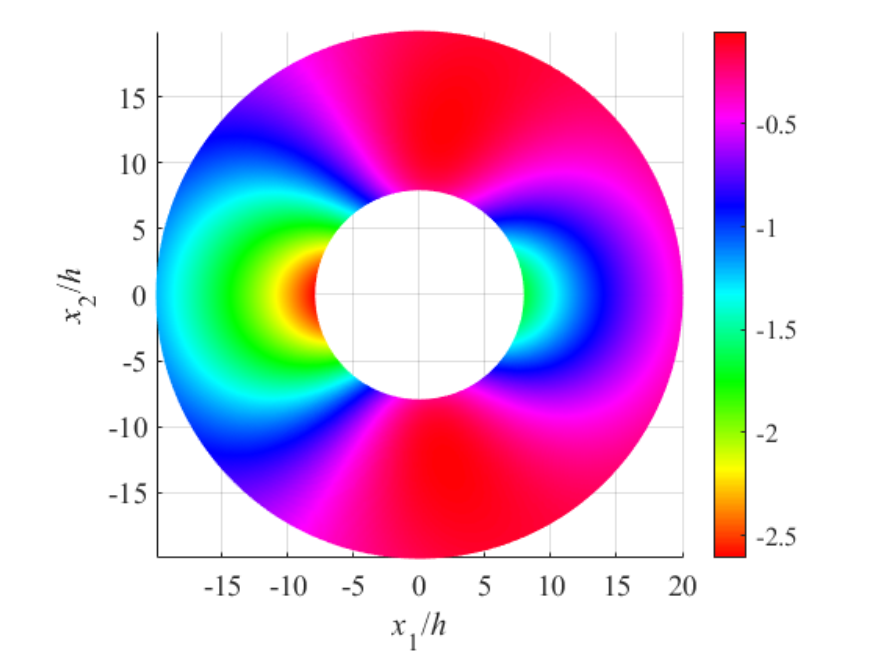}}
\subfigure[.]{\includegraphics[scale=0.36]{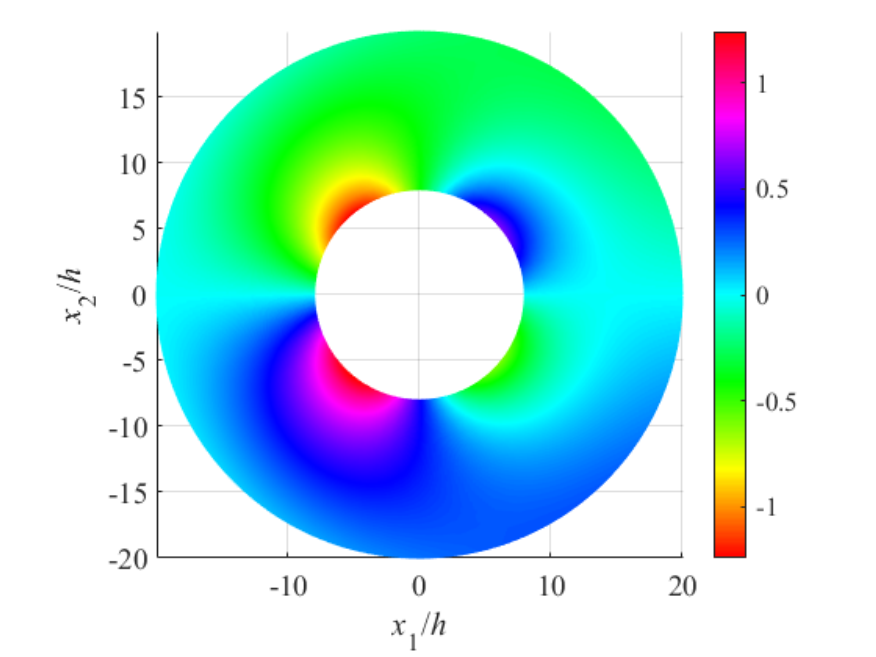}}
 \subfigure[.]{\includegraphics[scale=0.36]{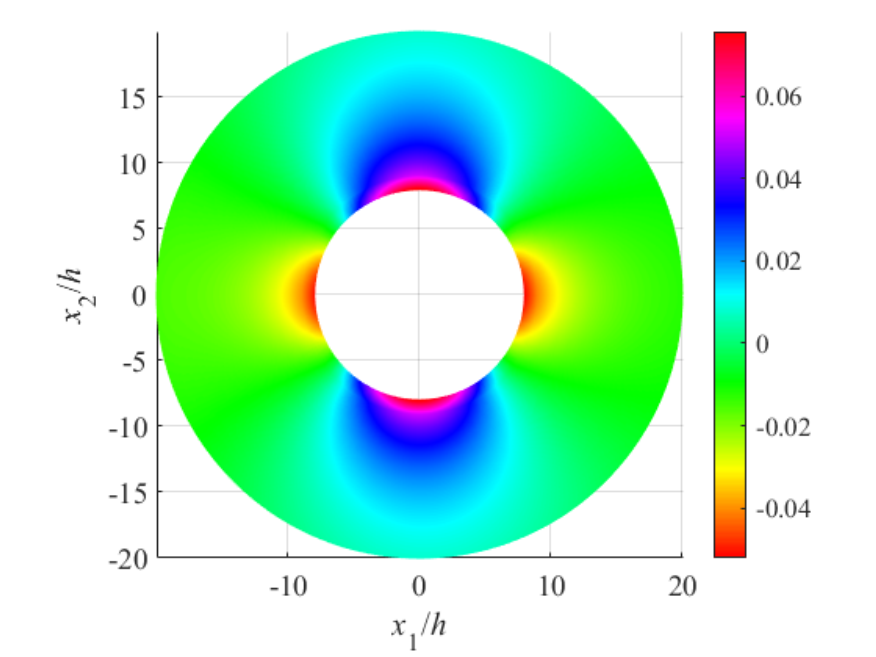}}
\caption{Error analysis of energy balance}
\label{fig12}
\end{figure}

\begin{figure}[h]
\centering
\subfigure[.]{\includegraphics[scale=0.36]{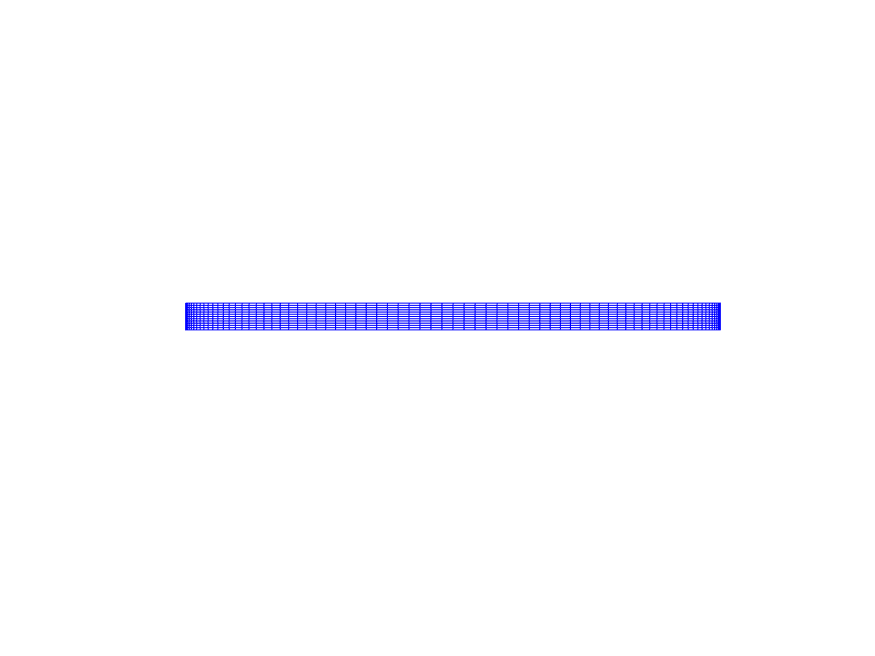}}
\subfigure[.]{\includegraphics[scale=0.36]{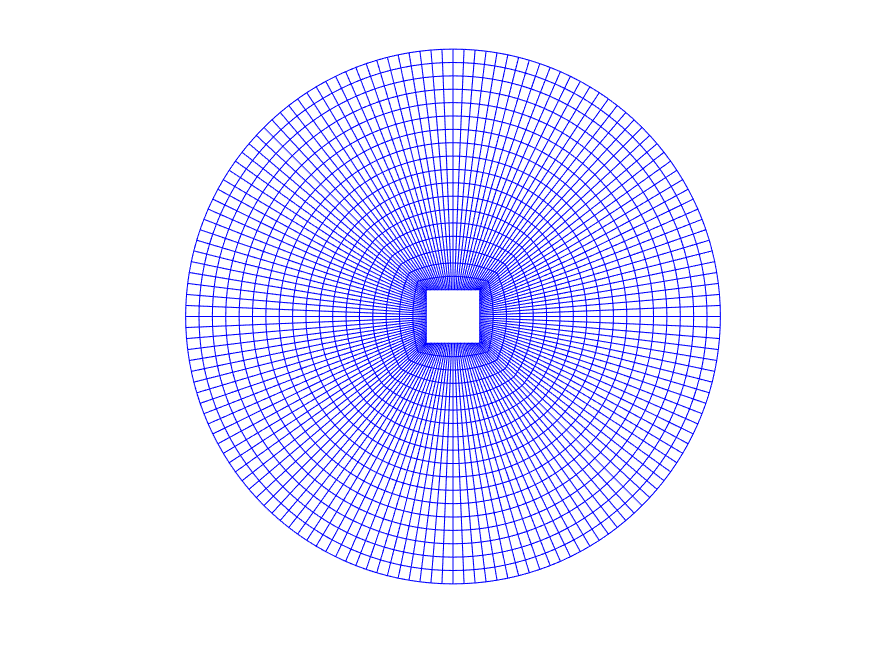}}
\subfigure[.]{\includegraphics[scale=0.36]{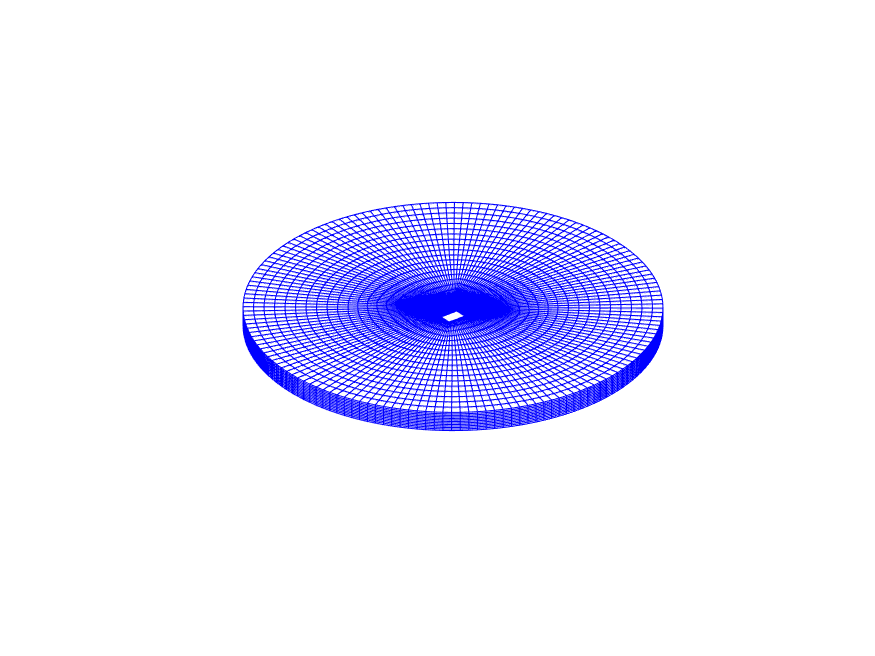}}
\subfigure[.]{\includegraphics[scale=0.36]{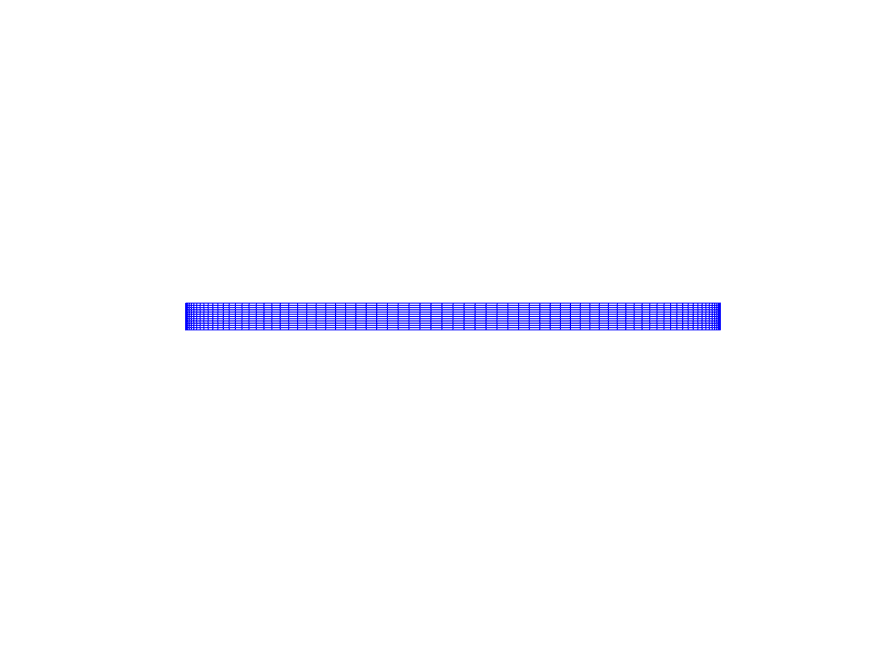}}
\subfigure[.]{\includegraphics[scale=0.36]{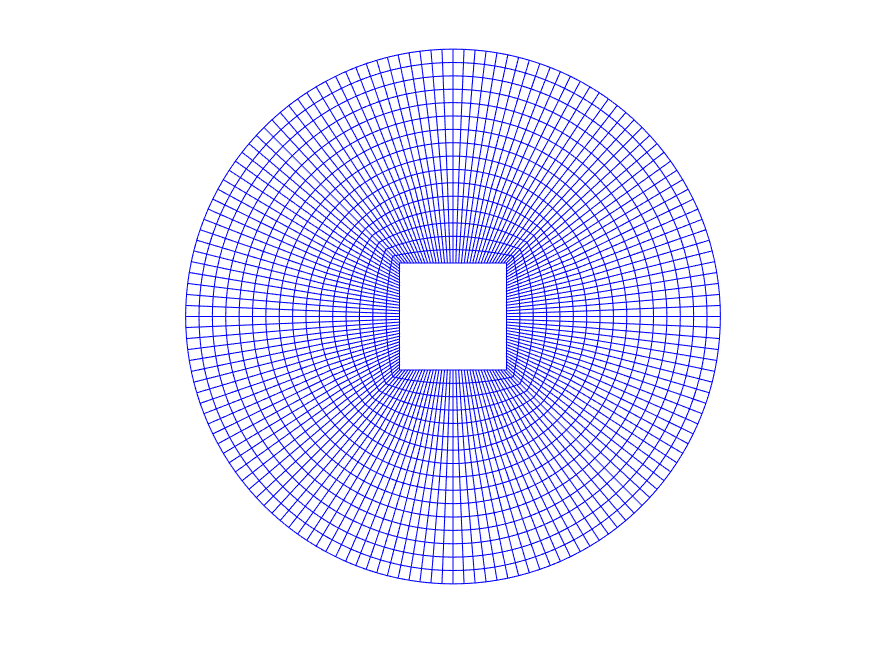}}
\subfigure[.]{\includegraphics[scale=0.36]{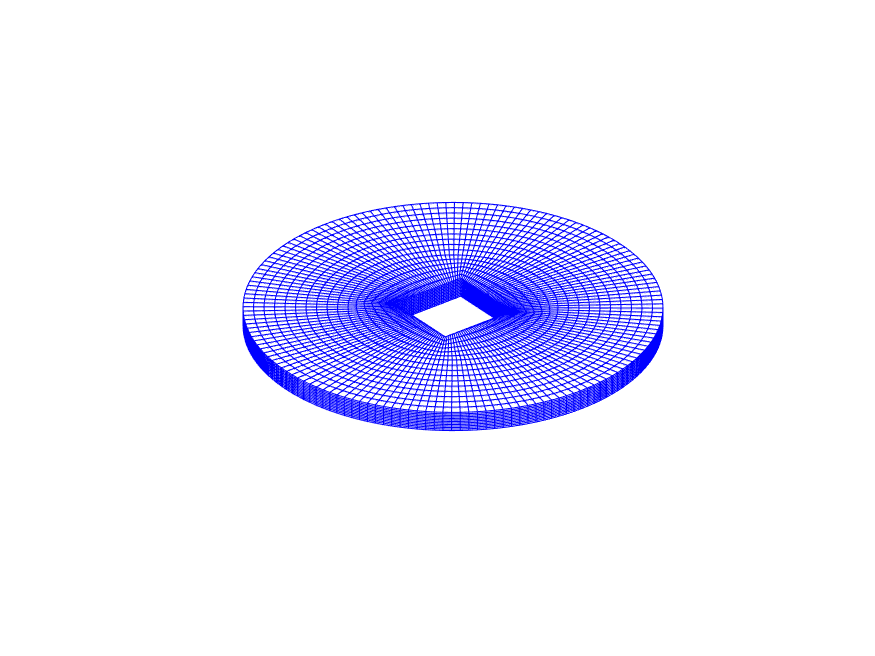}}
\subfigure[.]{\includegraphics[scale=0.36]{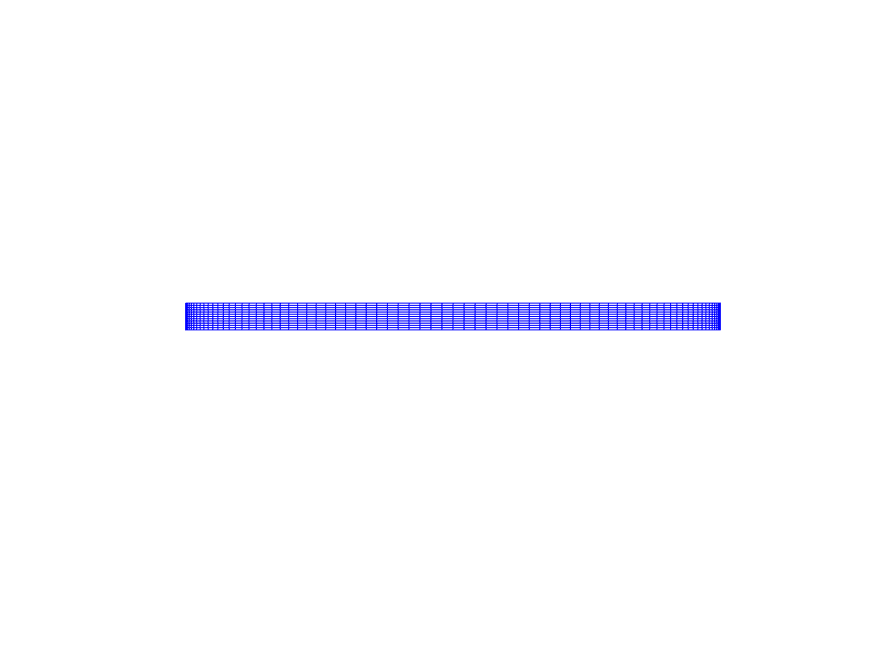}}
\subfigure[.]{\includegraphics[scale=0.36]{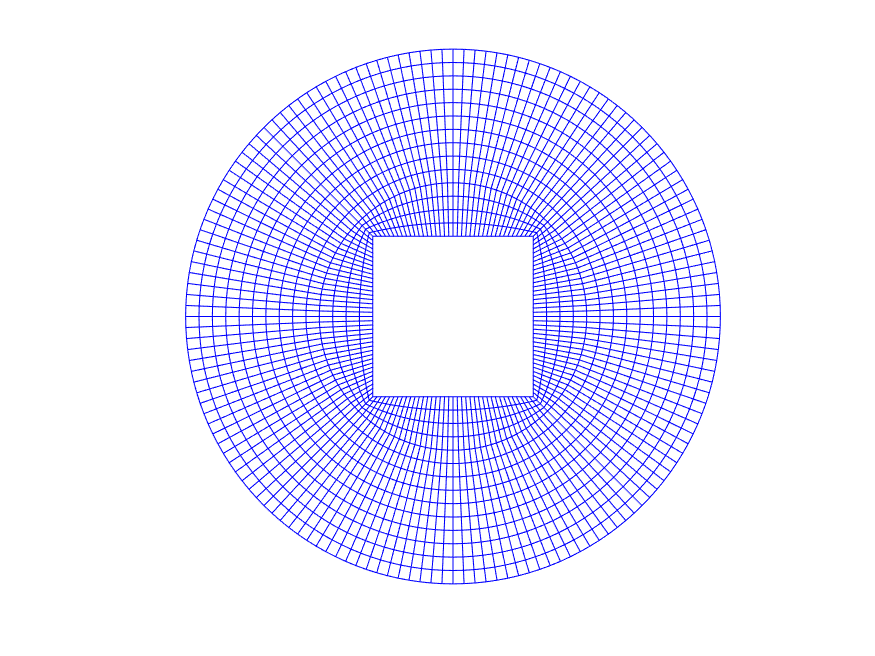}}
\subfigure[.]{\includegraphics[scale=0.36]{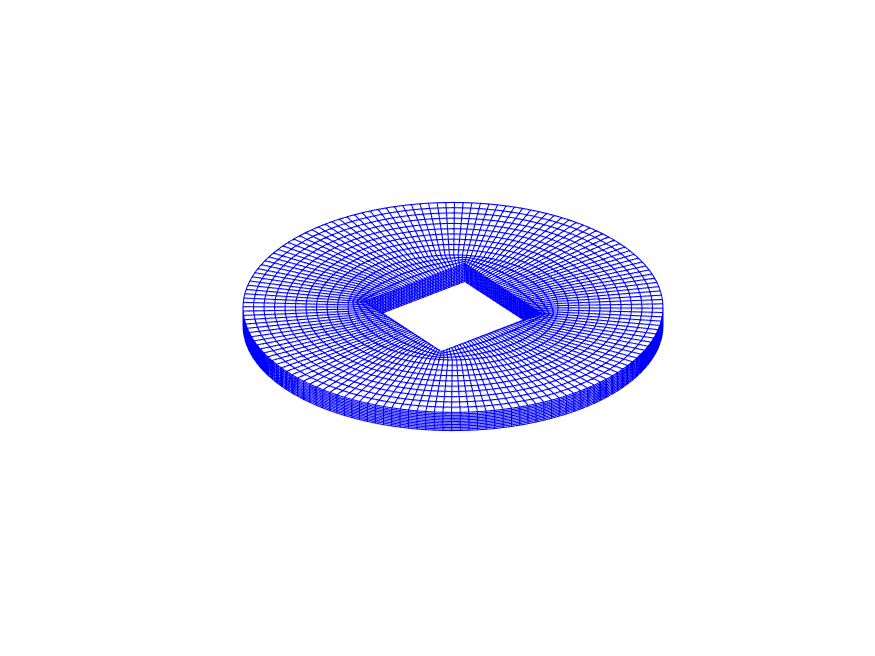}}
\caption{Error analysis of energy balance}
\label{fig13}
\end{figure}

\begin{figure}[h]
\centering
 \subfigure[.]{\includegraphics[scale=0.56]{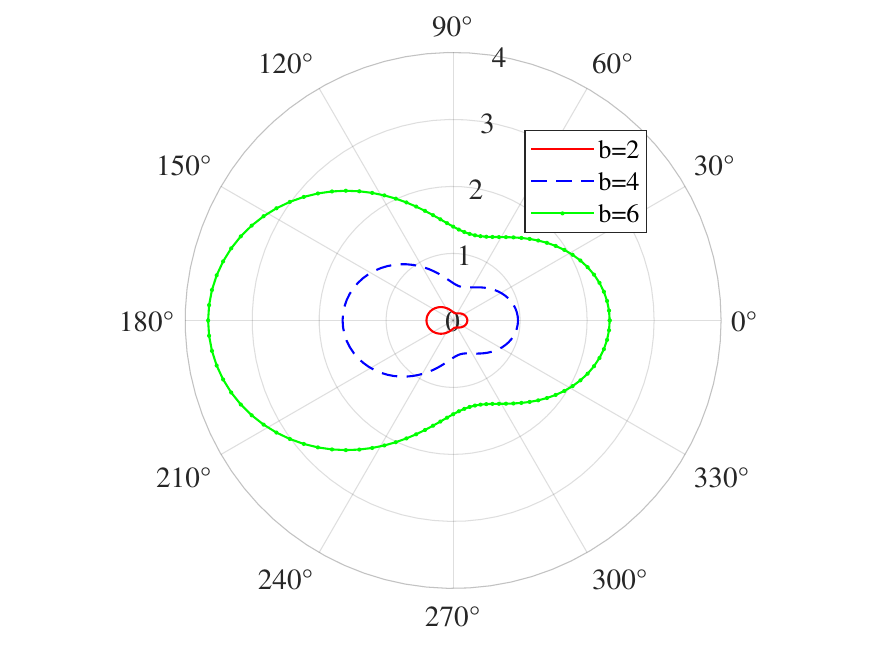}}
\subfigure[.]{\includegraphics[scale=0.56]{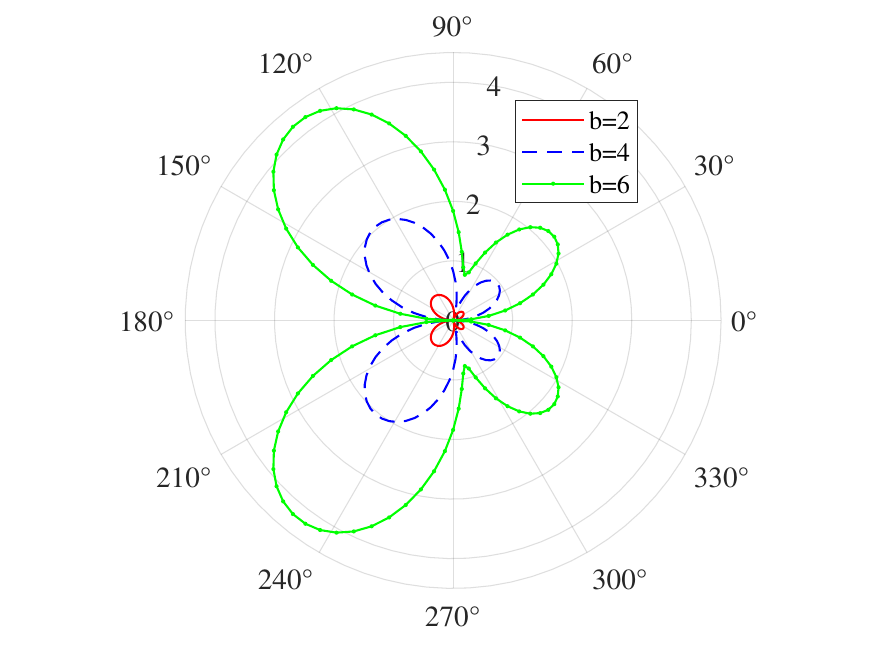}}
\caption{Error analysis of energy balance}
\label{fig14}
\end{figure}

\begin{figure}[h]
\centering
 \subfigure[.]{\includegraphics[scale=0.36]{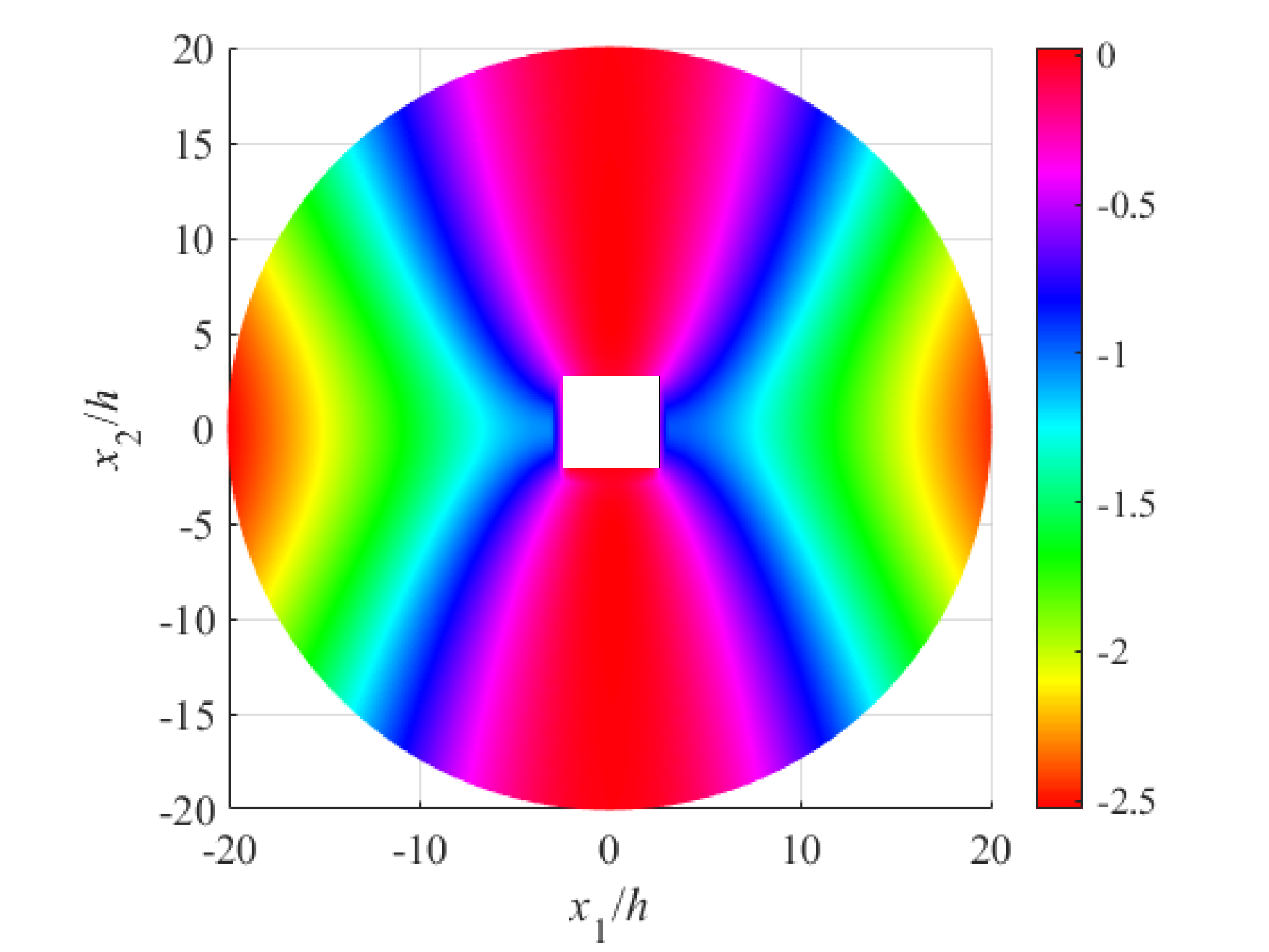}}
\subfigure[.]{\includegraphics[scale=0.36]{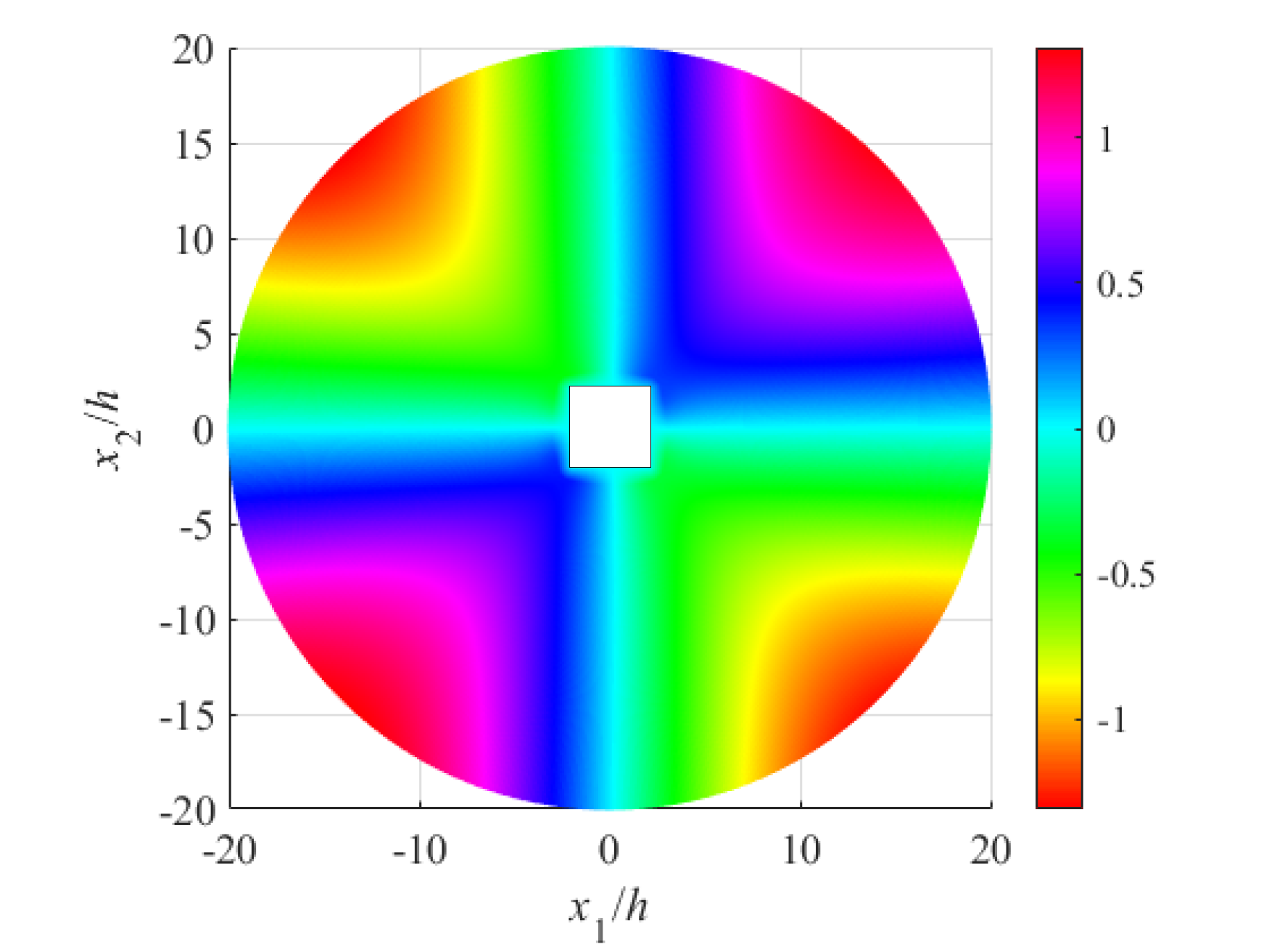}}
\subfigure[.]{\includegraphics[scale=0.36]{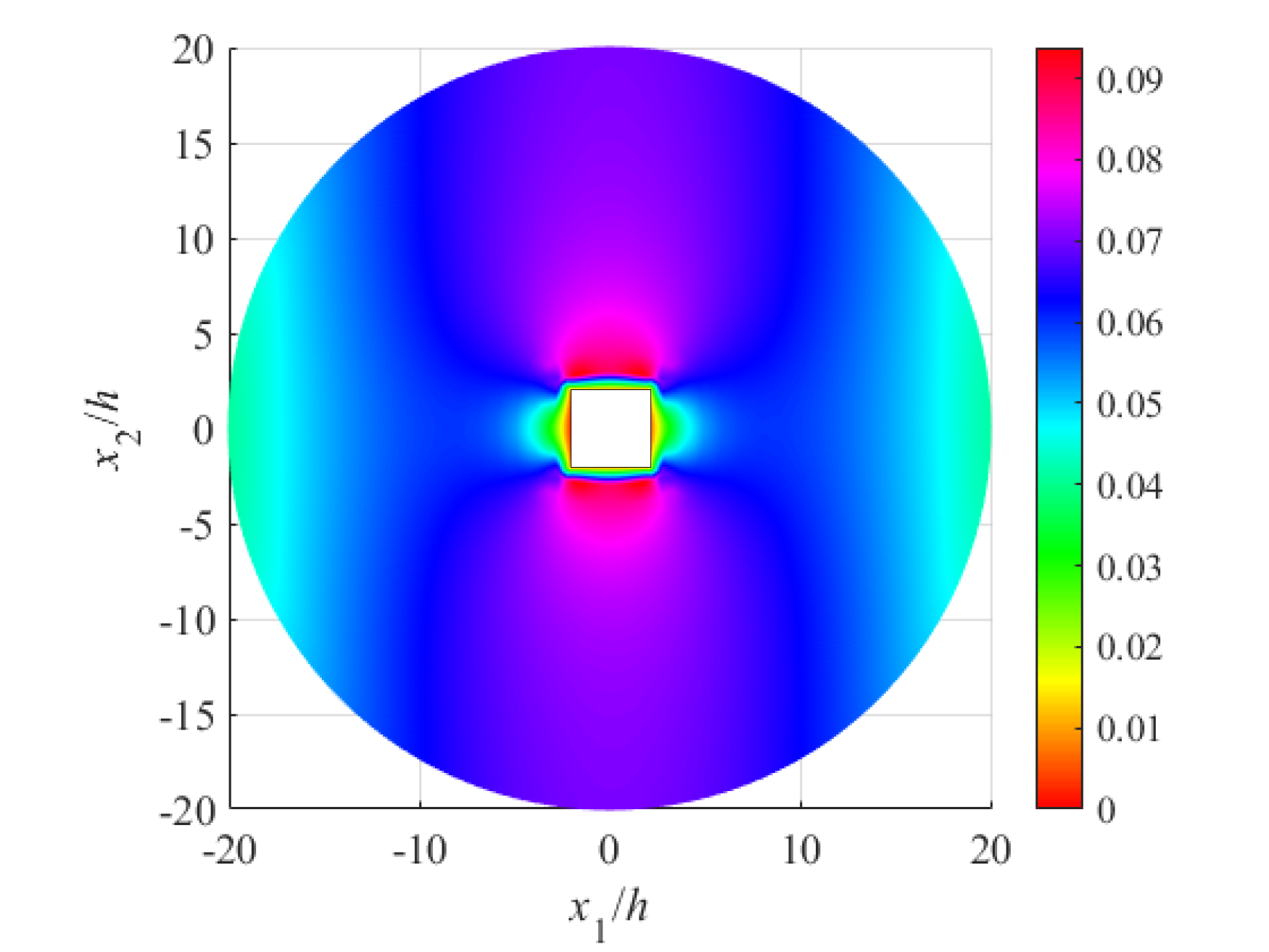}}
\subfigure[.]{\includegraphics[scale=0.36]{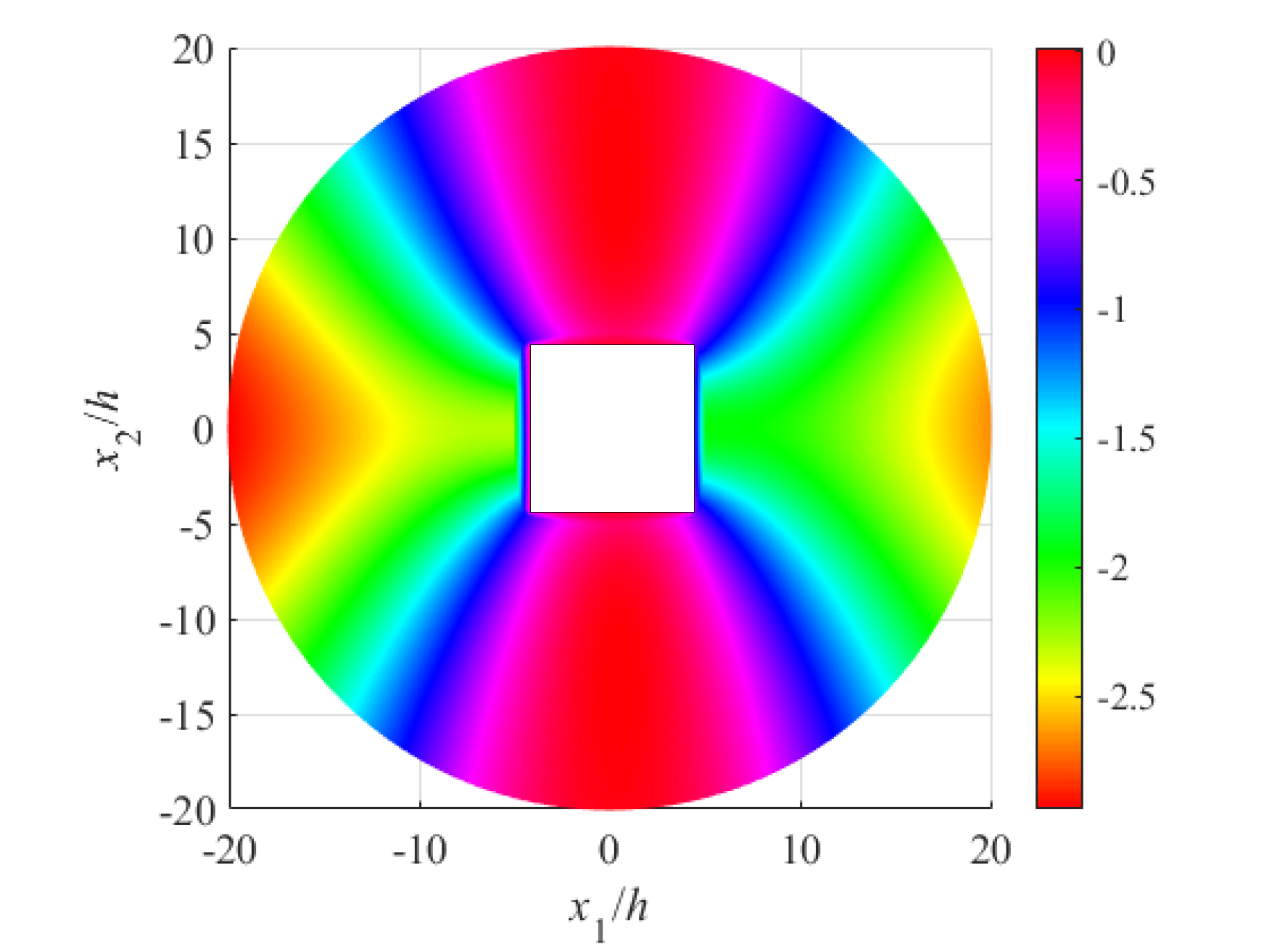}}
 \subfigure[.]{\includegraphics[scale=0.36]{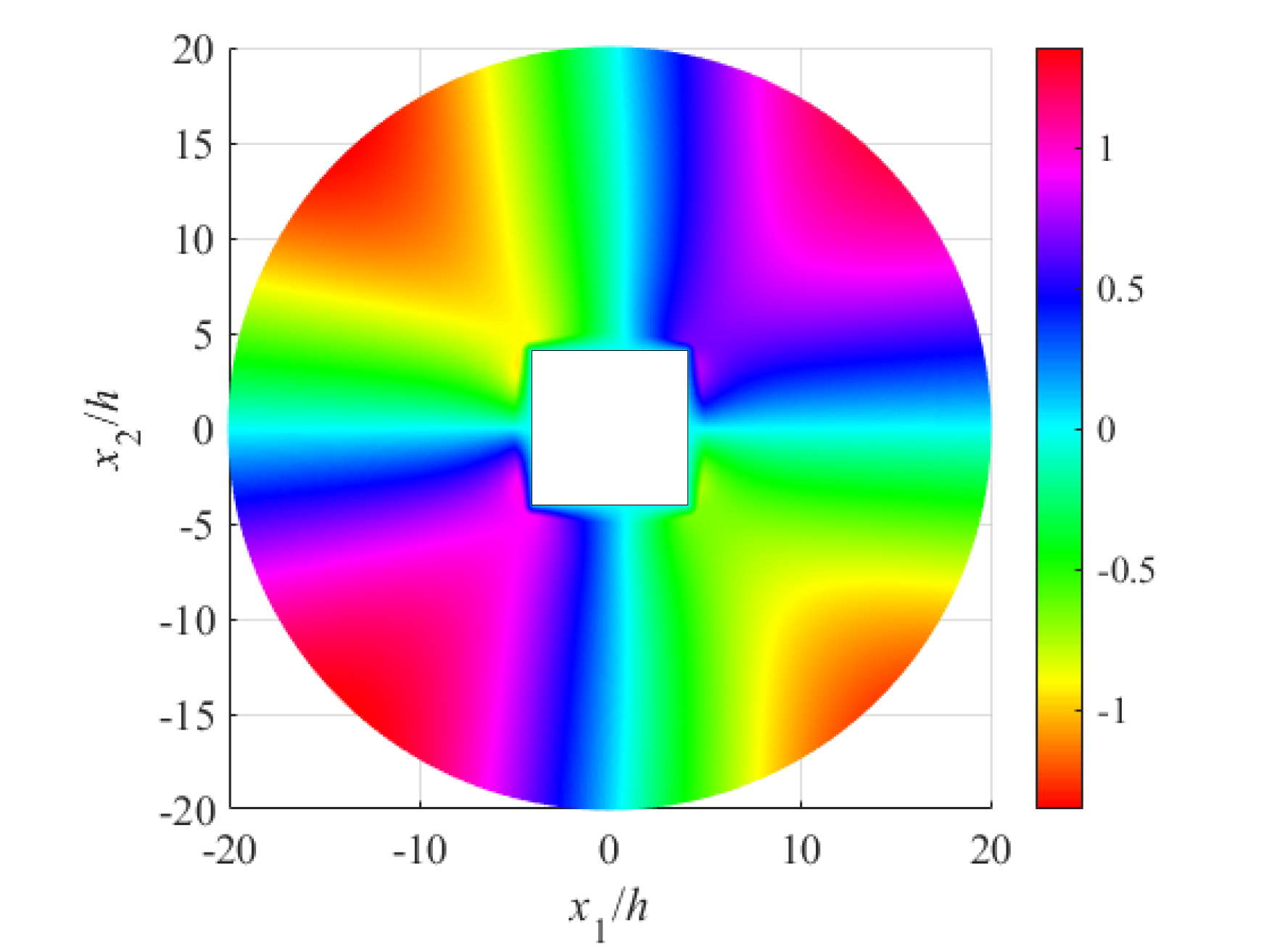}}
\subfigure[.]{\includegraphics[scale=0.36]{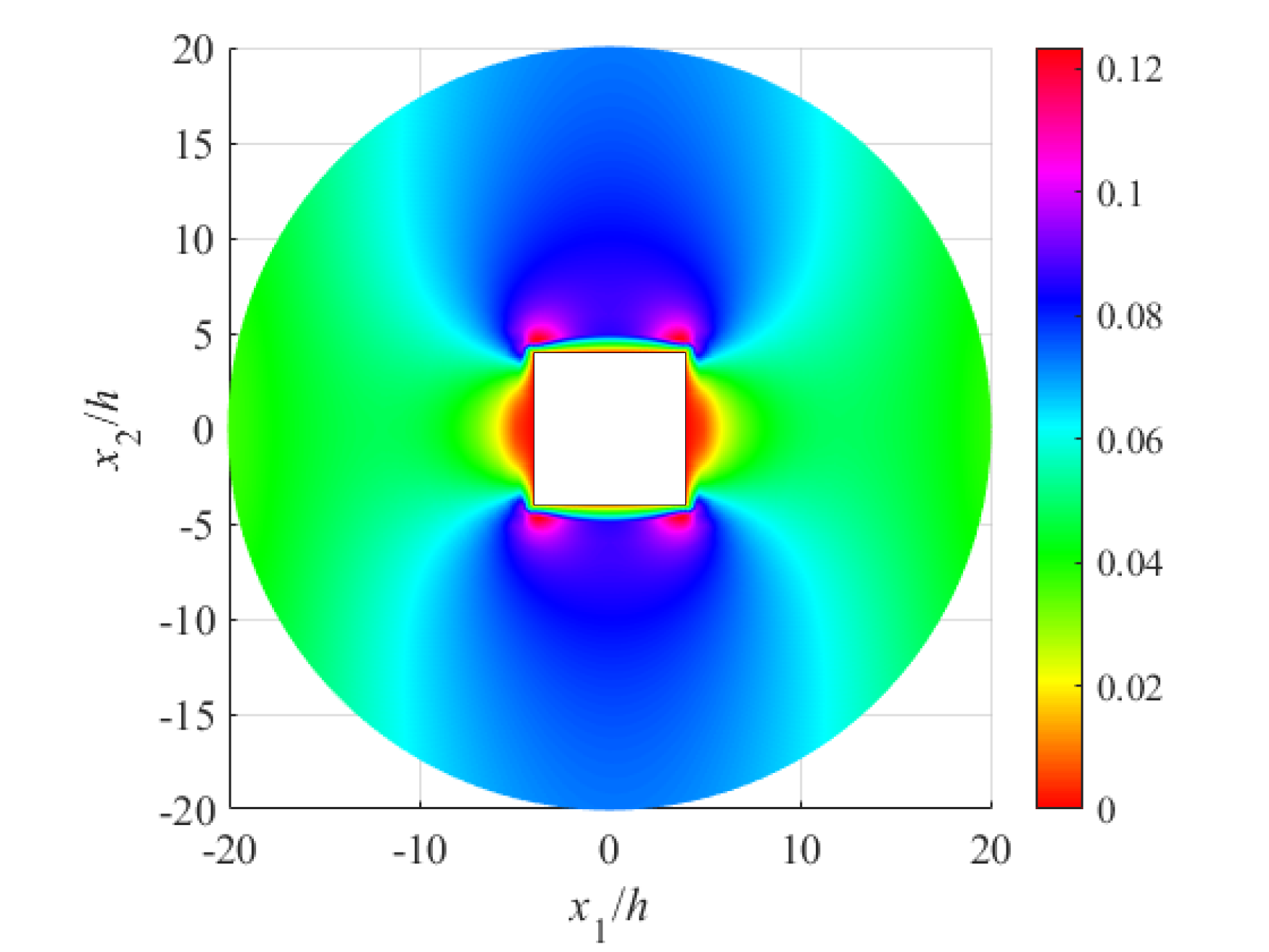}}
\subfigure[.]{\includegraphics[scale=0.36]{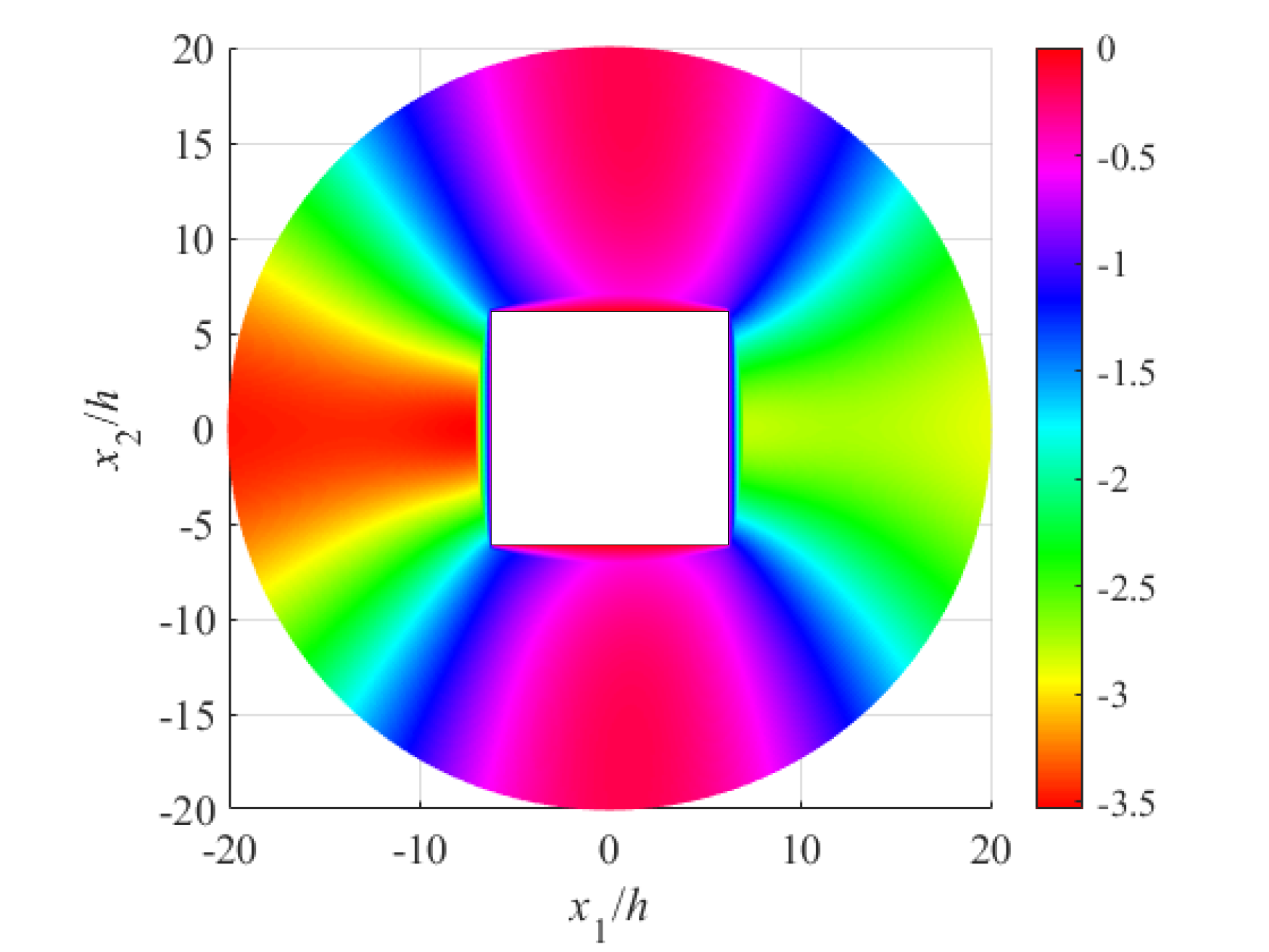}}
\subfigure[.]{\includegraphics[scale=0.36]{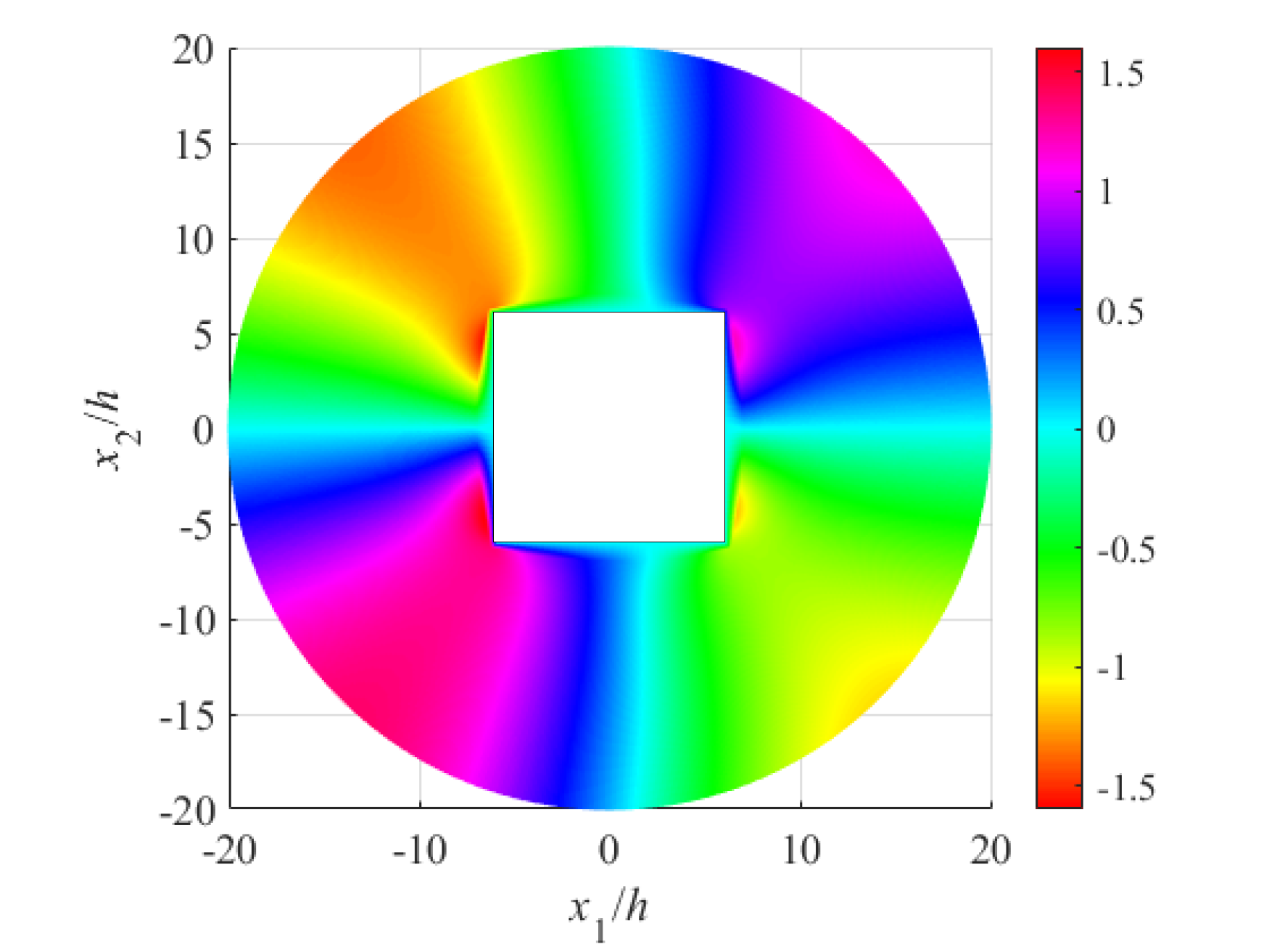}}
 \subfigure[.]{\includegraphics[scale=0.36]{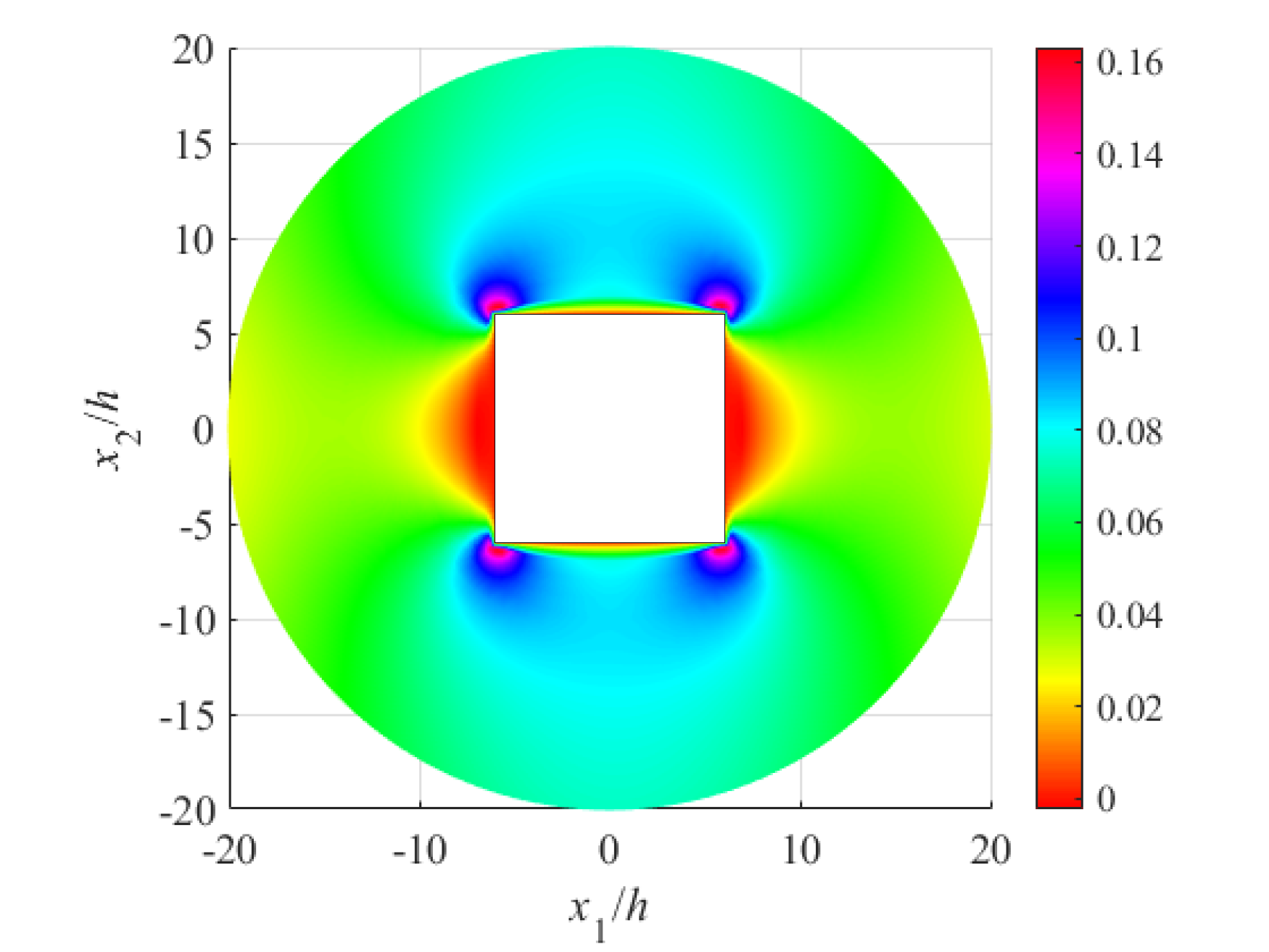}}
\caption{Error analysis of energy balance}
\label{fig15}
\end{figure}

\section{Conclusion}

This paper proposed a three-dimensional DtN FEM to analyze the scattering of Lamb and SH guided waves due to different shapes in an infinite 3D plate. This method is simple and elegant, which has advantages on dimension reduction and needs no absorption medium or perfectly matched layer to suppress the reflected waves compared to traditional FEM. For the future development, the scattering data from forward analysis by proposed 3D DtN FEM will subsequently be used for the inverse analysis of reconstructing both the location and shape of 3D defects.

\section{Conflicts of Interest}

The authors declare that they have no known competing financial interests or personal relationships that could have appeared to influence the work reported in this paper.

\section*{Acknowledgments}

This work was supported by the National Natural Science Foundation of China (Grant No. 12402098), the China Scholarships Council (Grant No. 202106830041) and the Fundamental Research Funds for the Central Universities (Grant No. 7100604733). 


\bibliography{mybibfile}

\end{document}